\documentclass[review,hidelinks,onefignum,onetabnum]{siamart220329}

\usepackage{lipsum}
\usepackage{amsfonts}
\usepackage{graphicx}
\usepackage{epstopdf}
\usepackage{algorithmic}
\ifpdf
  \DeclareGraphicsExtensions{.eps,.pdf,.png,.jpg}
\else
  \DeclareGraphicsExtensions{.eps}
\fi

\newsiamremark{remark}{Remark}
\newsiamremark{hypothesis}{Hypothesis}
\crefname{hypothesis}{Hypothesis}{Hypotheses}
\newsiamthm{claim}{Claim}

\headers{ReLU neural network approximation}{Z. Cai, J. Choi, and M. Liu}

\title{ReLU neural network approximation\\  to piecewise constant functions
\thanks{Submitted to the editors DATE.
\funding{This work was supported in part by the National Science Foundation under grant DMS-2110571.}}}

\author{Zhiqiang Cai\thanks{Department of Mathematics, Purdue University, 150 N. University Street, West Lafayette, IN 47907-2067 
  (\email{caiz@purdue.edu}, \email{choi508@purdue.edu}).}
\and Junpyo Choi\footnotemark[2]
\and Min Liu\thanks{School of Mechanical Engineering, Purdue University, 585 Purdue Mall,
West Lafayette, IN 47907-2088(\email{liu66@purdue.edu}). }}

\usepackage{amsopn}

\ifpdf
\hypersetup{
  pdftitle={ReLU Neural Network Approximation to Piecewise Constant Functions},
  pdfauthor={Z. Cai, J. Choi, and M. Liu}
}
\fi
\nolinenumbers
\usepackage[utf8]{inputenc}
\usepackage{bm}
\usepackage{mathtools}
\usepackage{indentfirst}
\usepackage{subfigure}
\usepackage{tcolorbox}
\usepackage{color}
\usepackage[export]{adjustbox}
\usepackage{scalerel}
\usepackage{tikz}
\usepackage{pgfplots}
\pgfplotsset{compat=1.18}
\usetikzlibrary{decorations.pathreplacing}

\usepackage{algorithmic}
\Crefname{ALC@unique}{Line}{Lines}

\newcommand{\R}{\mathbb{R}}

\usepackage{graphicx}
\usepackage{diagbox}
\usepackage{pifont}
\newcommand{\vertiii}[1]{{\left\vert\kern-0.25ex\left\vert\kern-0.25ex\left\vert #1 
    \right\vert\kern-0.25ex\right\vert\kern-0.25ex\right\vert}}

\usepackage{enumitem}
\setlist[itemize]{left=16pt} 

\def\bx{{\bf x}}

\begin{document}

\maketitle

\begin{abstract}
This paper studies the approximation property of ReLU neural networks (NNs) to piecewise constant functions with unknown interfaces in bounded regions in $\mathbb{R}^d$. Under the assumption that the discontinuity interface $\Gamma$ may be approximated by a connected series of hyperplanes with a prescribed accuracy $\varepsilon >0$, we show that a three-layer ReLU NN is sufficient to accurately approximate any piecewise constant function and establish its error bound. Moreover, if the discontinuity interface is convex, an analytical formula of the ReLU NN approximation with exact weights and biases is provided.
\end{abstract}



\section{Introduction}

For simplicity, consider the $d$-dimensional unit cube $\Omega=(0,1)^d$ with $d\ge 2$. Let $\left\{\Omega_1,\Omega_2\right\}$ be a partition of the domain $\Omega$; that is, $\Omega_1$ and $\Omega_2$ are open and connected subdomains of $\Omega$ such that 
\[
\Omega_1\cap  \Omega_2=\emptyset 
\quad\mbox{and}\quad
\bar{\Omega}=\bar{\Omega}_1\cup  \bar{\Omega}_2.
\]
Let $\chi(\bx)$ be a piece-wise constant function defined on $\Omega$ given by
\begin{equation}\label{step function}
 \chi(\bx)=\left\{\begin{array}{rl}
 0, & \bx \in \Omega_1,\\[2mm]
 1, & \bx \in \Omega_2.
 \end{array}
 \right.   
\end{equation}
Denote by $\Gamma=\partial \Omega_1\cap \partial\Omega_2$ the discontinuity interface of $\chi(\bx)$, where $\partial \Omega_i$ is the boundary of the subdomain $\Omega_i$. In this paper, we assume that the interface $\Gamma$ is in $C^0$ and that its $(d-1)$-dimensional measure $|\Gamma|$ is finite.

Functions of the form in \eqref{step function} are encountered in many applications such as classification tasks in data science and linear and nonlinear hyperbolic conservation laws with discontinuous solutions (see, e.g., \cite{anthony1999neural,IF22,Cai2023nonlinear, cai2023least}). Generally, a piecewise constant function has the form 
\begin{equation}\label{disconti class}
\chi(\bx)= \sum_{i=1}^m \alpha_i\chi_i (\bx),  
\end{equation}
where $\alpha_i$ is a real number, $\chi_i(\bx)={\bm 1}_{\Omega_i}(\bx)$ is the  indicator function of a subdomain $\Omega_i\subset \Omega$, and 
$\left\{\Omega_i\right\}_{i=1}^m$ forms a partition of the domain $\Omega$. The partition means that $\left\{\Omega_1,\ldots,\Omega_m\right\}$ are open, connected, and disjoint subdomains of $\Omega$ and that $\bar{\Omega}=\cup_{i=1}^m\bar{\Omega}_i$. Once we know how to approximate $\chi(\bx)$ in \eqref{step function} by neural networks (NNs), then approximating \eqref{disconti class} is a matter of concatenation or parallelization of the NNs (see, e.g., \cite{DeVore2021}).

A critical component of using NNs as a model is the use of a properly designed architecture (e.g., the number of layers), and carelessly chosen architectures could lead to poor performance regardless of the size of the network (see, e.g., \cite{cai2023least,Cai2023linearb}).
To efficiently approximate piecewise constant functions with unknown interface location, several practical guidelines on the architecture of NNs have been provided recently (see, e.g., \cite{imaizumi2019deep,petersen2018optimal,caragea2023neural,cai2023least}). The first notable work was done by Petersen and Voigtlaender in their 2018 paper \cite{petersen2018optimal}. 
For any prescribed accuracy $\varepsilon>0$, if the discontinuity interface $\Gamma$ is in $C^\beta$ with $\beta>0$, they showed that there exists a NN function $\mathcal{N}(\bx)$, generated by a ReLU NN with at most $(3+\lceil \log_2\beta\rceil)(11+2\beta/d)$ layers and at most $c\varepsilon^{-p(d-1)/\beta}$ nonzero weights for some constant $c>0$, such that 
\begin{equation}\label{petersen2018}
    \|\chi-\mathcal{N}\|_{L^p(\Omega)}\leq \varepsilon.
\end{equation}
In the case that $\Gamma$ can locally be parametrized by functions of Barron-type, it was proved in \cite{caragea2023neural} that for every $N\in\mathbb{N}$, there exists a NN function $\mathcal{N}(\bx)$, generated by a four-layer ReLU NN with a total of $\mathcal{O}(d+N)$ neurons, such that 
\[
\|\chi-\mathcal{N}\|_{L^p(\Omega)}\leq C\,d^{\frac{3}{2p}}N^{-\frac{\alpha}{2p}},
\]
where $C$ and $\alpha$ are positive constants independent of $N$. Here, the magnitude of the weights and biases can be chosen to be $\mathcal{O}(d+N^{1/2})$.

Recently, we studied this problem in \cite{cai2023least} through an explicit construction based on the {\it two-layer} ReLU approximation $p(\bx)$ in Lemma~3.2 of \cite{Cai2021linear}. Under the assumption that the interface $\Gamma$ may be approximated such that there exists a region of $\varepsilon$ width containing the interface, we were able to construct a continuous piecewise linear (CPWL) function with a sharp transition layer of $\varepsilon$ width whose approximation to the piecewise constant function $\chi(\bx)$ has the approximation accuracy $\varepsilon$. Combining with the main results in \cite{arora2016understanding}, this indicates that a ReLU NN with at most $\lceil \log_2(d+1)\rceil+1$ layers is sufficient to achieve the prescribed accuracy $\varepsilon$. However, \cite{cai2023least} does not provide an estimate of the minimum number of neurons in each layer.

The purpose of this paper is to address the following two questions: 
\begin{itemize}
    \item[(1)] What is the minimum number of hidden-layers of a ReLU NN in order to approximate a piecewise constant function with the prescribed accuracy?
    \item[(2)] How many neurons per each hidden-layer are needed?
\end{itemize}
Under the assumption that the interface $\Gamma$ may be approximated by a connected series of hyperplanes with a prescribed accuracy $\varepsilon >0$ (see Figure \ref{general22}), we show that a {\it three-layer} (two-hidden-layer) ReLU NN is sufficient and necessary to accurately approximate the piecewise constant function $\chi(\mathbf{x})$, in any dimensions, with an error bound of $\mathcal{O}(\varepsilon^{1/p})$ in the $L^p(\Omega)$ norm (see Theorem \ref{general theorem}). Again, this is done through an explicit construction based on a novel {\it three-layer} ReLU NN approximation (see, e.g, $\mathcal{N}(\bx)$ in \eqref{N-1} when the interface is a hyperplane). Moreover, the number of neurons at the first hidden-layer and their locations depend on the hyperplanes used for approximating the interface and the number of neurons of the second hidden-layer depends on convexity of the interface.

For classification problems or partial differential equations with a discontinuous solution, our approximation results would provide a guideline on the choice of ReLU NN architectures and on initialization for any training algorithm. It is well-known that initialization is critical for success of any optimization/iterative/training scheme when the resulting discrete problem is a non-convex optimization.

The remainder of the paper is organized as follows. Three-layer ReLU NN functions with relevant concepts and terminology are described in Section \ref{s:NN}. Then in Section \ref{main section}, we describe how to approximate the interface $\Gamma$ with necessary assumptions, and state the main result of the approximation theory by three-layer ReLU NN functions. The proof of a lemma for the theorem is provided in Section \ref{general proof}. Finally, multiple examples with $d\ge 2$ are given in Section \ref{s:example} to confirm our theoretical findings.

\section{Three-layer ReLU neural network functions}\label{s:NN}
In this paper, we will restrict our attention to three-layer (two-hidden-layer) neural network functions that are scalar-valued. A function $\mathcal{N}:\mathbb{R}^d\to\mathbb{R}$ is a three-layer neural network (NN) function if the function $\mathcal{N}$ has a representation as a composition of 3 functions $\mathbf{x}^{(l)}:\mathbb{R}^{n_{l-1}}\to\mathbb{R}^{n_{l}}$ ($n_0=d$, $n_3=1$) for $l=1,2,3$:
\begin{equation}\label{three nn}
    \mathcal{N}=\mathbf{x}^{(3)}\circ \mathbf{x}^{(2)}\circ \mathbf{x}^{(1)},
\end{equation}
where $\mathbf{x}^{(3)}$ is affine linear, and $\mathbf{x}^{(2)}$ and $\mathbf{x}^{(1)}$ are affine linear with a function $\sigma:\mathbb{R}\to\mathbb{R}$, called an activation function, applied to each component of the functions. Such a function is called a $d$--$n_1$--$n_2$--$1$ NN function.

As the activation function, we use the rectified linear unit (ReLU):
\begin{equation*}
\sigma(t)=\mathrm{ReLU}(t) \coloneqq \max\{0,t\}
=
\begin{cases}
0, & \text{if }  t\leq 0,\\
t, & \text{if } t >0,
\end{cases}
\end{equation*}
and refer to such a three-layer NN function as a three-layer (two-hidden-layer) ReLU NN function.
Therefore, the collection of all three-layer ReLU NN functions from $\mathbb{R}^d$ to $\mathbb{R}$ is the collection of all functions $\mathcal{N}:\mathbb{R}^d\to\mathbb{R}$ defined by
\[\mathcal{N}(\mathbf{x})=\bm{\omega}^{(3)}\sigma\left(\bm{\omega}^{(2)}\sigma\left(\bm{\omega}^{(1)}\mathbf{x}-\mathbf{b}^{(1)}\right)-\mathbf{b}^{(2)}\right)-\mathbf{b}^{(3)},\]
where for each $l=1,2,3$, $\bm{\omega}^{(l)}\in\mathbb{R}^{n_l\times n_{l-1}},\  \mathbf{b}^{(l)}\in\mathbb{R}^{n_l}$ for $n_l,n_{l-1}\in\mathbb{N}$. We may also assume that each row of the matrix $\bm{\omega}^{(1)}$ has unit length by adjusting the entries of $\bm{\omega}^{(2)}$ and $\mathbf{b}^{(1)}$ (see, e.g., \cite{DeVore2021}).

We will follow the same terminology in \cite{cai2023least}. Finally, in the numerical examples in this paper, as in \cite{cai2023least}, we will see the breaking hyperplanes of the first- and second-(hidden-) layers, which are defined as follows. For $l=1,2$, let
\[\bm{\omega}^{(l)}=(\mathbf{w}_1^{(l)},\ldots,\mathbf{w}_{n_l}^{(l)})^T\in\mathbb{R}^{n_l\times n_{l-1}},\quad\text{and}\quad \mathbf{b}^{(l)}=(b_1^{(l)},\ldots,b_{n_l}^{(l)})^T.\]
Then the first- (hidden-) layer breaking hyperplanes are for $i=1,\ldots,n_1$,
\[P_i^{(1)}=\left\{\mathbf{x}\in\mathbb{R}^d:\mathbf{w}_i^{(1)}\mathbf{x}-b_i^{(1)}=0\right\},\]
and the second- (hidden-) layer breaking (poly-) hyperplanes are for $i=1,\ldots,n_2$,
\[P_i^{(2)}=\left\{\mathbf{x}\in\mathbb{R}^d:\mathbf{w}_i^{(2)}\sigma\left(\bm{\omega}^{(1)}\mathbf{x}-\mathbf{b}^{(1)}\right)-b_i^{(2)}=0\right\}.\]

ReLU NN functions are continuous piecewise linear with respect to the partition of $\Omega\subset\mathbb{R}^d$ determined by the breaking hyperplanes. The constructions of approximations to piecewise constant functions in this paper will be better understood with the help of breaking hyperplanes.

\section{Main results}\label{main section}
\begin{figure}[htbp]
\centering
\subfigure[The interface $\Gamma$\label{general12}]{
\begin{minipage}[t]{0.4\linewidth}
\centering
\includegraphics[width=1.8in]{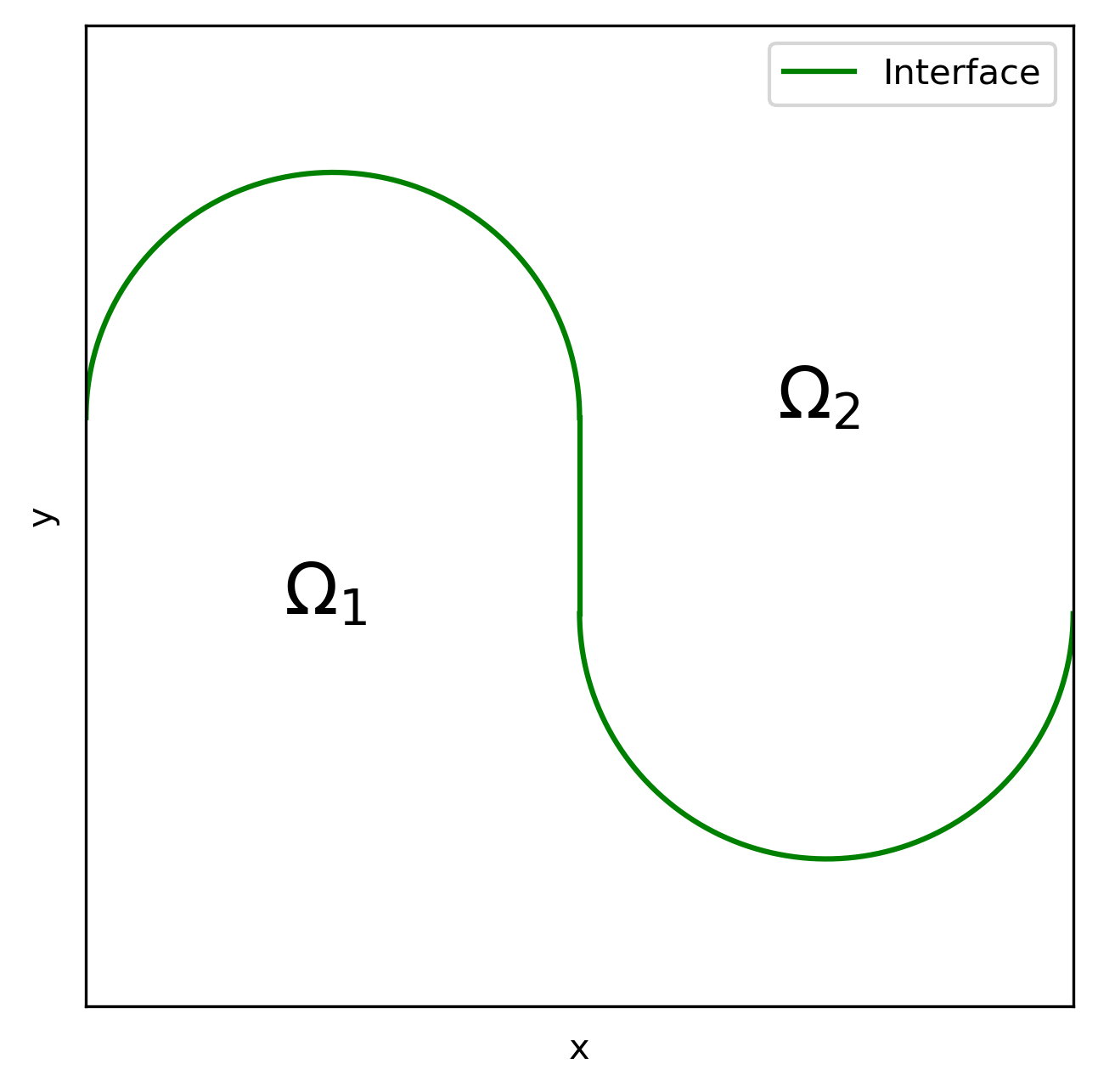}
\end{minipage}%
}%
\hspace{0.2in}
\subfigure[An approximation of the interface by connected series of hyperplanes\label{general22}]{
\begin{minipage}[t]{0.4\linewidth}
\centering
\includegraphics[width=1.8in]{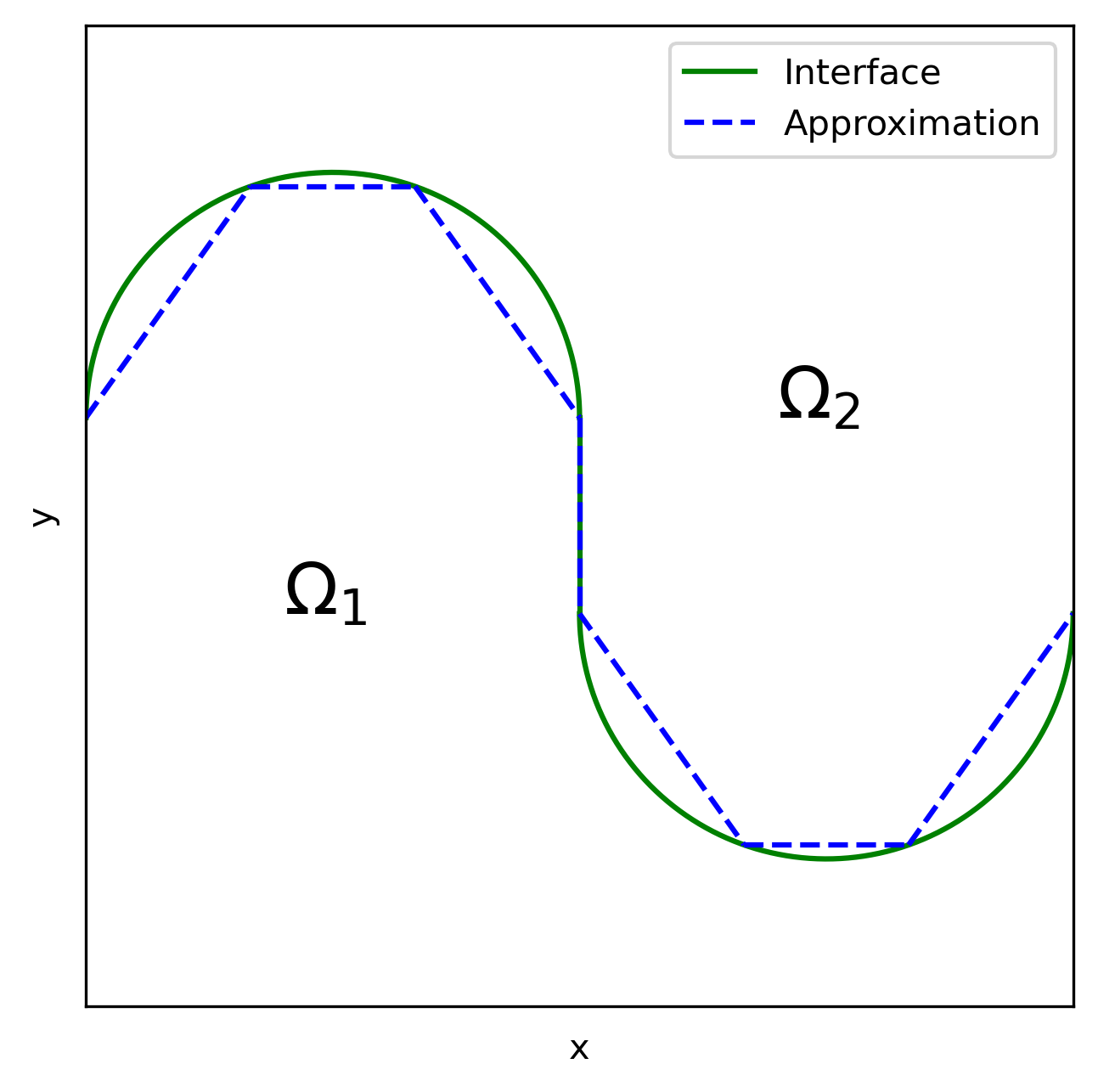}
\end{minipage}%
}%
\\
\subfigure[$\hat{\chi}$\label{general32}]{
\begin{minipage}[t]{0.4\linewidth}
\centering
\includegraphics[width=1.8in]{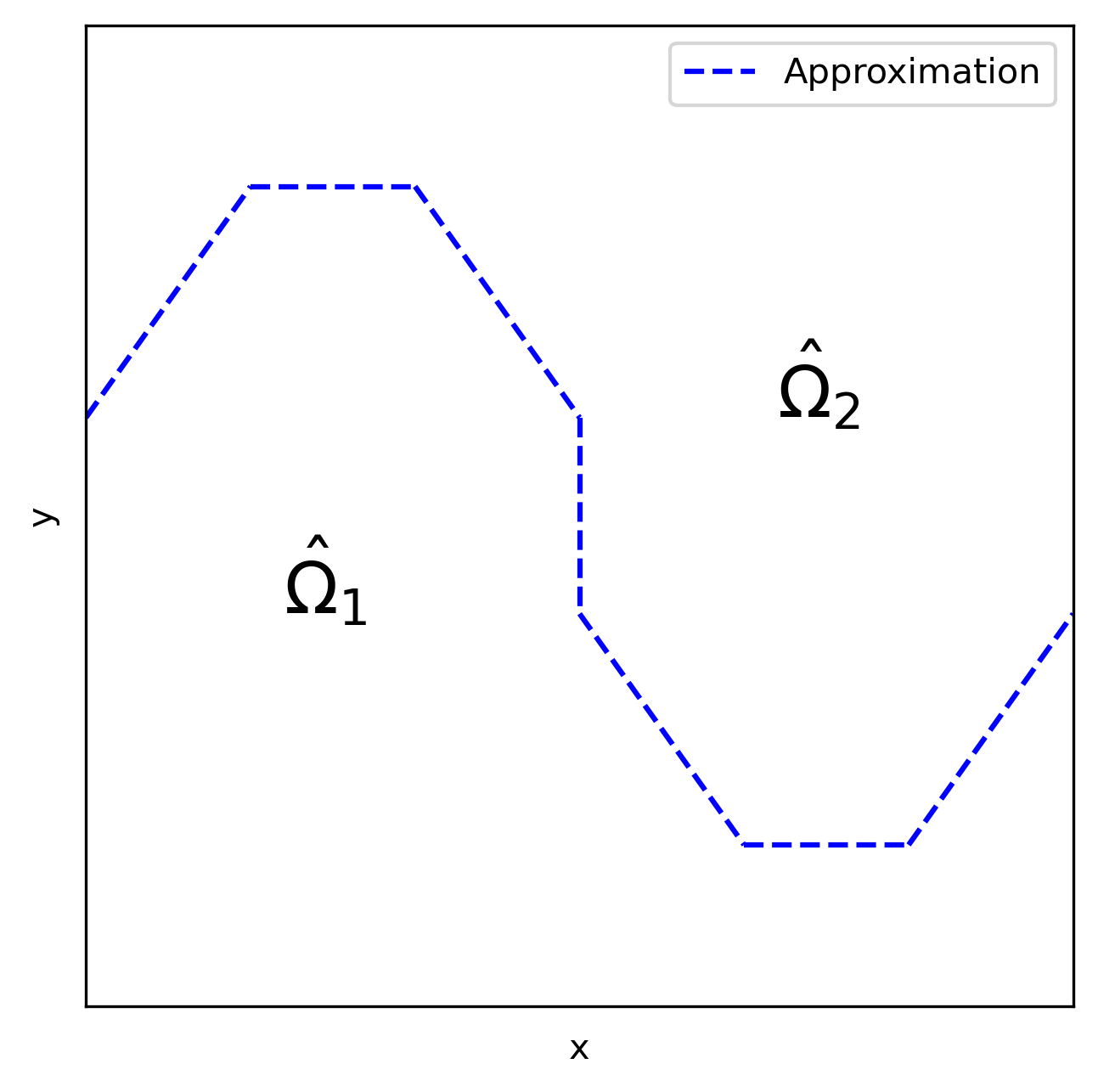}
\end{minipage}%
}%
\caption{An approximation of the interface $\Gamma$}
\end{figure}

Let $\Gamma=\partial\Omega_1\cap \partial\Omega_2$ be an interface in $C^0$. For any given $\varepsilon>0$, assume that there exists a connected series of hyperplanes 
\begin{equation}\label{HPs}
  \mathbf{a}_i\cdot\mathbf{x}-b_i=0\quad\mbox{for }\, i=1,\ldots, n  
\end{equation}
approximating the interface $\Gamma$ such that the hyperplanes divide the domain $\Omega$ by a partition $\left\{\hat{\Omega}_1,\hat{\Omega}_2\right\}$ (see Figures \ref{general12}, \ref{general22}, and \ref{general32}) and that 
\begin{equation}\label{epsilon}
    \big|\Omega_1\setminus \hat{\Omega}_1\big| + \big|\Omega_2\setminus \hat{\Omega}_2\big| \leq \varepsilon,
\end{equation}
where $\big|\Omega_i\setminus \hat{\Omega}_i\big|$ is the $d$-dimensional measure of $\Omega_i\setminus \hat{\Omega}_i$.
Let $\hat{\chi}$ be the indicator function of the subdomain $\hat{\Omega}_2$, i.e.,
\begin{equation}\label{step function hat}
 \hat{\chi}(\bx)=\left\{\begin{array}{rl}
 0, & \bx \in \hat{\Omega}_1,\\[2mm]
 1, & \bx \in \hat{\Omega}_2.
 \end{array}
 \right.   
\end{equation}
Then it is easy to see that \eqref{epsilon} implies that
\begin{equation}\label{chi chi-hat}
\|\chi-\hat{\chi}\|_{L^p(\Omega)} = \left(\big|\Omega_1\setminus \hat{\Omega}_1\big|  + \big|\Omega_2\setminus \hat{\Omega}_2\big|\right)^{1/p} \le \varepsilon^{1/p}.
\end{equation}

\begin{lemma}\label{general lemma}
Let $\hat{\Gamma}=\partial\hat{\Omega}_1\cap \partial\hat{\Omega}_2$. There exists a $d$--$n_1$--$n_2$--1 ReLU NN function $\mathcal{N}$ 
such that
\begin{equation}\label{general lp}
    \|\hat{\chi}-\mathcal{N}\|_{L^p(\Omega)}\le C(|\hat{\Gamma}|)\, \varepsilon^{1/p},
\end{equation}
where $n_1$ and $n_2$ are integers depending on, respectively, the number of the hyperplanes and convexity of $\hat{\Gamma}$, and $C(|\hat{\Gamma}|)$ is a positive constant depending on the $(d-1)$-dimensional measure of the interface $|\hat{\Gamma}|$.
\end{lemma}

\begin{proof}
The proof of the lemma is provided in Section \ref{general proof}. 
\end{proof}

\begin{theorem}\label{general theorem}
Under the assumption on the interface $\Gamma$, there exists a $d$--$n_1$--$n_2$--1 ReLU NN function $\mathcal{N}$ 
such that
\begin{equation}
    \|{\chi}-\mathcal{N}\|_{L^p(\Omega)}\le C(|\hat{\Gamma}|)\, \varepsilon^{1/p},
\end{equation}
where $n_1$ and $n_2$ are integers depending on, respectively, the number of the hyperplanes and convexity of $\hat{\Gamma}$, and $C(|\hat{\Gamma}|)$ is a positive constant depending on the $(d-1)$-dimensional measure of the interface $|\hat{\Gamma}|$.
\end{theorem}

\begin{proof}
It follows from \eqref{chi chi-hat}, Lemma \ref{general lemma} and the triangle inequality that there exists a $d$--$n_1$--$n_2$--1 ReLU NN function $\mathcal{N}$ for some $n_1,n_2\in\mathbb{N}$ such that
\begin{equation}
\|{\chi}-\mathcal{N}\|_{L^p(\Omega)}=\|{\chi}-\hat{\chi}+\hat{\chi}-\mathcal{N}\|_{L^p(\Omega)}\le\|{\chi}-\hat{\chi}\|_{L^p(\Omega)}+\|\hat{\chi}-\mathcal{N}\|_{L^p(\Omega)}\le \left(C(|\hat{\Gamma}|)+1\right)\, \varepsilon^{1/p},
\end{equation}
which completes the proof of the theorem.
\end{proof}

\section{Proof of Lemma \ref{general lemma} }\label{general proof}

This section proves Lemma \ref{general lemma} in Subsection \ref{convex section} and Subsection \ref{non-convex} when the subdomain $\hat{\Omega}_1$ is convex and non-convex, respectively.


\subsection{Convex $\hat{\Omega}_1$}\label{convex section}

\begin{figure}[htbp]
\centering
\subfigure[The interface $\hat{\Gamma}$ when $\hat{\Omega}_1$ is convex\label{convex interface approx1}]{
\begin{minipage}[t]{0.4\linewidth}
\centering
\includegraphics[width=1.8in]{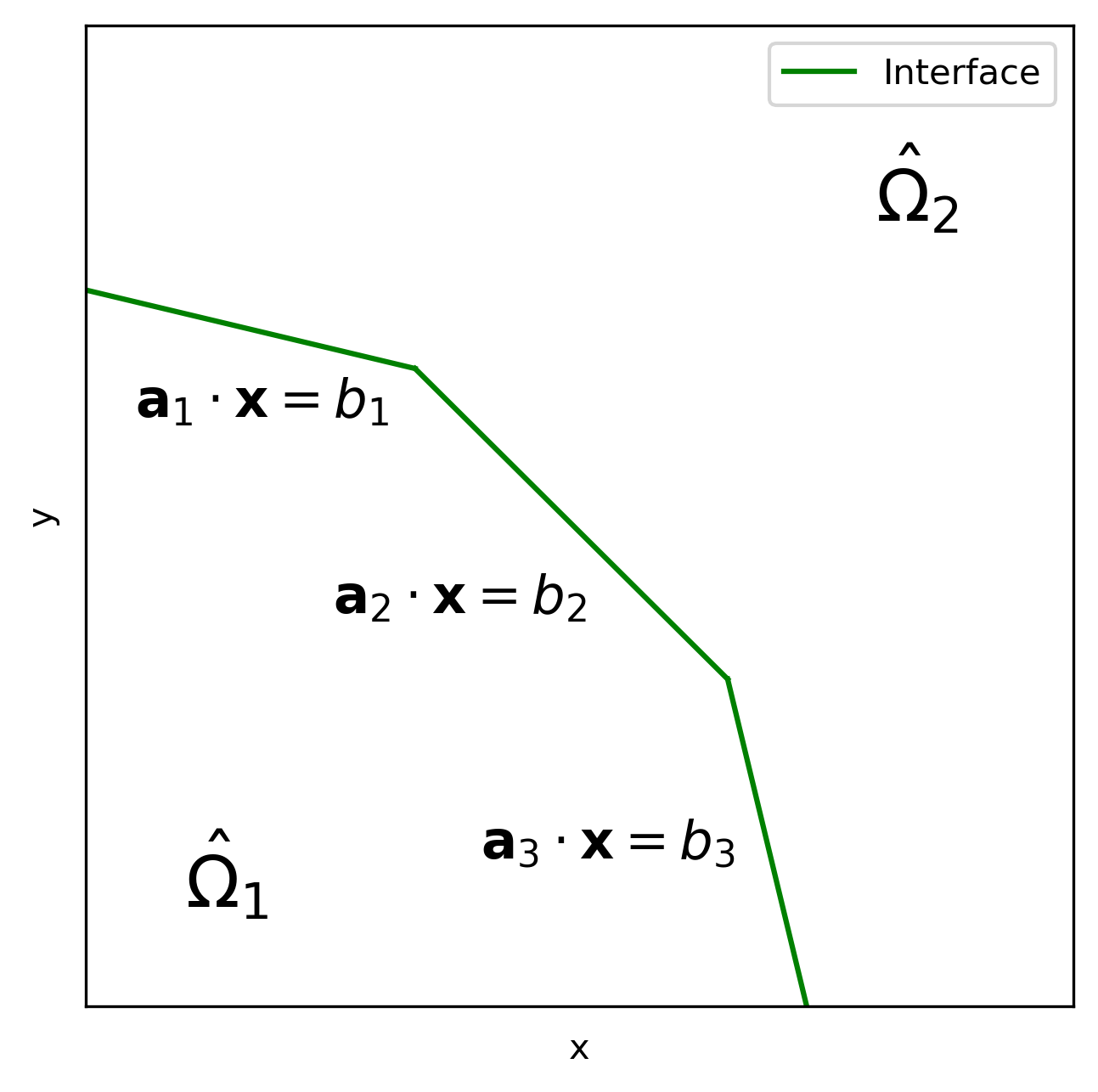}
\end{minipage}%
}%
\hspace{0.2in}
\subfigure[$\hat{\Omega}_{\varepsilon}$\label{convex interface approx12}]{
\begin{minipage}[t]{0.4\linewidth}
\centering
\includegraphics[width=1.8in]{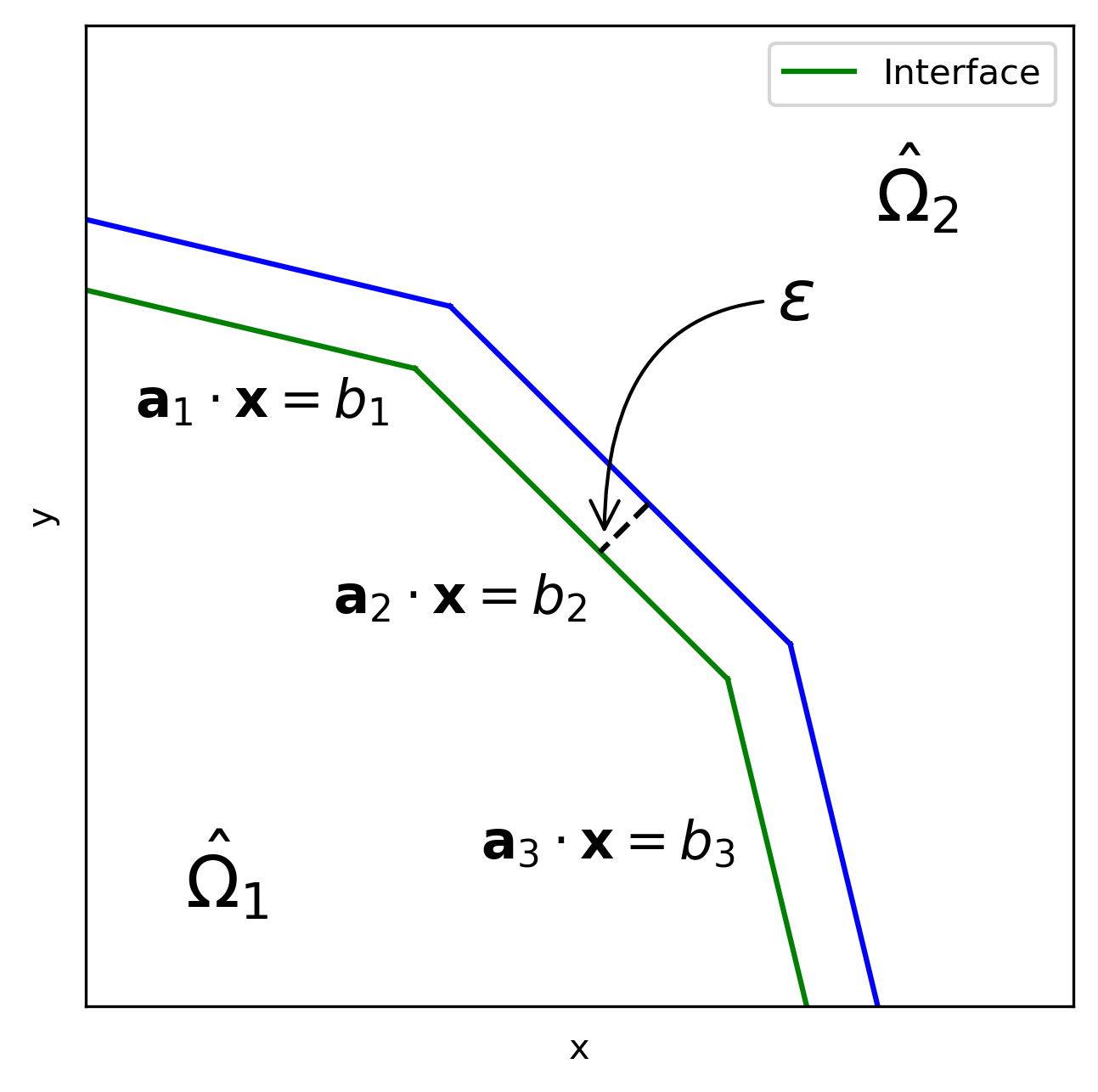}
\end{minipage}%
}%
\\
\subfigure[The region $\hat{\Omega}_{\varepsilon}$ is divided by the extension of $\mathbf{a}_i\cdot\mathbf{x}-b_i$.\label{extension of omega epsilon}]{
\begin{minipage}[t]{0.4\linewidth}
\centering
\includegraphics[width=1.8in]{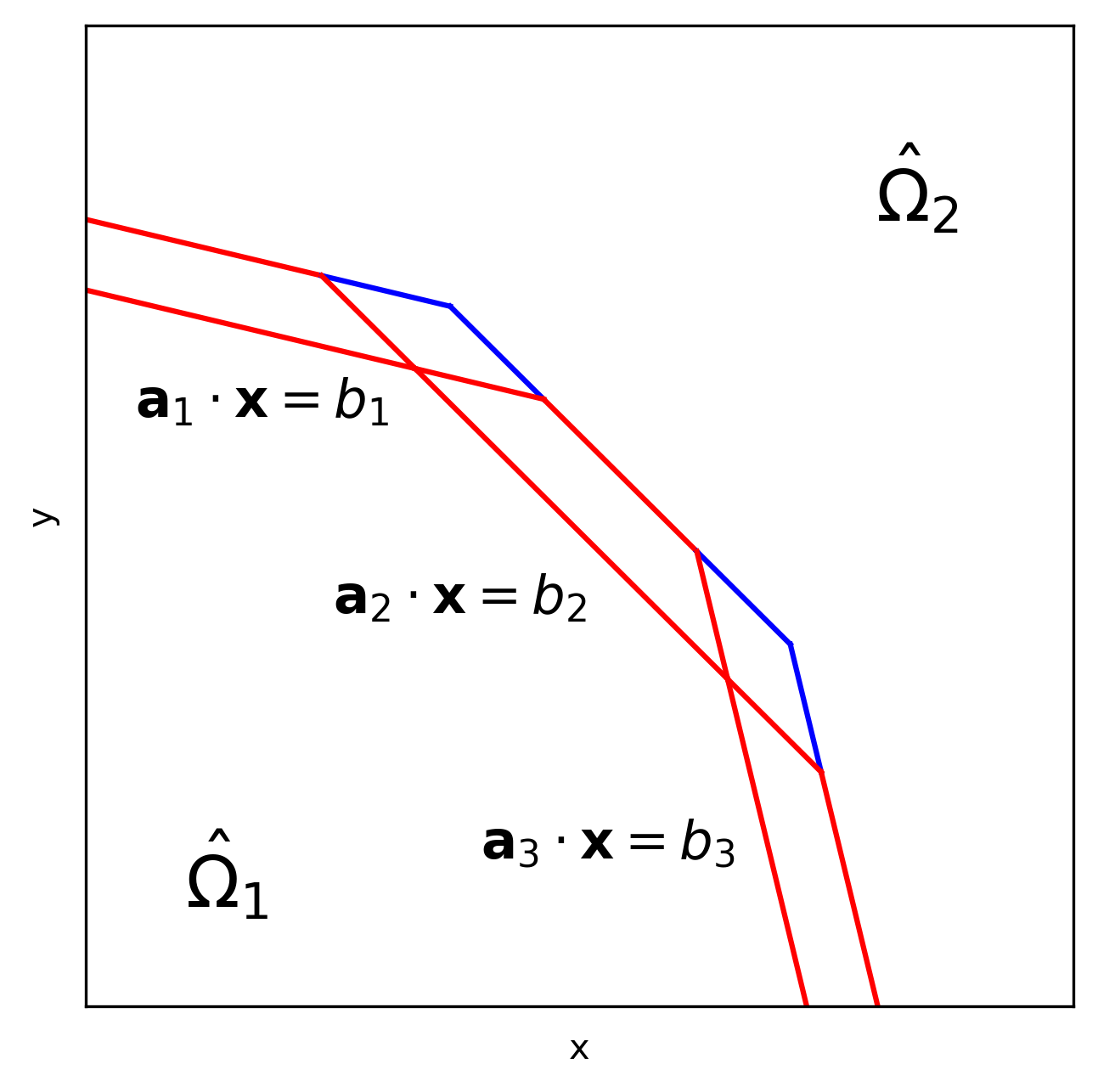}
\end{minipage}%
}%
\hspace{0.2in}
\subfigure[Subdividing each of the convex quadrilaterals with blue sides into two triangles\label{lines trough intersections}]{
\begin{minipage}[t]{0.4\linewidth}
\centering
\includegraphics[width=1.8in]{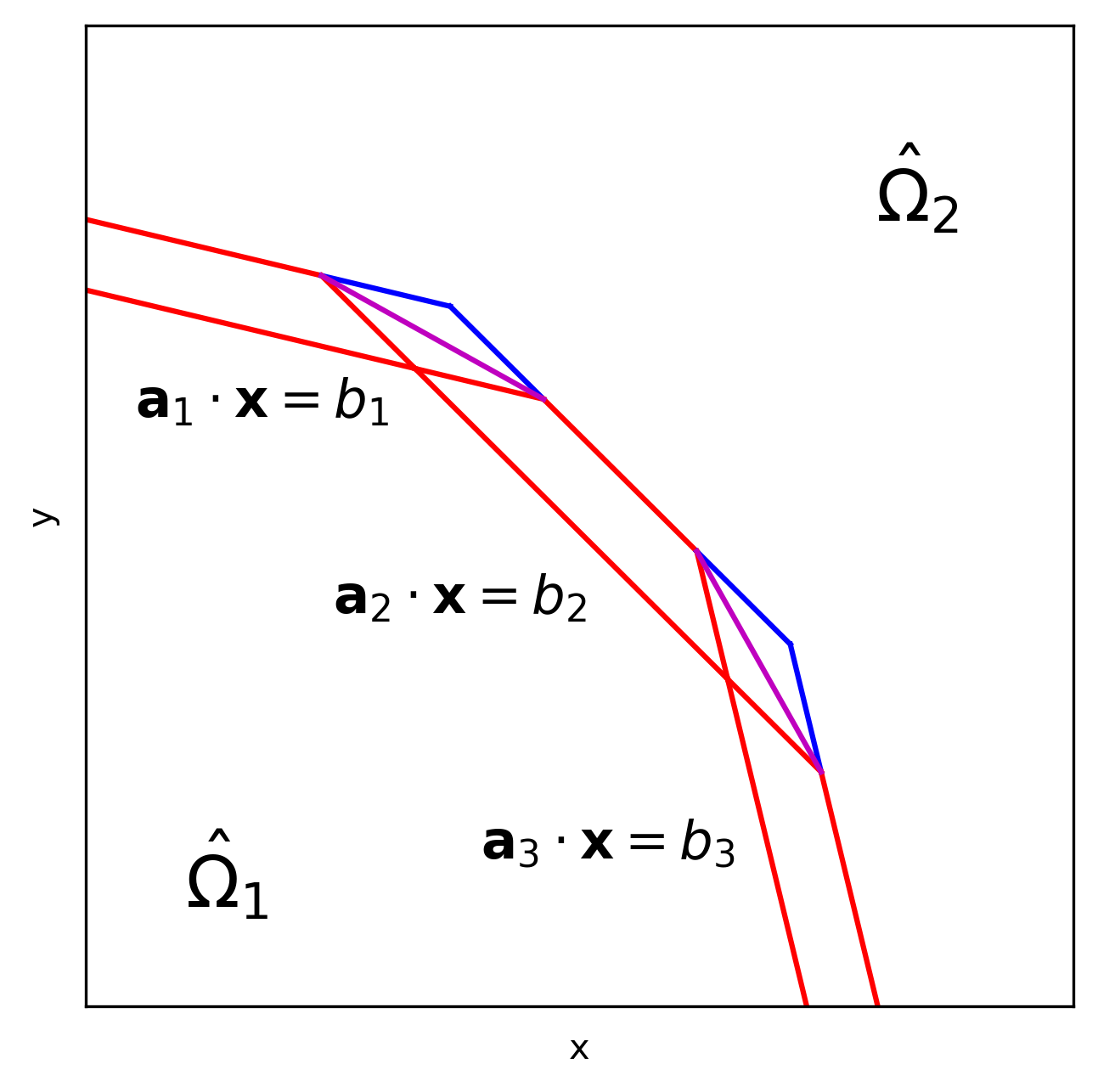}
\end{minipage}%
}%
\\
\subfigure[Removing the triangles not adjacent to the interface $\hat{\Gamma}$\label{remaing convex regions}]{
\begin{minipage}[t]{0.4\linewidth}
\centering
\includegraphics[width=1.8in]{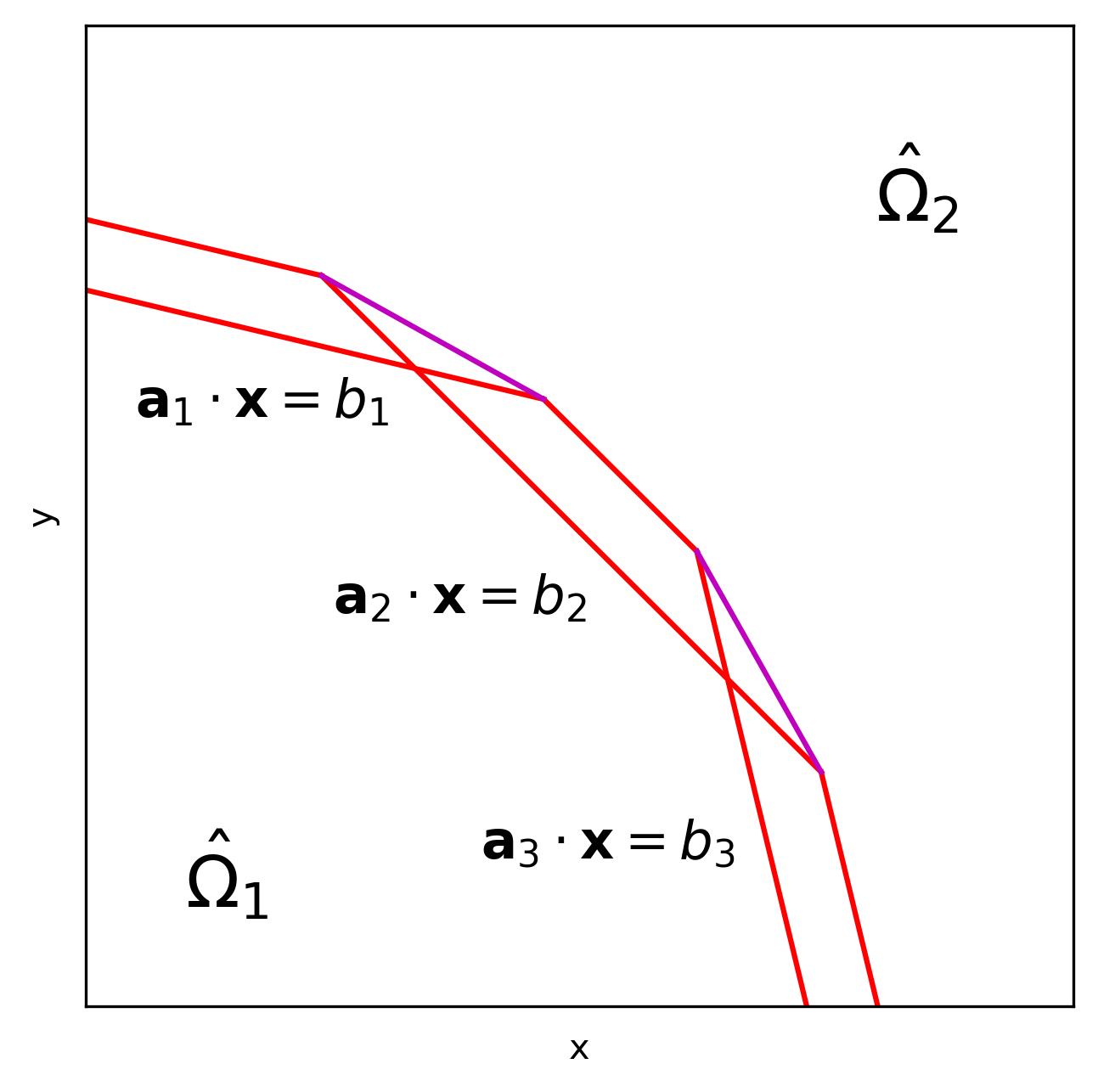}
\end{minipage}%
}%
\caption{The subdomain $\hat{\Omega}_1$ is convex.}
\end{figure}

This section shows the validity of Lemma \ref{general lemma} in a special case that the subdomain $\hat{\Omega}_1$ is convex (see Figure \ref{convex interface approx1}).

Without loss of generality, assume that the normal vectors $\mathbf{a}_i$ of the hyperplanes are the unit vectors and point toward $\hat{\Omega}_2$. Then we approximate the unit step function $\hat{\chi}(\bx)$ in \eqref{step function hat} by the following ReLU NN function
\begin{equation}\label{convex analytic}
\mathcal{N}(\mathbf{x})=1-\sigma\left(1-\frac{1}{\varepsilon}\sum_{i=1}^{n}\sigma(\mathbf{a}_i\cdot\mathbf{x}-b_i)\right).    
\end{equation}
The $\mathcal{N}(\mathbf{x})$ is a $d$--$n$--$1$--$1$ ReLU NN function.

When the interface $\hat{\Gamma}$ is a hyperplane $\mathbf{a}\cdot\mathbf{x}-b=0$ in $\R^d$, i.e., $n=1$, the $\mathcal{N}(\mathbf{x})$ has the form
\begin{equation}\label{N-1}
\mathcal{N}(\mathbf{x})=1-\sigma\left(1-\frac{1}{\varepsilon}\sigma(\mathbf{a}\cdot\mathbf{x}-b)\right). 
\end{equation}
The second term of $\mathcal{N}(\mathbf{x})$, a three-layer ReLU NN function, is a ramp function that equals negative one in $\hat{\Omega}_1$ and vanishes in $\hat{\Omega}_2\setminus Y_\varepsilon$, where $Y_\varepsilon=\left\{\bx\in \Omega :\, 0<\mathbf{a}\cdot\mathbf{x}-b< \varepsilon  \right\}$ is a strip with $\varepsilon$-width. It is then easy to see that
\[
 \hat{\chi}(\bx)-\mathcal{N}(\mathbf{x})= \left\{\begin{array}{ll}
 0, & \bx \in \hat{\Omega}_1 \cup \left(\hat{\Omega}_2\setminus Y_\varepsilon \right),\\[2mm]
 \sigma\left(1-\frac{1}{\varepsilon}\sigma(\mathbf{a}\cdot\mathbf{x}-b)\right), & \bx \in Y_\varepsilon ,
 \end{array}
 \right.   
\]
which, together with a simple calculation, implies the upper bound in \eqref{general lp}.

Now, we consider the case $n \ge 2$. For simplicity of presentation, the proof of the error bound in \eqref{general lp} is carried out in two dimensions $d=2$.
Denote by $\hat{\Omega}_{\varepsilon}\subset\Omega$ the region produced by translating $\mathbf{a}_i\cdot\mathbf{x}-b_i=0$ toward $\hat{\Omega}_2$ along $\mathbf{a}_i$ by $\varepsilon$ (see Figure \ref{convex interface approx12}). By extending the line segments $\mathbf{a}_i\cdot\mathbf{x}-b_i=0$, we partition the region $\hat{\Omega}_{\varepsilon}$ into convex subregions (see Figure \ref{extension of omega epsilon}).

The subregions of the first type are denoted by $\left\{\Upsilon_{1i}\right\}_{i=1}^n$, where $\Upsilon_{1i}$ is the subregion bounded by the line $\mathbf{a}_i\cdot\mathbf{x}-b_i=0$, its translated line $\mathbf{a}_i\cdot \mathbf{x}-b_i=\varepsilon$, and two neighboring lines or one neighboring line and the boundary of $\Omega$ (the convex quadrilaterals with red sides in Figure \ref{extension of omega epsilon}). More precisely, we have that
\begin{eqnarray*}
 \Upsilon_{11} &=&\left\{\bx\in\hat{\Omega}_{\varepsilon}: 
    0<\mathbf{a}_1\cdot \mathbf{x}-b_1<\varepsilon\, \mbox{ and }\,
    \mathbf{a}_{2}\cdot \mathbf{x}-b_{2}<0\right\}, \\[2mm]
  \Upsilon_{1n} &=&\left\{\bx\in\hat{\Omega}_{\varepsilon}: 
    0<\mathbf{a}_n\cdot \mathbf{x}-b_n<\varepsilon\, \mbox{ and }\,
    \mathbf{a}_{n-1}\cdot \mathbf{x}-b_{n-1}<0\right\}, 
\end{eqnarray*}
and that for $i=2,\ldots,n-1$ 
\[ 
    \Upsilon_{1i}=\left\{\bx\in\hat{\Omega}_{\varepsilon}: 
    0<\mathbf{a}_i\cdot \mathbf{x}-b_i<\varepsilon,\,
    \mathbf{a}_{i-1}\cdot \mathbf{x}-b_{i-1}<0,\, \mbox{ and }\,
    \mathbf{a}_{i+1}\cdot \mathbf{x}-b_{i+1}<0\right\}.
\]

Notice that $\hat{\Omega}_\varepsilon\setminus \left(\cup_{i=1}^{n} \Upsilon_{1i}\right)$ consists of $n-1$ convex quadrilaterals (with blue sides in Figure \ref{extension of omega epsilon}). Subdividing each of these convex quadrilaterals into two triangles (see Figure \ref{lines trough intersections}) and removing the triangles not adjacent to the interface $\hat{\Gamma}$ (see Figure \ref{remaing convex regions}), the remaining triangles are denoted by $\left\{\Upsilon_{2i}\right\}_{i=1}^{n-1}$, where $\Upsilon_{2i}$ is given by
\[
\Upsilon_{2i}=\left\{\bx\in\hat{\Omega}_\varepsilon : \mathbf{a}_j\cdot \mathbf{x}-b_j>0 \mbox{ for } j=i,i+1,\, \mbox{ and }
\left(\mathbf{a}_i + \mathbf{a}_{i+1}\right) \cdot \mathbf{x}-\left(b_i+b_{i+1}\right)<\varepsilon \right\}.
\]
We then have the following lemma.

\begin{lemma}\label{error lemma}
Let $\mathcal{N}(\mathbf{x})$ be the three-layer ReLU NN function defined in \eqref{convex analytic}, then we have
    \begin{equation}\label{error}
    \hat{\chi}(\bx)-\mathcal{N}(\mathbf{x})= \left\{\begin{array}{ll}
 0, & \bx \in \Omega \setminus \left(\bigcup\limits_{j=1}^2\bigcup\limits_{i=1}^{n+1-j}\Upsilon_{ji}\right),\\[5mm]
 \hat{\chi}(\mathbf{x})-\dfrac{1}{\varepsilon}(\mathbf{a}_i\cdot\mathbf{x}-b_i), & \bx \in \Upsilon_{1i} \mbox{ for } i=1,\ldots n,\\[4mm]
  \hat{\chi}(\mathbf{x})-\dfrac{1}{\varepsilon}\left[(\mathbf{a}_i+\mathbf{a}_{i+1})\cdot\mathbf{x}-(b_i+b_{i+1})\right], 
  & \bx \in \Upsilon_{2i} \mbox{ for } i=1,\ldots n-1.
 \end{array}
 \right.   
\end{equation}
\end{lemma}

\begin{proof}
Let 
\[
\hat{\Omega}_3=\hat{\Omega}_2 \setminus  \left(\bigcup\limits_{j=1}^2\bigcup\limits_{i=1}^{n+1-j}\Upsilon_{ji}\right).
\]
Since $\mathbf{a}_i$ points toward $\hat{\Omega}_2$, clearly, we have $\sigma(\mathbf{a}_i\cdot\mathbf{x}-b_i)=0$ for all $\bx\in \hat{\Omega}_1$ and $i=1,\ldots,n$. This implies
\begin{equation}\label{Omega3}
    \mathcal{N}(\bx)=1-\sigma(1)=0, \quad\forall \,\,\bx\in \hat{\Omega}_1.
\end{equation}
Clearly, we have
\[
1-\frac{1}{\varepsilon}\sum_{i=1}^{n}\sigma(\mathbf{a}_i\cdot\mathbf{x}-b_i) = \left\{\begin{array}{ll}
 1-\dfrac{1}{\varepsilon}(\mathbf{a}_i\cdot\mathbf{x}-b_i), & \bx\in \Upsilon_{1i},\\[4mm]
1 -\dfrac{1}{\varepsilon}\left[(\mathbf{a}_i+\mathbf{a}_{i+1})\cdot\mathbf{x}-(b_i+b_{i+1})\right], & \bx\in \Upsilon_{2i}.
\end{array}
 \right.   
\]
It is easy to see that
\[
1-\dfrac{1}{\varepsilon}(\mathbf{a}_i\cdot\mathbf{x}-b_i) \left\{\begin{array}{ll}
>0, & 0<\mathbf{a}_i\cdot\mathbf{x}-b_i< \varepsilon,\\ [2mm]
\leq 0, & \varepsilon\leq \mathbf{a}_i\cdot\mathbf{x}-b_i
\end{array}
 \right. 
\]
and that similar inequalities hold for $1 -\dfrac{1}{\varepsilon}\left[(\mathbf{a}_i+\mathbf{a}_{i+1})\cdot\mathbf{x}-(b_i+b_{i+1})\right]$; furthermore, by the definition of $\hat{\Omega}_3$, we have
\[
1-\frac{1}{\varepsilon}\sum_{i=1}^{n}\sigma(\mathbf{a}_i\cdot\mathbf{x}-b_i) <0, \quad\forall \,\, \bx\in \hat{\Omega}_3.
\]
Now, applying the activation function $\sigma$, multiplying by $-1$, and adding 1 imply 
\[ 
     \mathcal{N}(\bx) = \left\{\begin{array}{ll}
 (\mathbf{a}_i\cdot\mathbf{x}-b_i)/\varepsilon, & \bx\in \Upsilon_{1i},\\[2mm]
\left[(\mathbf{a}_i+\mathbf{a}_{i+1})\cdot\mathbf{x}-(b_i+b_{i+1})\right]/\varepsilon, & \bx\in \Upsilon_{2i},
 \\[2mm]
 1 & \bx\in\hat{\Omega}_3,
 \end{array}
 \right.   
\] 
which, together with \eqref{Omega3}, leads to \eqref{error}. This completes the proof of the lemma.
\end{proof}

\begin{proof}[Proof of {\em Lemma \ref{general lemma}} for convex $\hat{\Omega}_1$]
When $\hat{\Omega}_1$ is convex, to show the validity of Lemma \ref{general lemma}, notice that for all $p\in [1,\infty)$, we have by Lemma \ref{error lemma},
\begin{eqnarray*}
\left|\hat{\chi}(\mathbf{x})-\mathcal{N}({\mathbf{x}})\right|^p 
&=& \left|\hat{\chi}(\mathbf{x})-\frac{1}{\varepsilon}(\mathbf{a}_i\cdot\mathbf{x}-b_i)\right|^p\le 1, \quad\forall \,\bx\in \Upsilon_{1i},\\ [2mm]
\mbox{and }\, \left|\hat{\chi}(\mathbf{x})-\mathcal{N}({\mathbf{x}})\right|^p
&=& \left|\hat{\chi}(\mathbf{x})-\dfrac{1}{\varepsilon}\left[(\mathbf{a}_i+\mathbf{a}_{i+1})\cdot\mathbf{x}-(b_i+b_{i+1})\right]\right|^p \le 1, \quad\forall \,\bx\in \Upsilon_{2i},
\end{eqnarray*}
which implies
\begin{equation}\label{a}
 \|\hat{\chi}-\mathcal{N}\|^p_{L^p(\Upsilon_{1i})}\le \big|\Upsilon_{1i}\big| \quad\mbox{and}\quad  \|\hat{\chi}-\mathcal{N}\|^p_{L^p(\Upsilon_{2i})}\le \big|\Upsilon_{2i}\big|,
\end{equation}
where $\big|\Upsilon_{ji}\big|$ denotes the area of the quadrilateral $\Upsilon_{ji}$. It follows from \eqref{error} and \eqref{a} that
\begin{equation}\label{4.6}
    \|\hat{\chi}-\mathcal{N}\|_{L^p(\Omega)}^p
=\sum_{i=1}^{n}\|\hat{\chi}-\mathcal{N}\|_{L^p(\Upsilon_{1i})}^p+\sum_{i=1}^{n-1}\|\hat{\chi}-\mathcal{N}\|_{L^p(\Upsilon_{1i})}^p
\le \sum_{i=1}^{n}\big|\Upsilon_{1i}\big|+\sum_{i=1}^{n-1}\big|\Upsilon_{2i}\big| \leq \big|\hat{\Omega}_\varepsilon\big|,
\end{equation}
which, together with the fact that $\big|\hat{\Omega}_\varepsilon\big|\leq C\, \big|\hat{\Gamma}\big|\, \varepsilon$ for a positive constant $C$, implies the error bound in \eqref{general lp}. This completes the proof of Lemma \ref{general lemma}.
\end{proof}

\subsection{Non-Convex $\hat{\Omega}_1$}\label{non-convex}

\begin{figure}[htbp]
\centering
\subfigure[A non-convex $\hat{\Omega}_1$\label{nonconvex}]{
\begin{minipage}[t]{0.4\linewidth}
\centering
\includegraphics[width=1.8in]{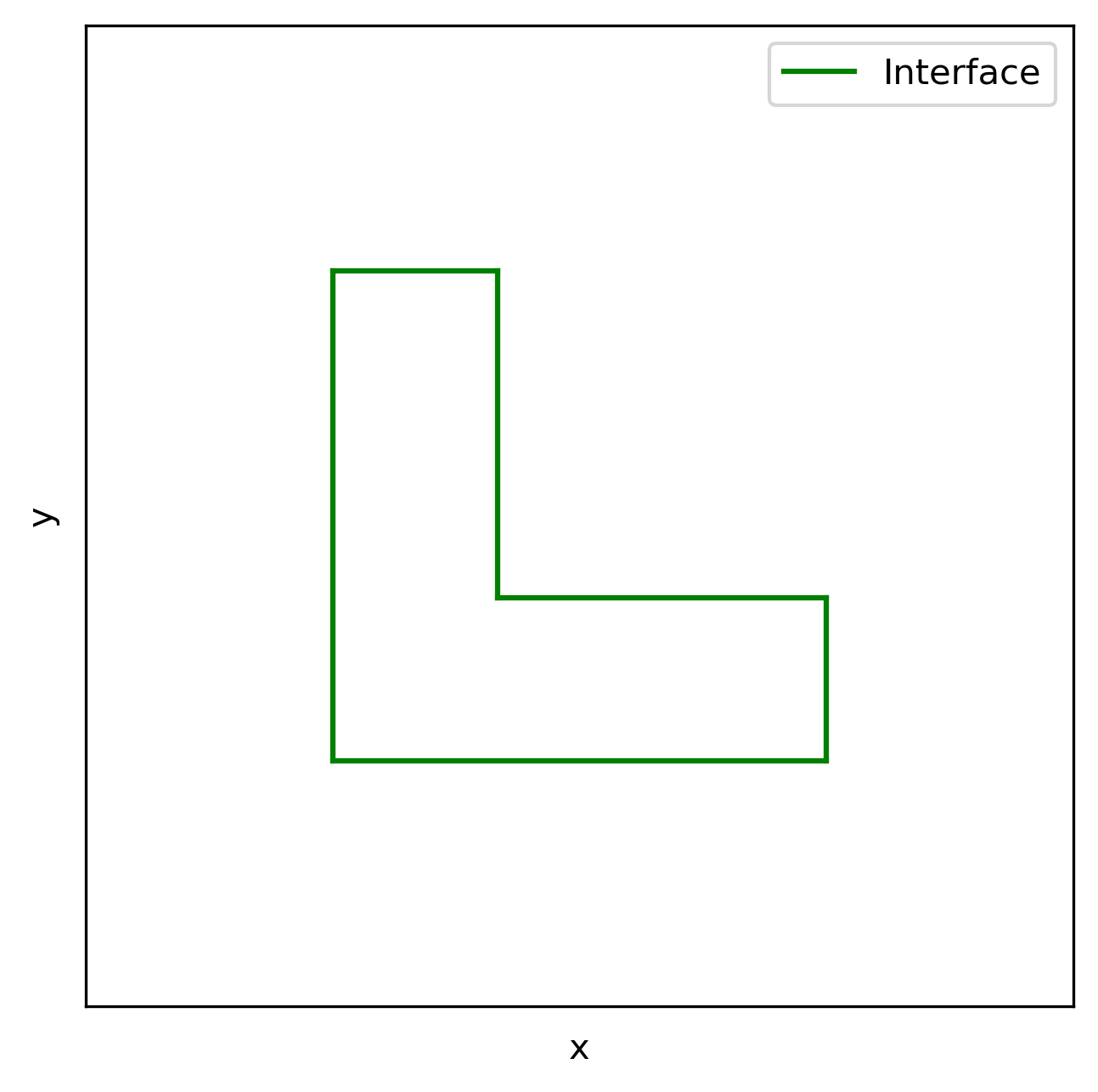}
\end{minipage}%
}%
\hspace{0.2in}
\subfigure[Thue convex hull of $\hat{\Omega}_1$\label{convex hull}]{
\begin{minipage}[t]{0.4\linewidth}
\centering
\includegraphics[width=1.8in]{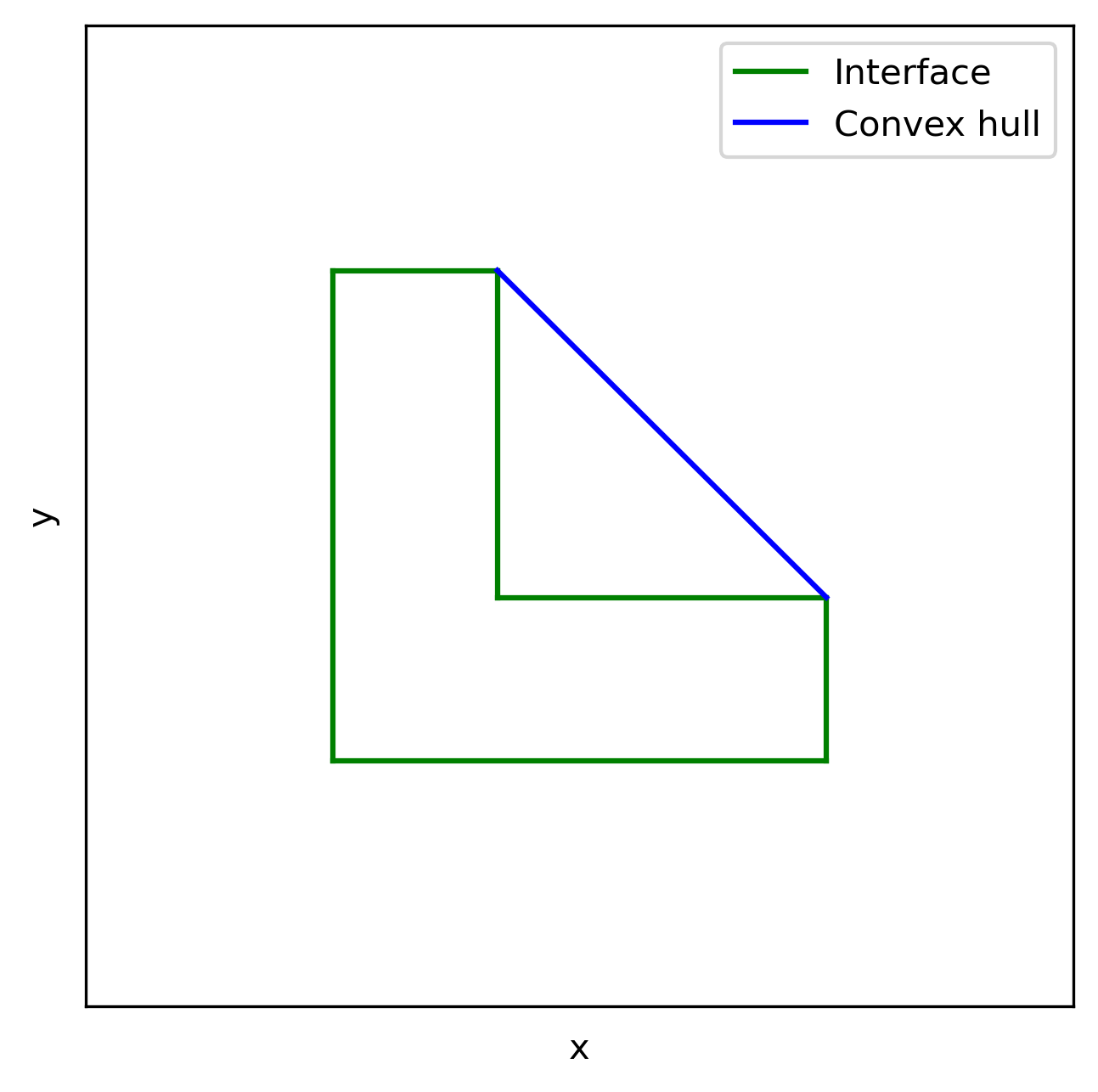}
\end{minipage}%
}%
\\
\subfigure[A convex decomposition of $\hat{\Omega}_1$\label{convex decompistion}]{
\begin{minipage}[t]{0.4\linewidth}
\centering
\includegraphics[width=1.8in]{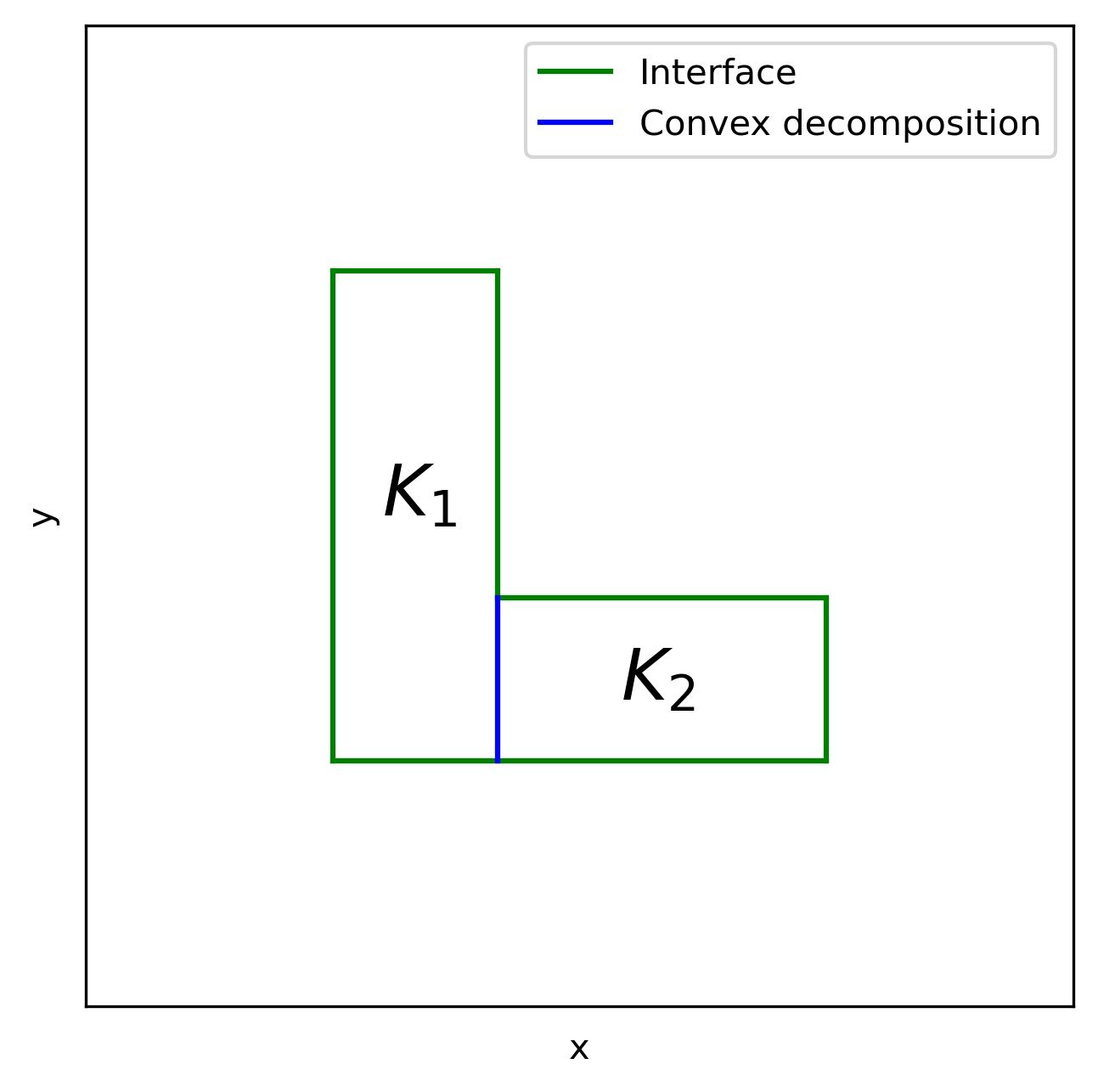}
\end{minipage}%
}%
\hspace{0.2in}
\subfigure[A subset of $K_2$\label{convex decompistion2}]{
\begin{minipage}[t]{0.4\linewidth}
\centering
\includegraphics[width=1.8in]{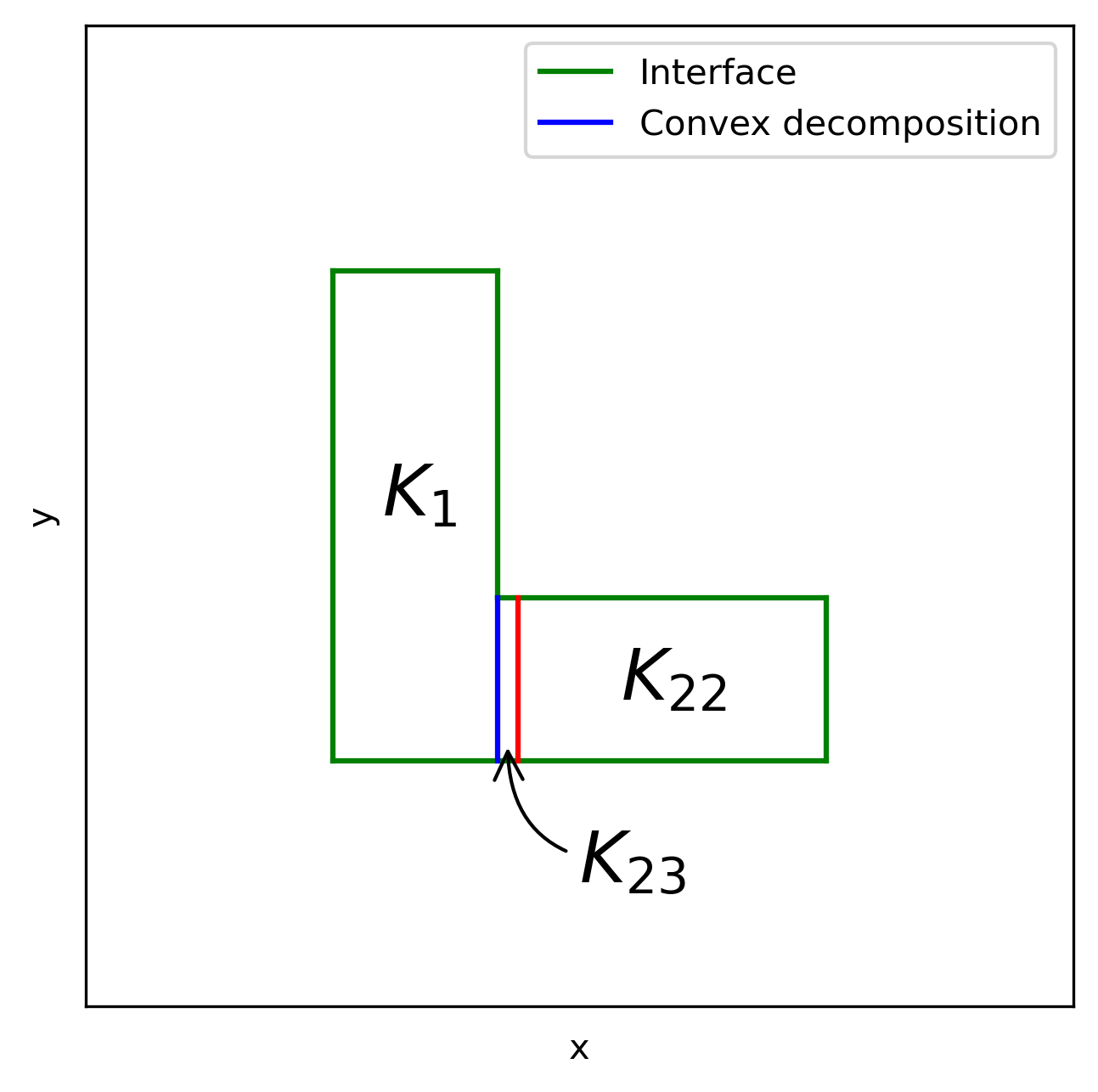}
\end{minipage}%
}%
\caption{The subdomain $\hat{\Omega}_1$ is non-convex.}
\end{figure}

This section shows the validity of Lemma \ref{general lemma} when $\hat{\Omega}_1$ is non-convex (see, e.g., Figure \ref{nonconvex}).  Our proof is again through explicit constructions. Specifically, we 
present two approaches: one is based on the convex hull of $\hat{\Omega}_1$ (see Subsubsection \ref{convex hull subsubsection}) and the other uses a convex decomposition of $\hat{\Omega}_1$ (see Subsubsection \ref{convex decomposition subsubsection}). 

\subsubsection{Convex hull}\label{convex hull subsubsection}
Let
\begin{equation}\label{Ki}
   \Omega_1^{(1)}=\hat{\Omega}_1\cup \left(\cup_{i=1}^k K_i\right), 
\end{equation}
be the convex hull of $\hat{\Omega}_1$ (see Figure \ref{convex hull}) generated by a convex hull algorithm (see, e.g., \cite{mount2002cmsc,barber1996quickhull,avis1995good}),
where $K_i$ are polytopes and pairwise disjoint. Without loss of generality, we assume that all $K_i$ ($i=1,\ldots,k$) are convex. Otherwise, the procedure presented in this section may be applied to non-convex $K_i$s for the indicator functions ${\bm 1}_{\Omega\setminus K_i}(\bx)$ of the subdomains $\Omega\setminus K_i$. Note that the procedure may be needed for several times recursively.

Let $\hat{\chi}_0(\bx)$ be the unit step function defined on the convex hull $\Omega_1^{(1)}$ of the non-convex subdomain $\hat{\Omega}_1$:
\begin{equation}\label{chi0}
 \hat{\chi}_0(\bx)=\left\{\begin{array}{rl}
 0, & \bx \in \Omega_1^{(1)}\subset\Omega,\\[2mm]
 1, & \bx \in \Omega\setminus \Omega_1^{(1)},
 \end{array}
 \right.
\end{equation}
then its discontinuity interface is $\hat{\Gamma}_0=\partial \Omega_1^{(1)} \cap \partial \left(\Omega\setminus \Omega_1^{(1)}\right)$ consisting of $n_{1,0}$ faces. As proved in Subsection \ref{convex section} (see \eqref{4.6}), there exists a $d$--$n_{1,0}$--$1$--1 ReLU NN function approximation $\mathcal{N}_0(\bx)$ such that \begin{equation}\label{chi 0}
    \|\hat{\chi}_0-\mathcal{N}_0\|_{L^p(\Omega)}=\|\hat{\chi}_0-\mathcal{N}_0\|_{L^p(\hat{\Omega}_{\varepsilon, 0})}\le \big|\hat{\Omega}_{\varepsilon, 0} \big|^{1/p},
\end{equation}
where $\hat{\Omega}_{\varepsilon, 0}$ is the region with $\varepsilon$-width containing the interface $\hat{\Gamma}_0$ by translating the faces of $\hat{\Gamma}_0$ towards the subdomain $\Omega\setminus\Omega_1^{(1)}$.

For each convex polytope $K_i$ ($i=1,\ldots, k$) in \eqref{Ki} having $n_{1,i}$ faces, let ${\hat{\chi}}_i(\bx)$ be the unit step function defined on ${K}_i$:
\[
 {\hat{\chi}}_i(\bx)=\left\{\begin{array}{rl}
 1, & \bx \in {K}_i\subset\Omega,\\[2mm]
 0, & \bx \in \Omega\setminus {K}_i.
 \end{array}
 \right.
\]
Define the following $d$--$n_{1,i}$--$1$--1 ReLU NN function
\begin{equation}\label{Ni}
    \mathcal{N}_i(\bx)=\sigma\left(1-\frac{1}{\varepsilon}\sum_{j=1}^{n_{1,i}}\sigma(\mathbf{a}_{i,j}\cdot\mathbf{x}-b_{i,j})\right),
\end{equation}
where the hyperplanes $\mathbf{a}_{i,j}\cdot\mathbf{x}-b_{i,j}=0$ ($j=1,\ldots, n_{1,i}$) are the faces of $\partial {K}_i$ with $\mathbf{a}_{i,j}$ pointing toward $\Omega\setminus {K}_i$. In a similar fashion as in Subsection \ref{convex section}, it is easy to check that
\begin{equation}\label{chi i}
    \|{\hat{\chi}}_i-\mathcal{N}_i\|_{L^p(\Omega)} = \|{\hat{\chi}}_i-\mathcal{N}_i\|_{L^p(\hat{\Omega}_{\varepsilon, i})} \le \big| \hat{\Omega}_{\varepsilon, i} \big|^{1/p},
\end{equation}
where $\hat{\Omega}_{\varepsilon, i}$ is a region having $\varepsilon$-width.

Now, we are ready to define the following $d$--$n_{1}$--$n_{2}$--1 ReLU NN function:
\begin{equation}\label{N}
    \mathcal{N}(\bx)=\mathcal{N}_0(\bx)+\sum_{i=1}^k\mathcal{N}_i(\bx),
\end{equation}
where $\mathcal{N}_0(\bx)$ is given in a similar fashion as in  \eqref{convex analytic}.

\begin{proof}[Proof of {\em Lemma \ref{general lemma}} for non-convex $\hat{\Omega}_1$]
Note that
\[\hat{\chi}=\sum_{i=0}^k\hat{\chi}_i.\]
Then it follows from \eqref{N}, the triangle inequality, \eqref{chi 0}, and \eqref{chi i} that
\begin{equation}\label{convex hull inequality}
\|\hat{\chi}-\mathcal{N}\|_{L^p(\Omega)}\le \sum_{i=0}^k\|\hat{\chi}_i-\mathcal{N}_i\|_{L^p\left(\Omega\right)}\le \sum_{i=0}^k\big| \hat{\Omega}_{\varepsilon, i} \big|^{1/p},
\end{equation}
which, together with the fact that $\big|\hat{\Omega}_{\varepsilon,i}\big|\leq C_i\, \big|\hat{\Gamma}_i\big|\, \varepsilon$ for a positive constant $C_i$, implies the error bound in \eqref{general lp}. Here $\hat{\Gamma}_i=\partial {K}_i \cap \partial \left(\Omega\setminus {K}_i\right)$ for $i=1,\ldots, k$ This completes the proof of Lemma \ref{general lemma}.
\end{proof}

\subsubsection{Convex decomposition}\label{convex decomposition subsubsection}

Assume that $\hat{\Omega}_1$ has a convex decomposition (see, e.g., \cite{hertel1983fast}) given by
\[ 
   \hat{\Omega}_1=\bigcup_{i=1}^l K_i,
\] 
where all $K_i$ are convex polytopes. For simplicity of presentation, assume that $l=2$, i.e., the decomposition has only two convex polytopes: $\hat{\Omega}_1=K_1\cup K_2$ (see Figure \ref{convex decompistion}). 

Denote the indicator function of the subdomain $K_1$ by
\[
 {\hat{\chi}}_1(\bx)=\left\{\begin{array}{rl}
 1, & \bx \in {K}_1\subset\Omega,\\[2mm]
 0, & \bx \in \Omega\setminus {K}_1.
 \end{array}
 \right.
\]
Let $\mathbf{a}\cdot\mathbf{x}-b=0$ be the hyperplane that divides $\hat{\Omega}_1$ into $K_1$ and $K_2$ (blue line in Figure \ref{convex decompistion}). Assume that $\mathbf{a}$ points toward $K_2$. Translate $\mathbf{a}\cdot\mathbf{x}-b=0$ toward $K_2$ by $\varepsilon$ to obtain the hyperplane $\mathbf{a}\cdot\mathbf{x}-b-\varepsilon=0$ (red line in Figure \ref{convex decompistion2}). Partition $K_2$ by $\{K_{22}, K_{23}\}$ (see Figure \ref{convex decompistion2}), where 
\[
K_{22}=\{\mathbf{x}\in\hat{\Omega}_1:\varepsilon<\mathbf{a}\cdot\mathbf{x}-b\} \quad\mbox{and}\quad K_{23}=\{\mathbf{x}\in\hat{\Omega}_1:0<\mathbf{a}\cdot\mathbf{x}-b<\varepsilon\}.
\]
Denote by ${\hat{\chi}}_{22}(\bx)$ and ${\hat{\chi}}_{23}(\bx)$ the respective indicator functions of $K_{22}$ and $K_{23}$:
\[
 {\hat{\chi}}_{22}(\bx)=\left\{\begin{array}{rl}
 1, & \bx \in {K}_{22}\subset\Omega,\\[2mm]
 0, & \bx \in \Omega\setminus {K}_{22}.
 \end{array}
 \right.
 \quad\mbox{and}\quad {\hat{\chi}}_{23}(\bx)=\left\{\begin{array}{rl}
 1, & \bx \in {K}_{23}\subset\Omega,\\[2mm]
 0, & \bx \in \Omega\setminus {K}_{23}.
 \end{array}
 \right.
\]

Assume that polygonal domains $K_1$ and $K_{22}$ have $n_1$ and $n_2$ faces, respectively. In a similar fashion as in \eqref{Ni} and \eqref{chi i}, there exist $d$--$n_{1}$--$1$--1 and $d$--$n_{2}$--$1$--1 ReLU NN functions $\mathcal{N}_1$ and $\mathcal{N}_{22}$ such that
\begin{equation}\label{chi 1}
\left\{\begin{array}{l}
    \qquad\,\, \|{\hat{\chi}}_1-\mathcal{N}_1\|_{L^p(\Omega)} = \|{\hat{\chi}}_1-\mathcal{N}_1\|_{L^p(\hat{\Omega}_{\varepsilon, 1})} \le \big| \hat{\Omega}_{\varepsilon, 1} \big|^{1/p} \\[2mm]
     \mbox{and }\,\,   \|{\hat{\chi}}_{22}-\mathcal{N}_{22}\|_{L^p(\Omega)} = \|{\hat{\chi}}_{22}-\mathcal{N}_{22}\|_{L^p(\hat{\Omega}_{\varepsilon, 22})} \le \big| \hat{\Omega}_{\varepsilon, 22} \big|^{1/p},
    \end{array}
    \right.
\end{equation}
where $\hat{\Omega}_{\varepsilon, 1}$ and $\hat{\Omega}_{\varepsilon, 22}$ are regions having $\varepsilon$-width. Clearly, there exist positive constants $C_1$, $C_{22}$, and $C_{23}$ such that
\begin{equation}\label{4.15}
\big|\hat{\Omega}_{\varepsilon,1}\big|\leq C_1\, \big|\hat{\Gamma}_1\big|\, \varepsilon, \quad \big|\hat{\Omega}_{\varepsilon,22}\big|\leq C_{22}\, \big|\hat{\Gamma}_{22}\big|\, \varepsilon, \quad\mbox{and}\quad \big|K_{23}\big|\leq C_{23}\, \big|\hat{\Gamma}_{23}\big|\, \varepsilon  
\end{equation}
where $\hat{\Gamma}_1$, $\hat{\Gamma}_{22}$, and $\hat{\Gamma}_{23}$ are the boundaries of $K_1$, $K_{22}$, and $K_{23}$.

\begin{proof}[Proof of {\em Lemma \ref{general lemma}} for non-convex $\hat{\Omega}_1$]

Let
\[ 
    \mathcal{N}(\bx)=1-\mathcal{N}_1(\bx)-\mathcal{N}_{22}(\bx).
\] 
Note that
\[
\hat{\chi}(\bx)-\mathcal{N}(\bx)=\left(\hat{\chi}_1-\mathcal{N}_1\right) + \left(\hat{\chi}_{22}-\mathcal{N}_{22}\right) + \hat{\chi}_{23}.
\]
It follows from the triangle inequality, \eqref{chi 1}, and \eqref{4.15} that
\begin{eqnarray}\nonumber
\|\hat{\chi}-\mathcal{N}\|_{L^p(\Omega)}&\le & \|\hat{\chi}_1-\mathcal{N}_1\|_{L^p\left(\Omega\right)}+\|\hat{\chi}_{22}-\mathcal{N}_{22}\|_{L^p\left(\Omega\right)}+\|\hat{\chi}_{23}\|_{L^p\left(\Omega\right)}\\ [2mm] \label{4.16}
&\le & \big| \hat{\Omega}_{\varepsilon, 1} \big|^{1/p}+\big| \hat{\Omega}_{\varepsilon, 22} \big|^{1/p}+|K_{23}|^{1/p} \leq C\, \varepsilon. 
\end{eqnarray}
This completes the approximation.
\end{proof}

\begin{remark}
The $\varepsilon$-width region $K_{23}$ is a subdomain of $\hat{\Omega}_1$, and its boundary contains the hyperplanes $\mathbf{a}\cdot\mathbf{x}-b=0$ (blue line in Figure \ref{convex decompistion}) and $\mathbf{a}\cdot\mathbf{x}-b-\varepsilon=0$ (red line in Figure \ref{convex decompistion2}) which are not part of the interface $\hat{\Gamma}$. For any $\bx\in K_{23}$, it is easy to see that 
\[
\hat{\chi}(\bx)-\mathcal{N}(\bx)= \mathcal{N}_1(\bx) +\mathcal{N}_{22} (\bx) =\left(1-\frac{1}{\varepsilon}(\mathbf{a}\cdot\mathbf{x}-b)\right) +\left(1-\frac{1}{\varepsilon}(-\mathbf{a}\cdot\mathbf{x}+b+\varepsilon)\right) =0,
\]
which, together with (\ref{4.16}), implies that the ReLU NN function $\mathcal{N}(\bx)$ approximates the discontinuous step function $\hat{\chi}(\bx)$ without overshooting. Moreover, it is clear from the construction that $\mathcal{N}(\bx)$ has no oscillation. No overshooting and no oscillation remain true for the ReLU NN approximation $\mathcal{N}(\bx)$ constructed in Subsubsection \ref{convex hull subsubsection}. 
\end{remark}

\section{Examples}\label{s:example}

This section validates our theoretical findings with several examples in $d\ge 2$ dimensions. The first three examples demonstrate Theorem \ref{general theorem} for convex $\hat{\Omega}_1$, with the third example extending to the case of $d=10000$. The final example illustrates a non-convex case using the twp decomposition procedures outlined in Subsection \ref{non-convex}.

\subsection{A two-dimensional circular interface}\label{convex 2d example}

\begin{figure}[htbp]
\centering
\subfigure[The piecewise constant function $\chi(\mathbf{x})$\label{2d example function}]{
\begin{minipage}[t]{0.4\linewidth}
\centering
\includegraphics[width=1.8in]{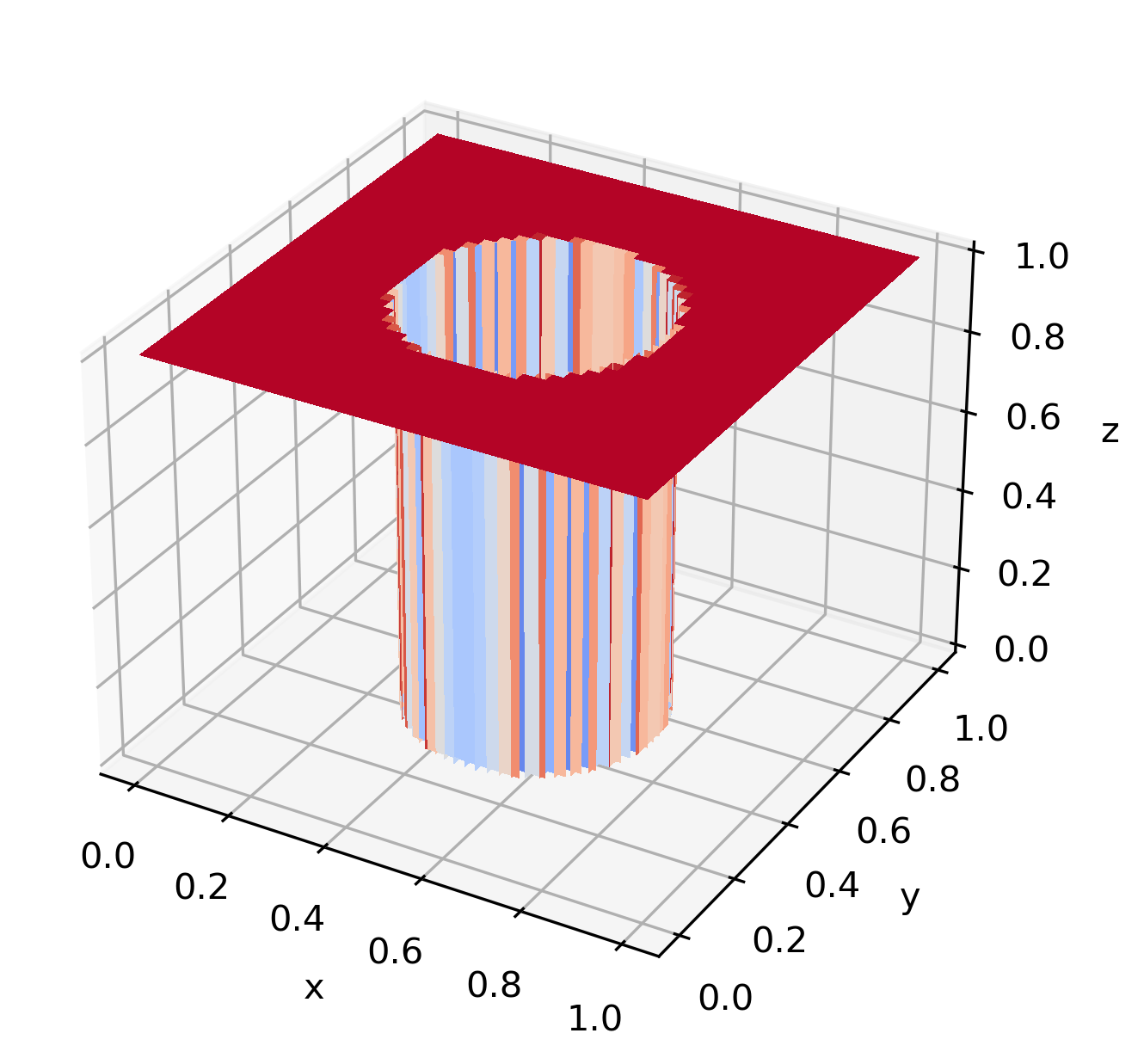}
\end{minipage}%
}%
\hspace{0.2in}
\subfigure[The circular interface\label{2d interface}]{
\begin{minipage}[t]{0.4\linewidth}
\centering
\includegraphics[width=1.8in]{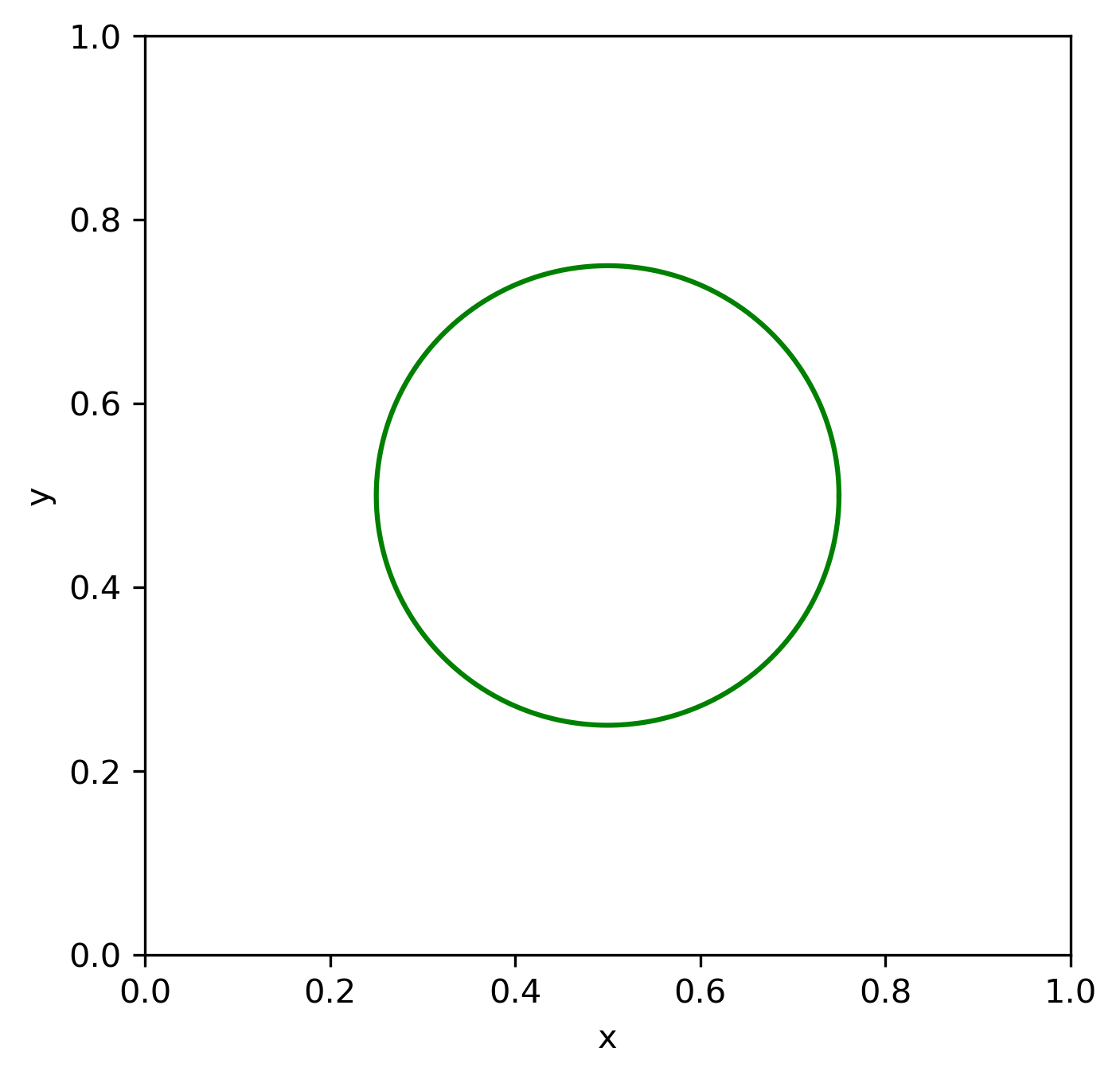}
\end{minipage}%
}%
\\
\subfigure[An approximation of the interface by $n=6$ line segments\label{2d_approx_n6}]{
\begin{minipage}[t]{0.4\linewidth}
\centering
\includegraphics[width=1.8in]{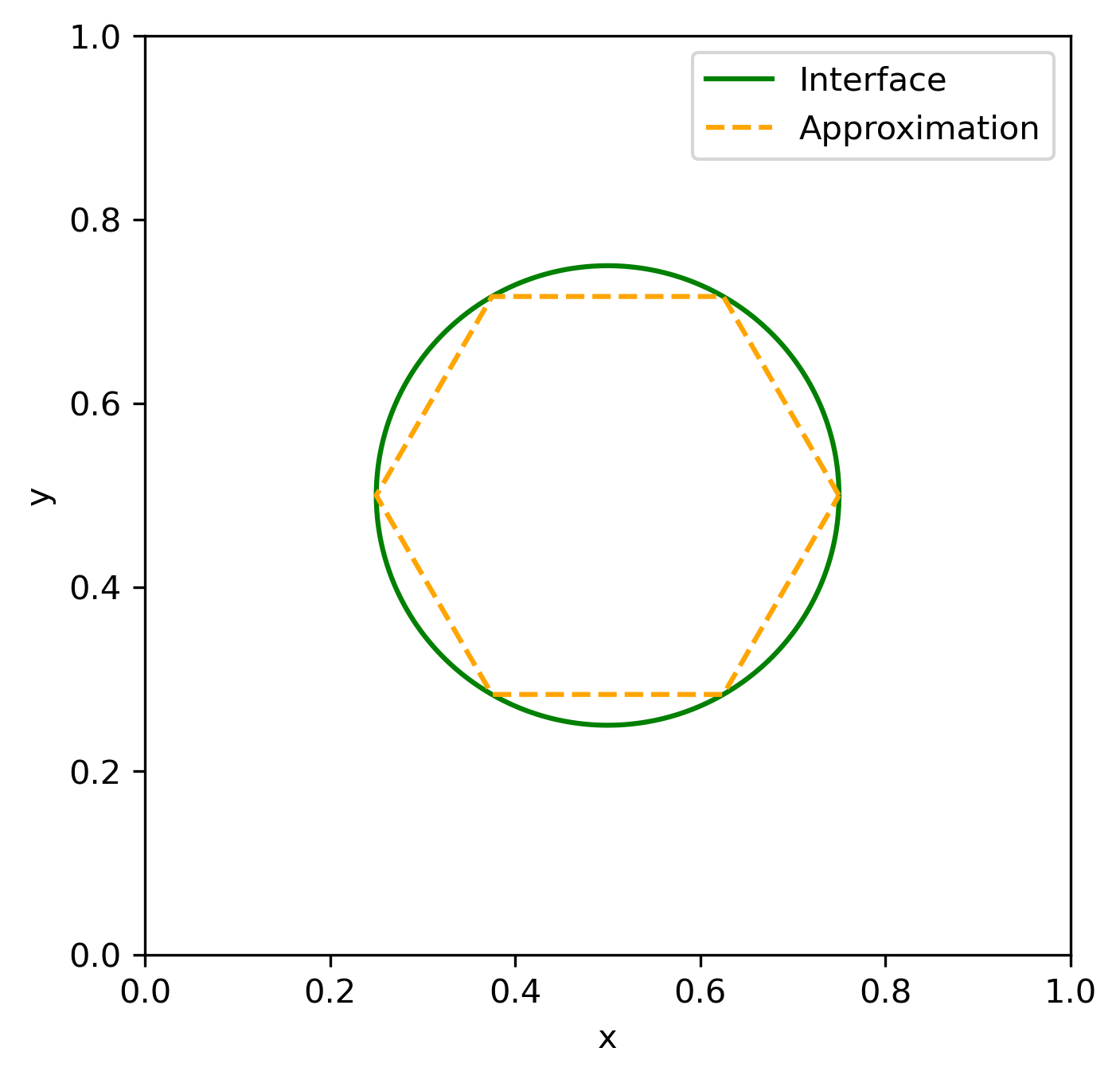}
\end{minipage}%
}%
\hspace{0.2in}
\subfigure[An approximation of the interface by $n=50$ line segments\label{2d_approx_n50}]{
\begin{minipage}[t]{0.4\linewidth}
\centering
\includegraphics[width=1.8in]{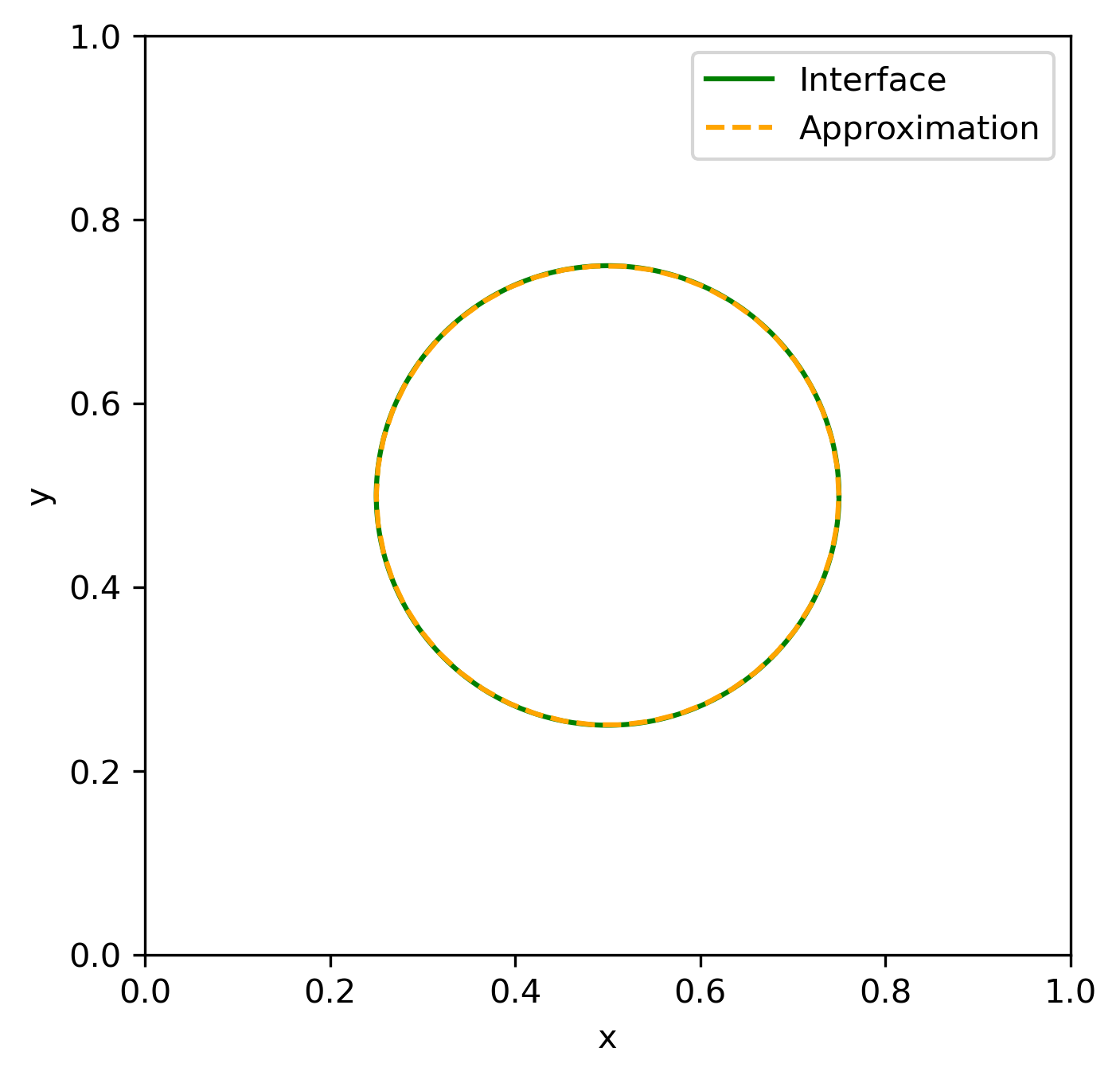}
\end{minipage}%
}%
\\
\subfigure[An approximation of $\chi(\mathbf{x})$ by the $2$--$6$--1--1 ReLU NN function in \eqref{convex analytic} with $\varepsilon=1/25$\label{2d approx graph n6}]{
\begin{minipage}[t]{0.4\linewidth}
\centering
\includegraphics[width=1.8in]{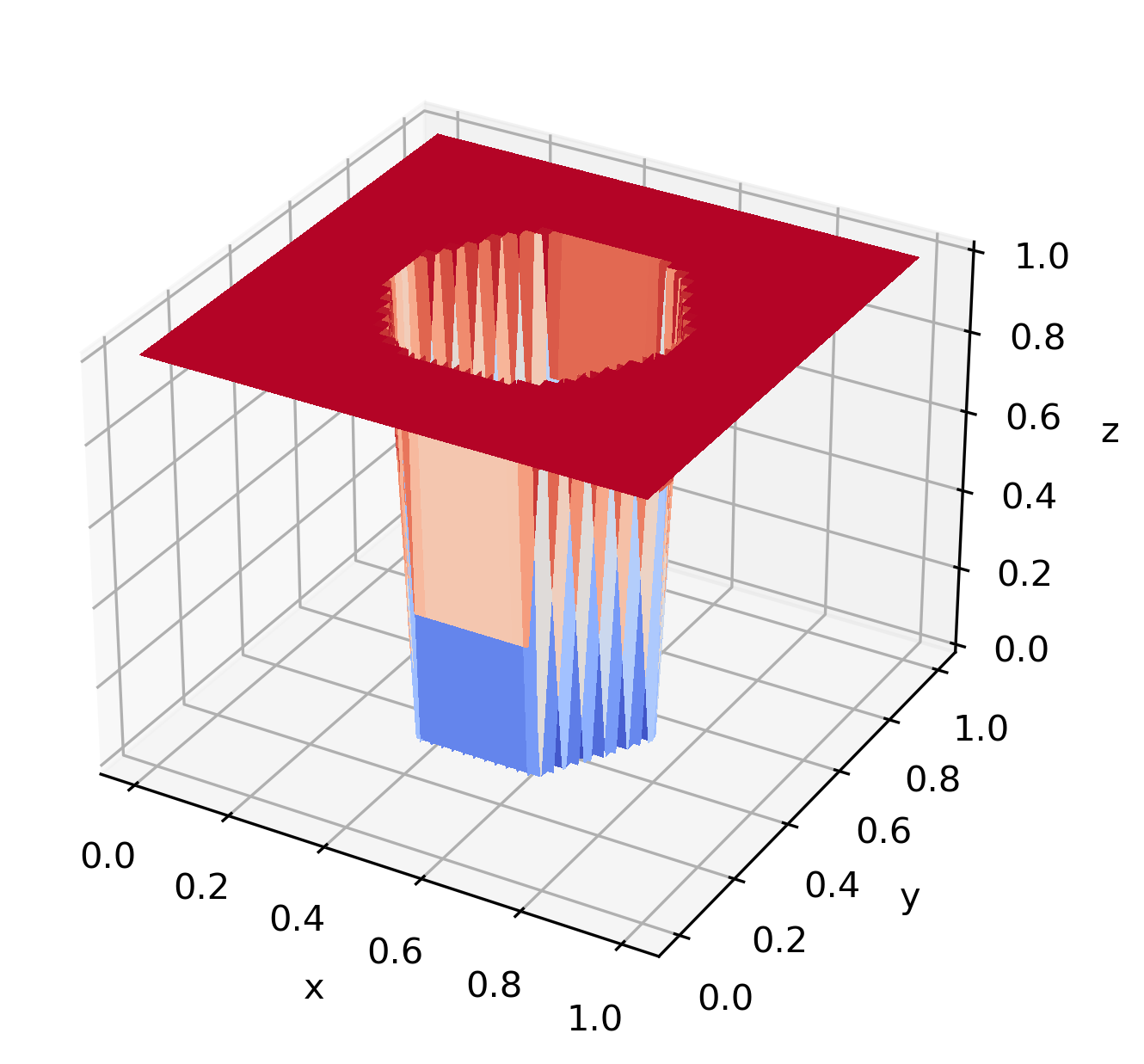}
\end{minipage}%
}%
\hspace{0.2in}
\subfigure[An approximation of $\chi(\mathbf{x})$ by the $2$--$50$--1--1 ReLU NN function in \eqref{convex analytic} with $\varepsilon=1/2000$\label{2d approx graph n50}]{
\begin{minipage}[t]{0.4\linewidth}
\centering
\includegraphics[width=1.8in]{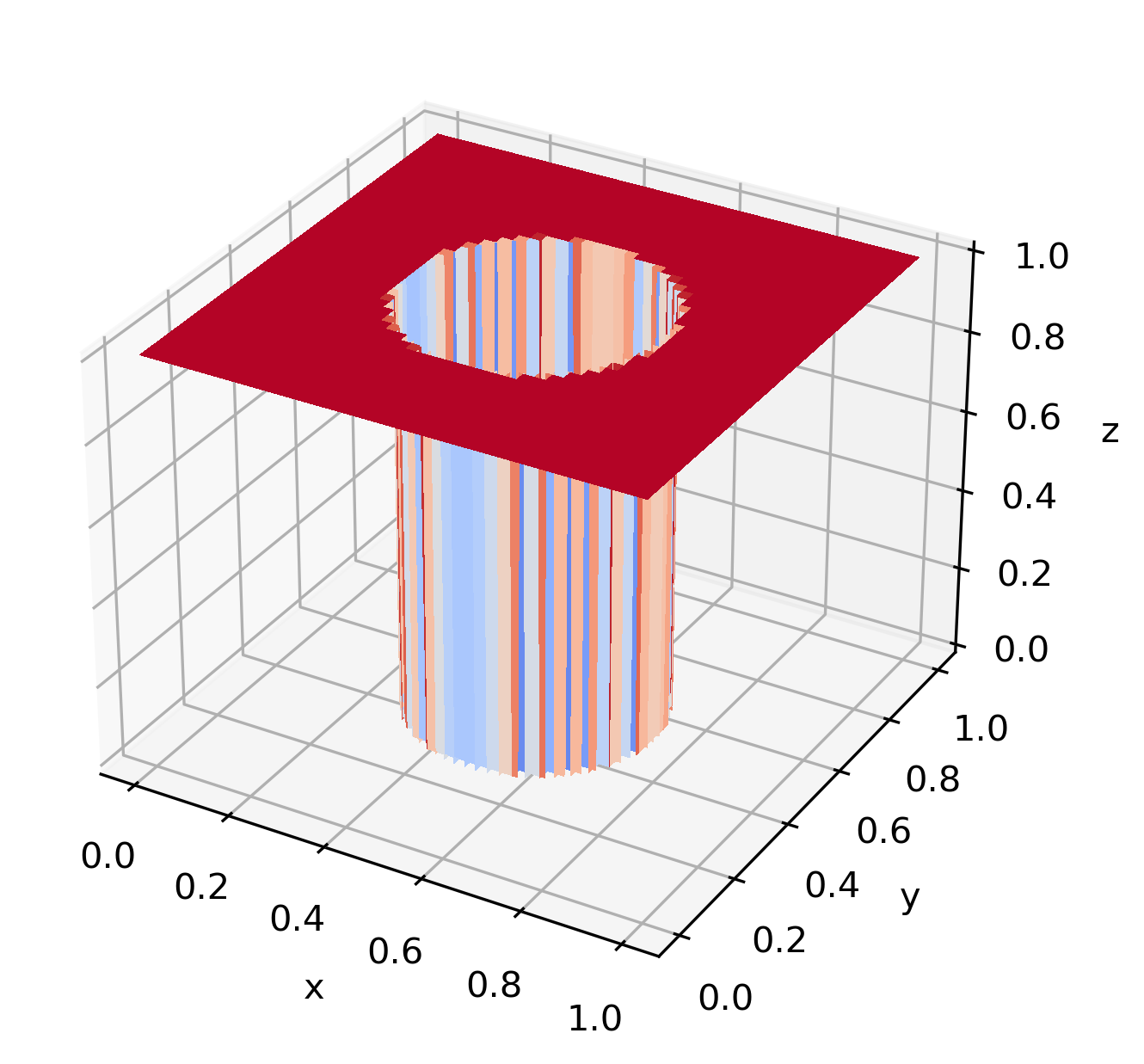}
\end{minipage}%
}%
\\
\subfigure[The breaking hyperplanes of the approximation in Figure \ref{2d approx graph n6}\label{2d_breaking_n6}]{
\begin{minipage}[t]{0.4\linewidth}
\centering
\includegraphics[width=1.8in]{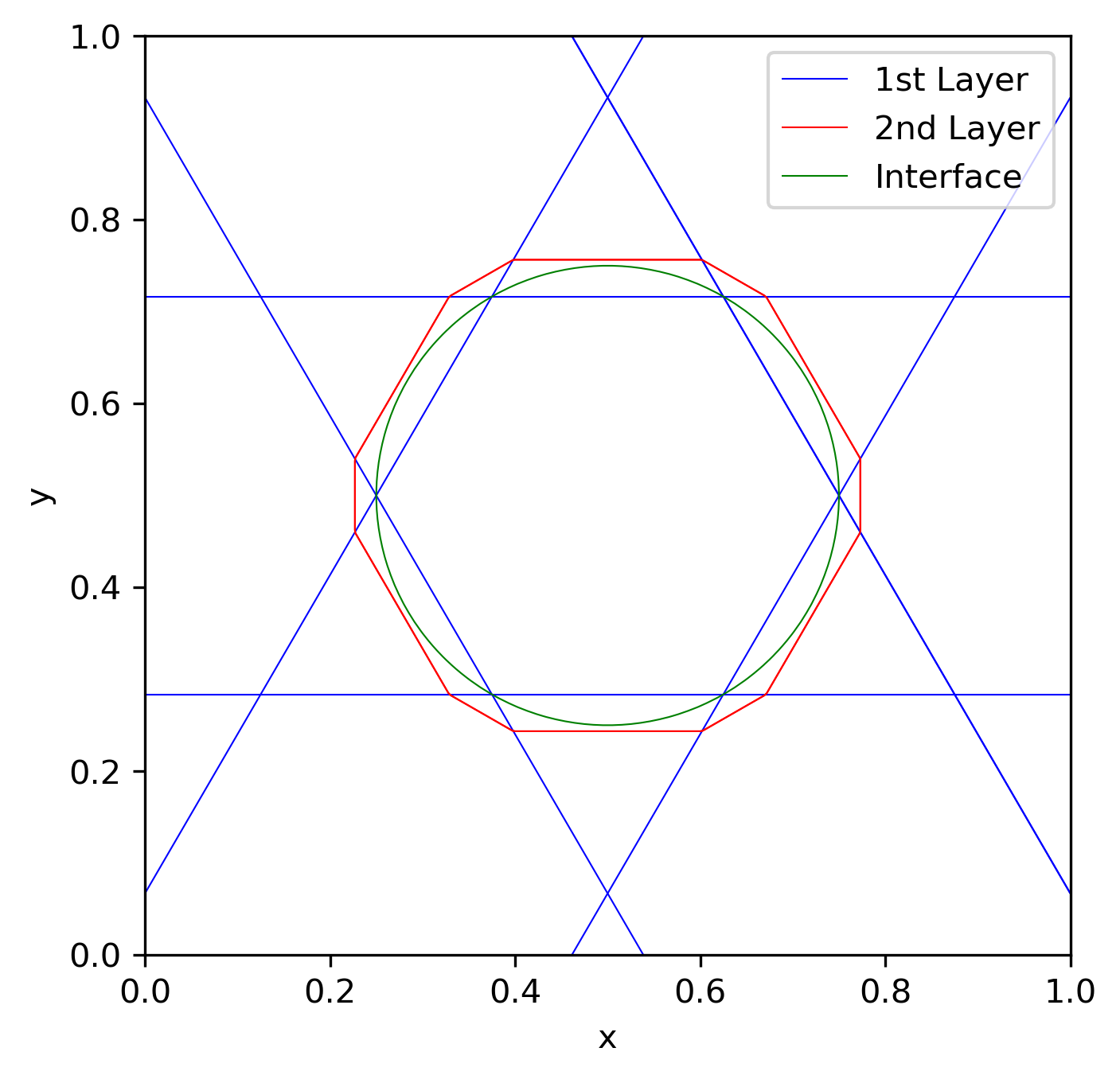}
\end{minipage}%
}%
\hspace{0.2in}
\subfigure[The breaking hyperplanes of the approximation in Figure \ref{2d approx graph n50}\label{2d_breaking_n50}]{
\begin{minipage}[t]{0.4\linewidth}
\centering
\includegraphics[width=1.8in]{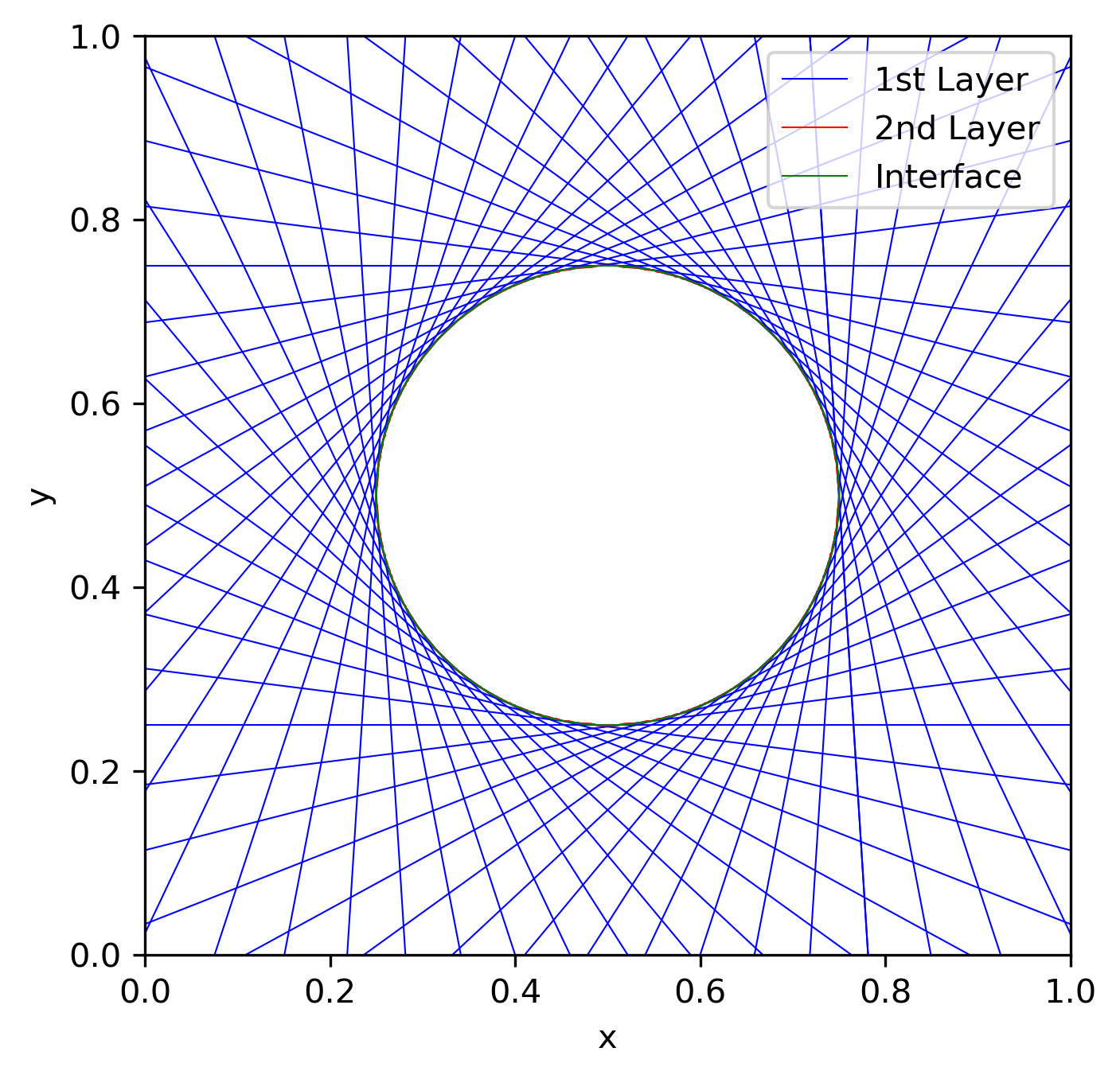}
\end{minipage}%
}%
\caption{A convex example to illustrate Theorem \ref{general theorem} for the case $d=2$}
\label{2d example}
\end{figure}

Let $\Omega=(0,1)^2$, 
\[
\Omega_1=\{(x,y)\in \Omega:(x-0.5)^2+(y-0.5)^2<0.25^2\}, \quad\mbox{and}\quad \Omega_2=\Omega\setminus\Omega_1.
\]
The piecewise constant function $\chi(\mathbf{x})$ is shown in Figure \ref{2d example function}. The interface $\Gamma$ is a circle centered at $(0.5,0.5)$ with a radius of $0.25$ (see Figure \ref{2d interface}):
\[
\Gamma=\{(x,y)\in \Omega:(x-0.5)^2+(y-0.5)^2=0.25^2\}.
\]
Consider approximations of the interface $\Gamma$ by $n=6$ and $50$ line segments (see Figures \ref{2d_approx_n6} and \ref{2d_approx_n50}), respectively. The $2$--$6$--1--1 and $2$--$50$--1--1 ReLU NN approximations given in \eqref{convex analytic} with $\varepsilon=1/25$ and $1/2000$ are shown in Figures \ref{2d approx graph n6} and \ref{2d approx graph n50}, respectively.  Figures \ref{2d_breaking_n6} and \ref{2d_breaking_n50} illustrate the breaking lines of the first and second layers, with the distances between them equal to $\varepsilon$.

\subsection{A three-dimensional spherical interface}\label{convex 3d example}

\begin{figure}[htbp]
\centering
\subfigure[The piecewise constant function $\chi(\mathbf{x})$ on $z=0.205$\label{3d example function}]{
\begin{minipage}[t]{0.4\linewidth}
\centering
\includegraphics[width=1.8in]{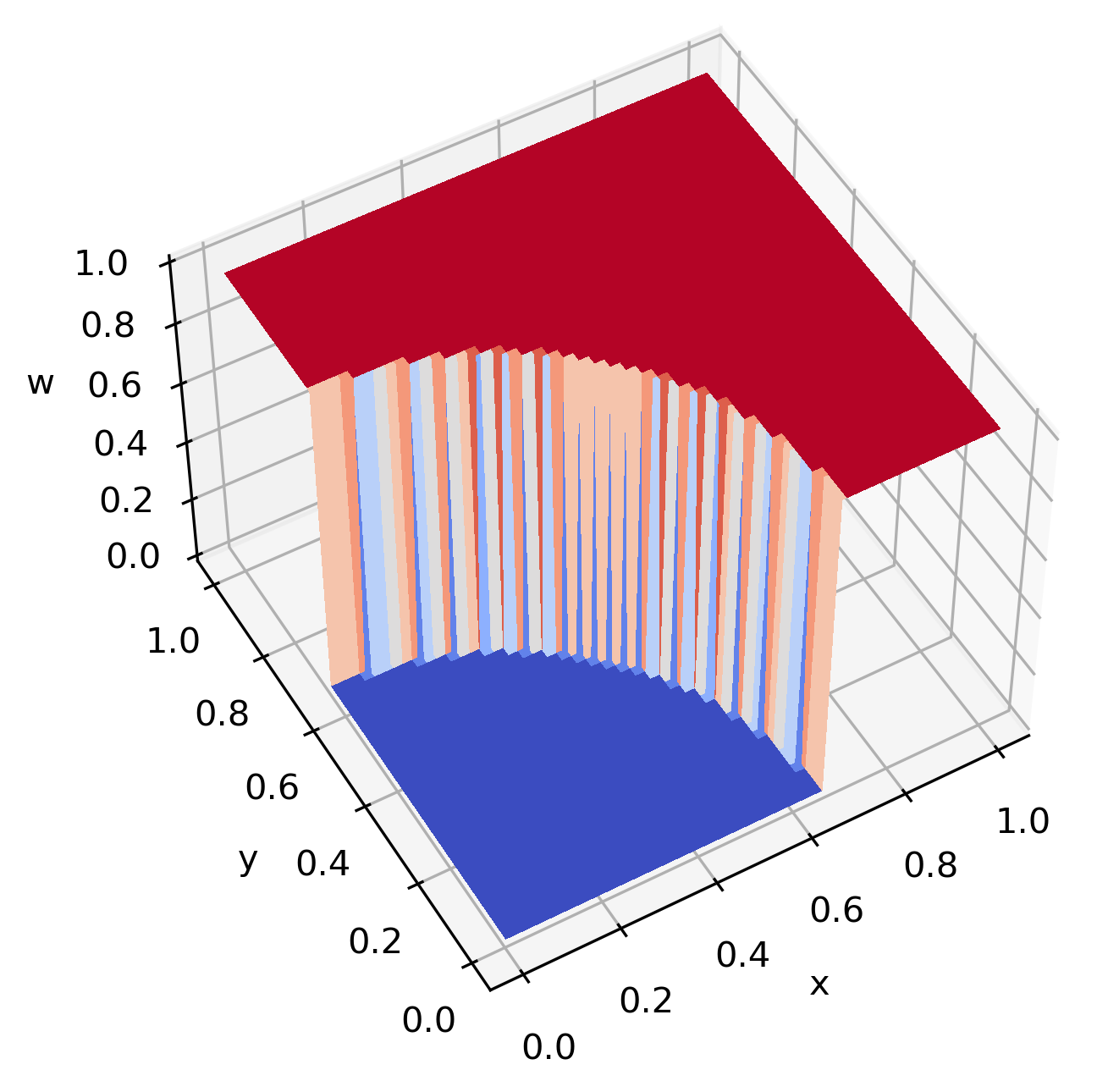}
\end{minipage}%
}%
\hspace{0.2in}
\subfigure[The spherical interface\label{3d interface}]{
\begin{minipage}[t]{0.4\linewidth}
\centering
\includegraphics[width=1.8in]{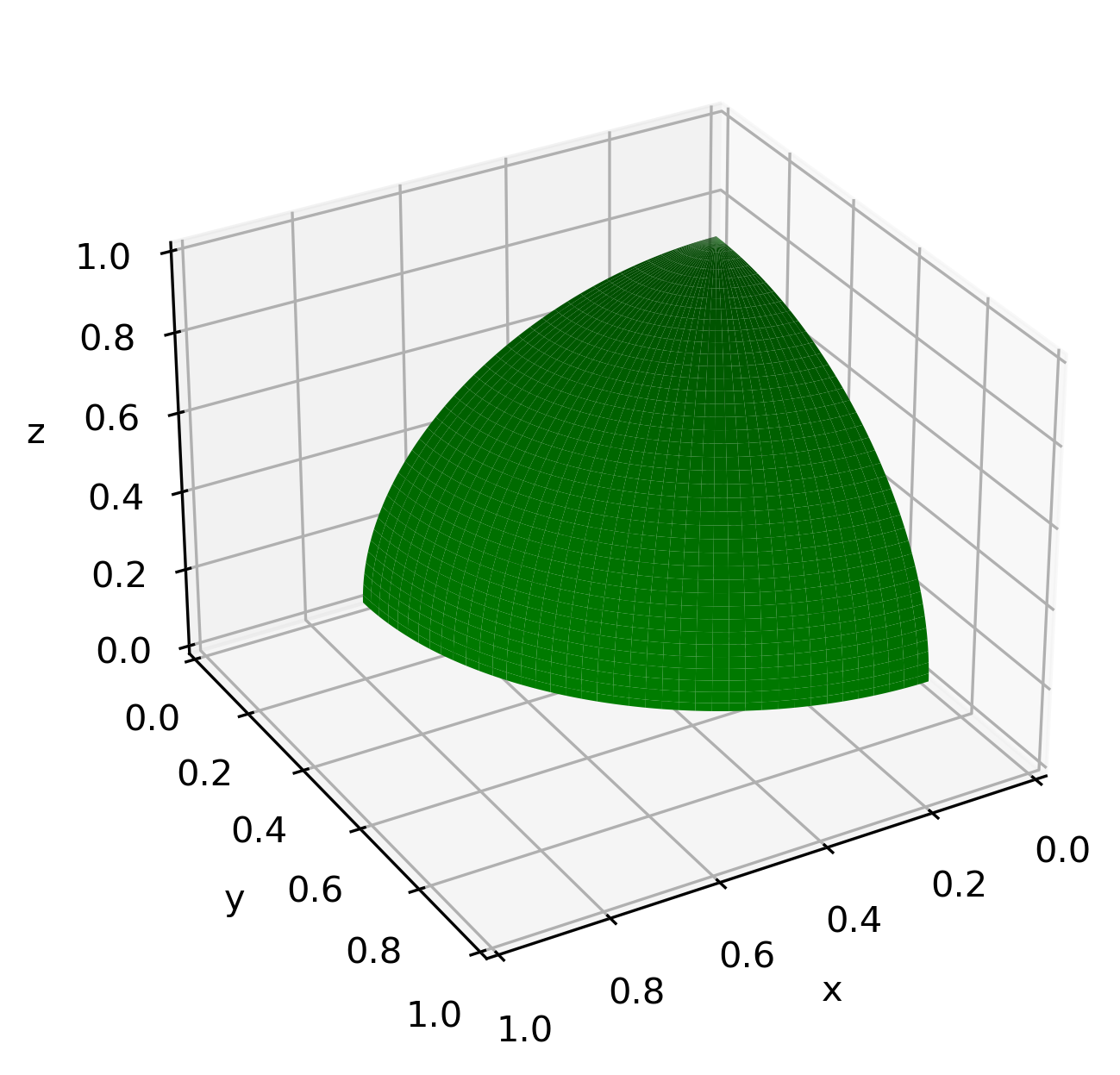}
\end{minipage}%
}%
\\
\subfigure[An approximation of the interface by $n=9$ plane segments\label{3d_approx_n9}]{
\begin{minipage}[t]{0.4\linewidth}
\centering
\includegraphics[width=1.8in]{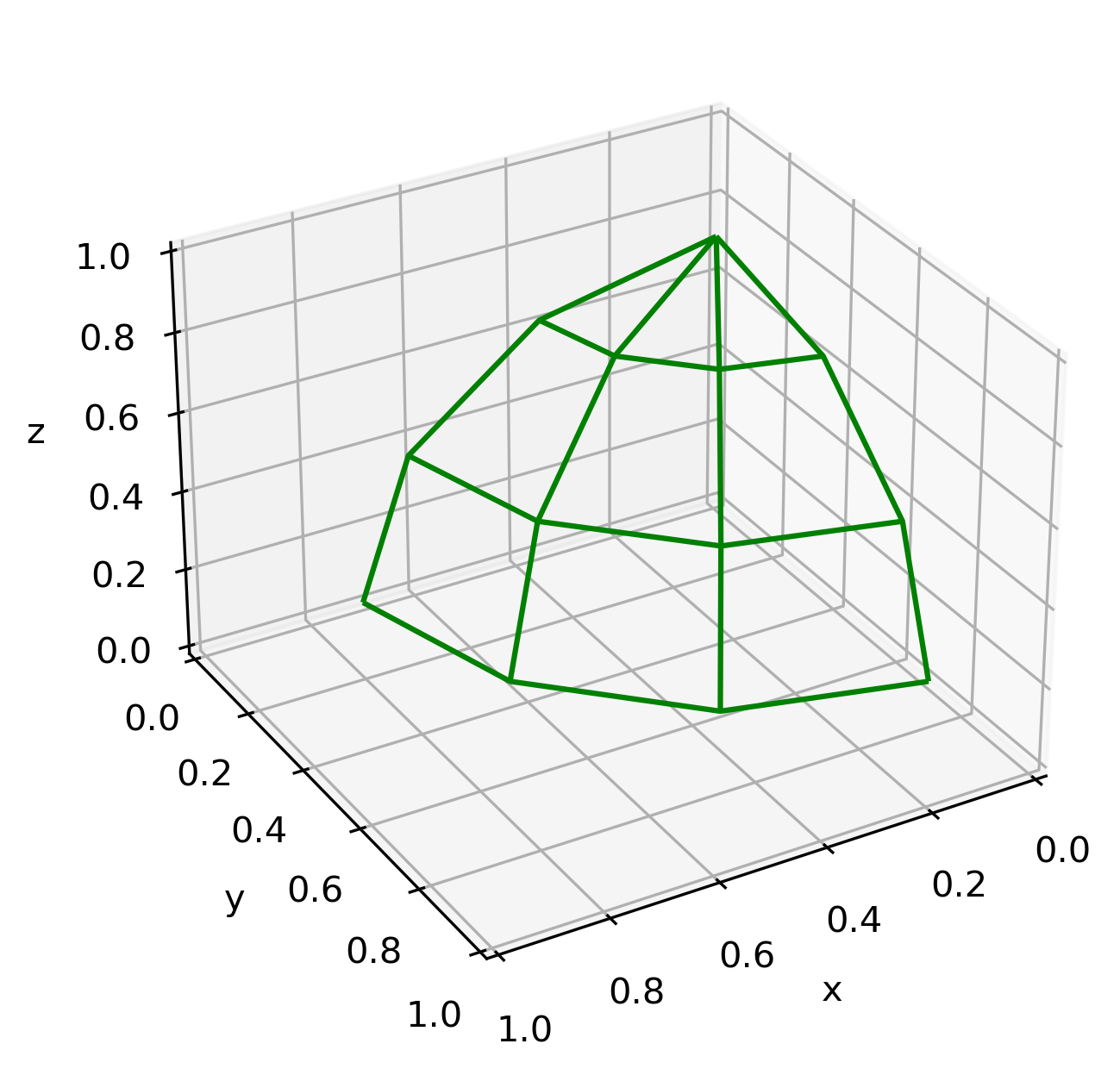}
\end{minipage}%
}%
\hspace{0.2in}
\subfigure[An approximation of the interface by $n=100$ plane segments\label{3d_approx_n100}]{
\begin{minipage}[t]{0.4\linewidth}
\centering
\includegraphics[width=1.8in]{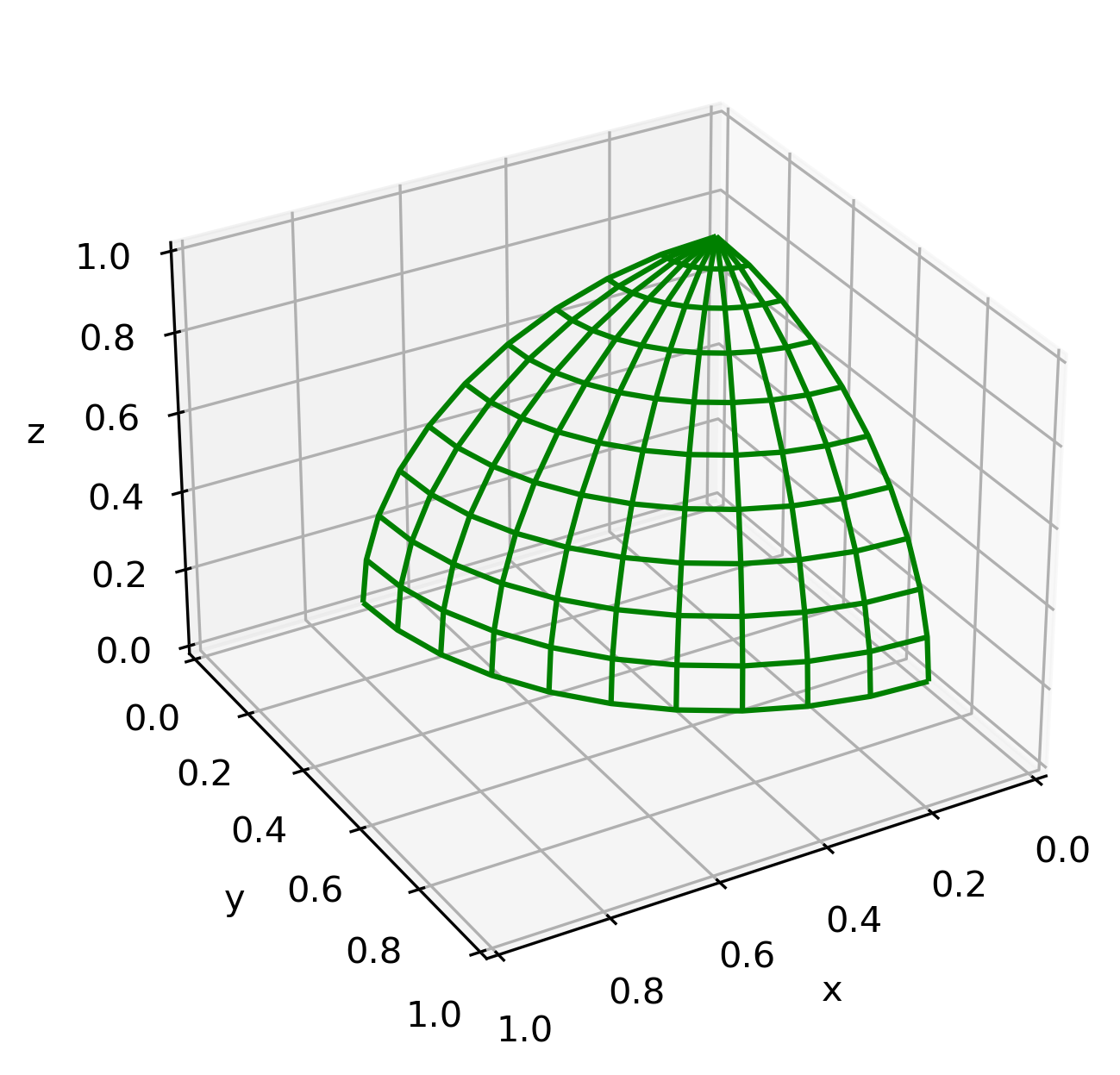}
\end{minipage}%
}%
\\
\subfigure[An approximation of $\chi(\mathbf{x})$ by the $3$--$9$--1--1 ReLU NN function in \eqref{convex analytic} with $\varepsilon=1/15$ on $z=0.205$\label{3d approx graph n9}]{
\begin{minipage}[t]{0.4\linewidth}
\centering
\includegraphics[width=1.8in]{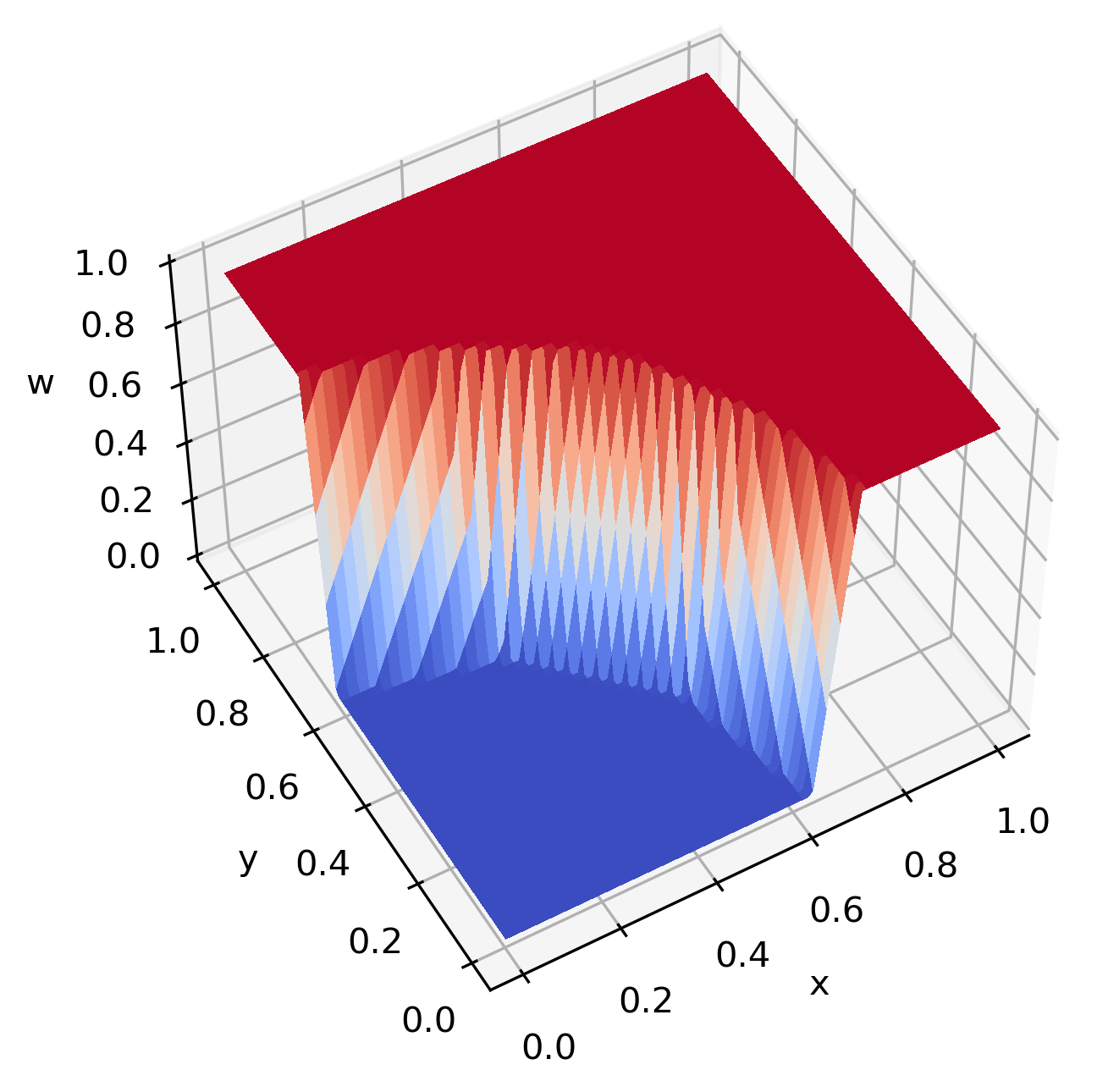}
\end{minipage}%
}%
\hspace{0.2in}
\subfigure[An approximation of $\chi(\mathbf{x})$ by the $3$--$100$--1--1 ReLU NN function in \eqref{convex analytic} with $\varepsilon=1/100$ on $z=0.205$\label{3d approx graph n100}]{
\begin{minipage}[t]{0.4\linewidth}
\centering
\includegraphics[width=1.8in]{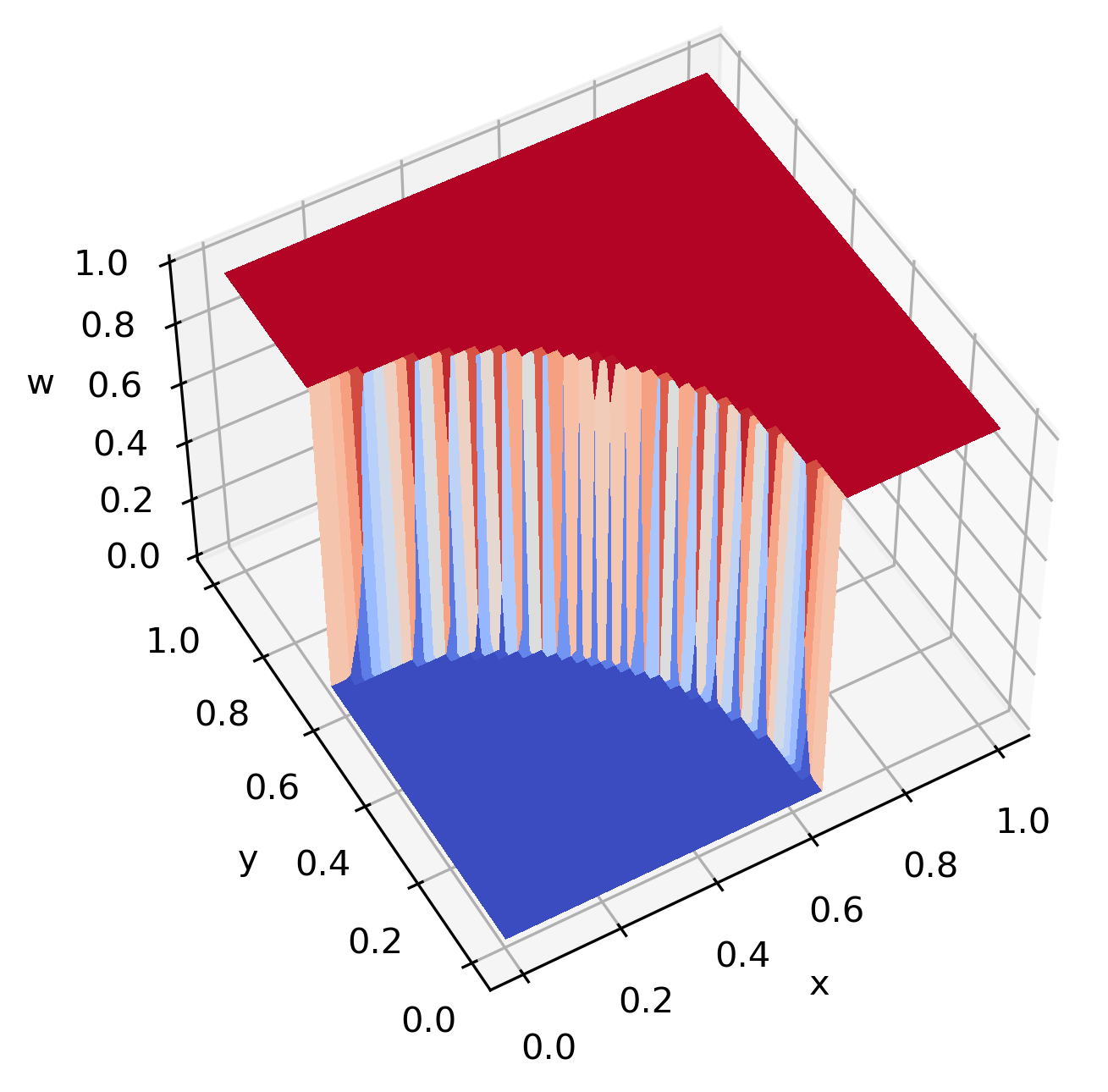}
\end{minipage}%
}%
\\
\subfigure[The breaking hyperplanes of the approximation in Figure \ref{3d approx graph n9} on $z=0.205$\label{3d_breaking_n9}]{
\begin{minipage}[t]{0.4\linewidth}
\centering
\includegraphics[width=1.8in]{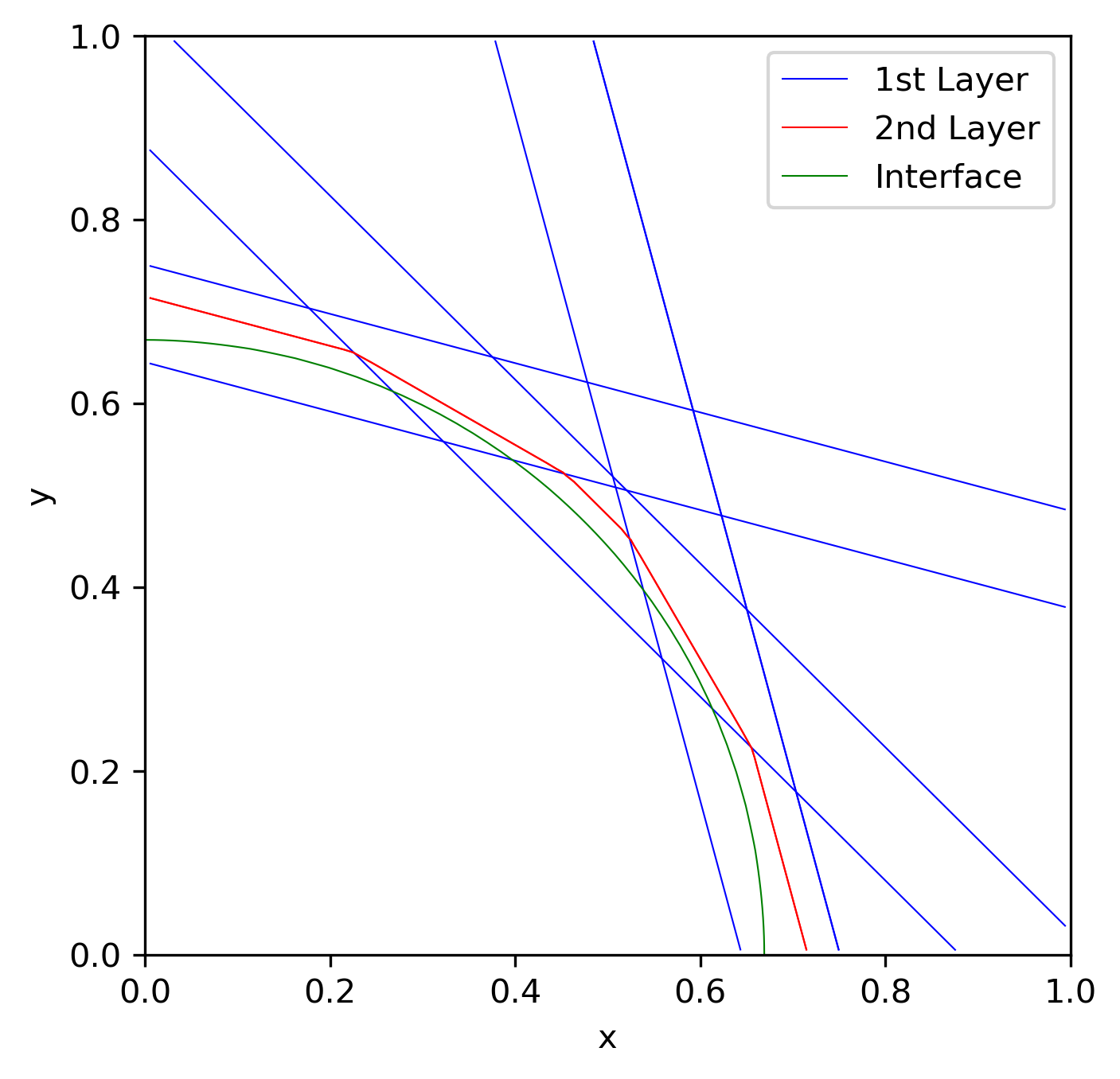}
\end{minipage}%
}%
\hspace{0.2in}
\subfigure[The breaking hyperplanes of the approximation in Figure \ref{3d approx graph n100} on $z=0.205$\label{3d_breaking_n100}]{
\begin{minipage}[t]{0.4\linewidth}
\centering
\includegraphics[width=1.8in]{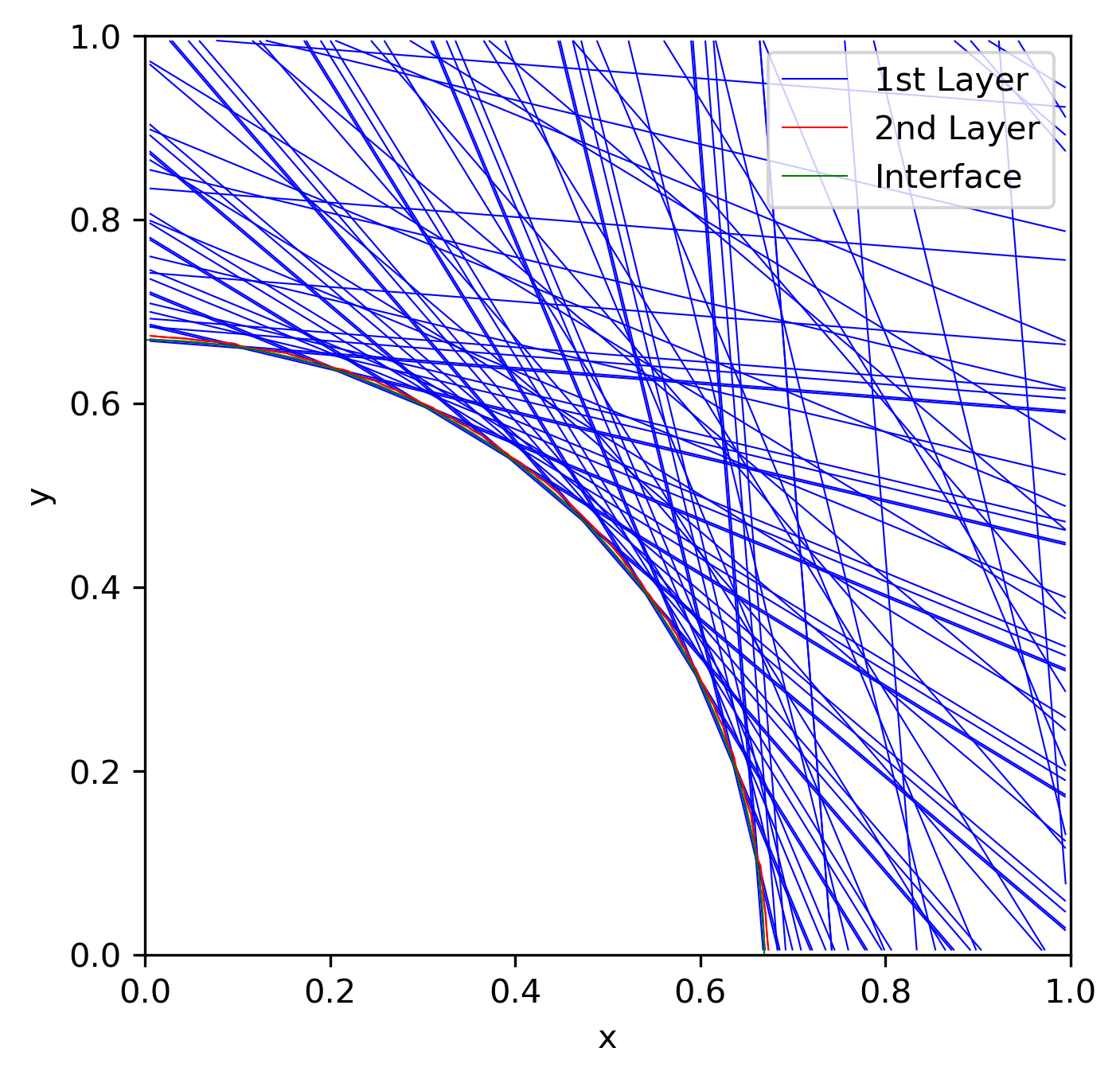}
\end{minipage}%
}%
\caption{A convex example to illustrate Theorem \ref{general theorem} for the case $d=3$}
\label{3d example}
\end{figure}

Let $\Omega=(0,1)^3$, 
\[
\Omega_1=\{(x,y,z)\in \Omega: z<\sqrt{0.7^2-x^2-y^2}\}, \quad\mbox{and}\quad \Omega_2=\Omega\setminus\Omega_1. 
\]
The intersection between the piecewise constant function $\chi(x,y,z)$ and the hyperplane $z=0.205$ is shown in Figure \ref{3d example function}. The interface $\Gamma$ is part of a sphere centered at $(0,0,0)$ with a radius of $0.7$ (see Figure \ref{3d interface}):
\[
\Gamma=\{(x,y,z)\in \Omega:x^2+y^2+z^2=0.7^2\},
\]
which is approximated by $n=9$ and $100$ plane segments (see Figures \ref{3d_approx_n9} and \ref{3d_approx_n100}), respectively. The $3$--$9$--1--1 and $3$--$100$--1--1 ReLU NN approximations given in \eqref{convex analytic} with $\varepsilon=1/15$ and $1/100$ are depicted in Figures \ref{3d approx graph n9} and \ref{3d approx graph n100}, respectively. Figures \ref{3d_breaking_n9} and \ref{3d_breaking_n100} illustrate the first- and second-layer breaking hyperplanes on $z=0.205$.

\subsection{A $10000$-dimensional hypercube interface}\label{convex nd example}
\begin{figure}[htbp]
\centering
\subfigure[The piecewise constant function $\chi(\mathbf{x})$ on $x_i=0.255$ for $i= 3,\ldots, 10000$\label{nd example function}]{
\begin{minipage}[t]{0.4\linewidth}
\centering
\includegraphics[width=1.8in]{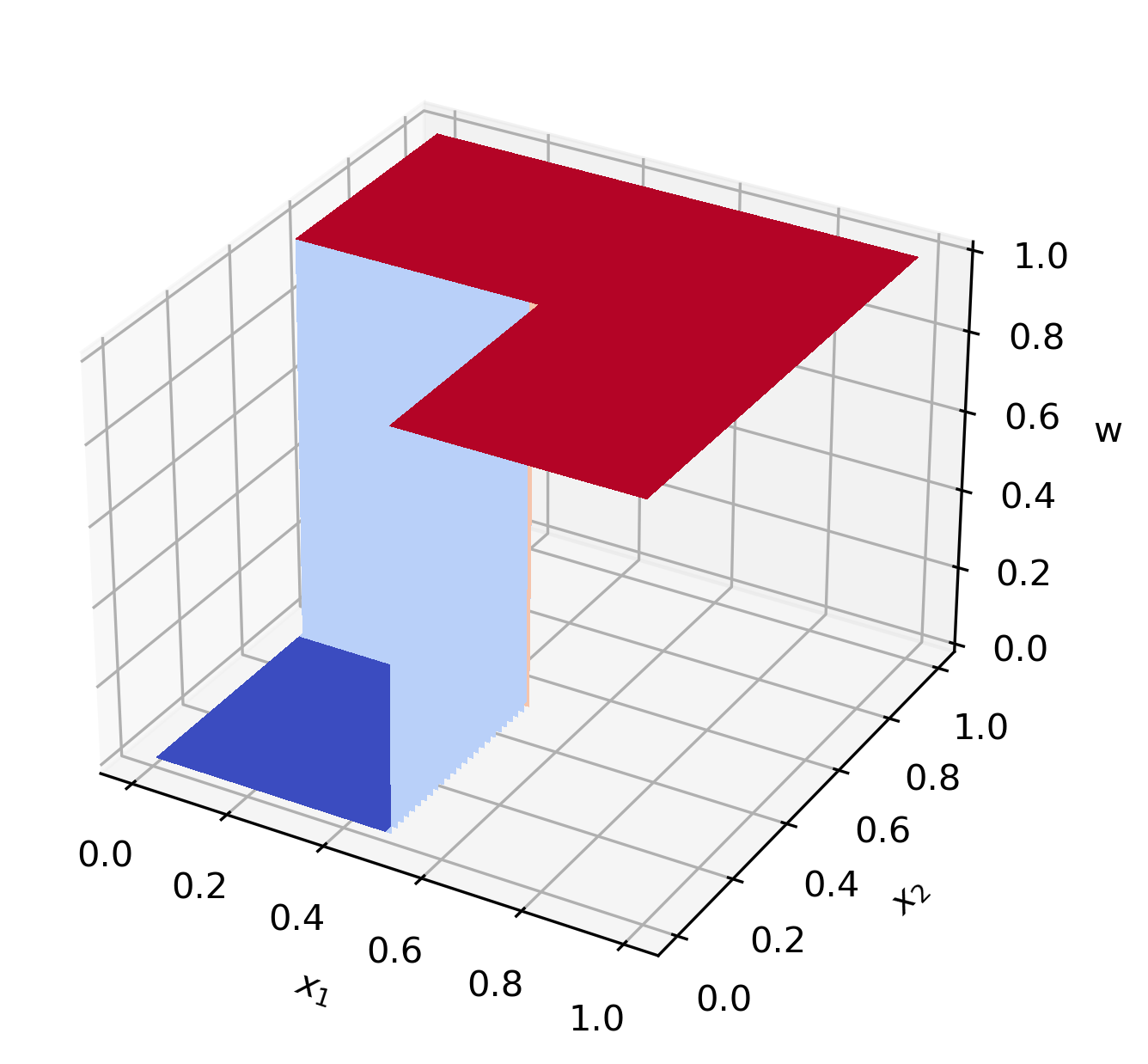}
\end{minipage}%
}%
\hspace{0.2in}
\subfigure[The interface on $x_i=0.255$ for $i=4,\ldots, 10000$\label{nd interface}]{
\begin{minipage}[t]{0.4\linewidth}
\centering
\includegraphics[width=1.8in]{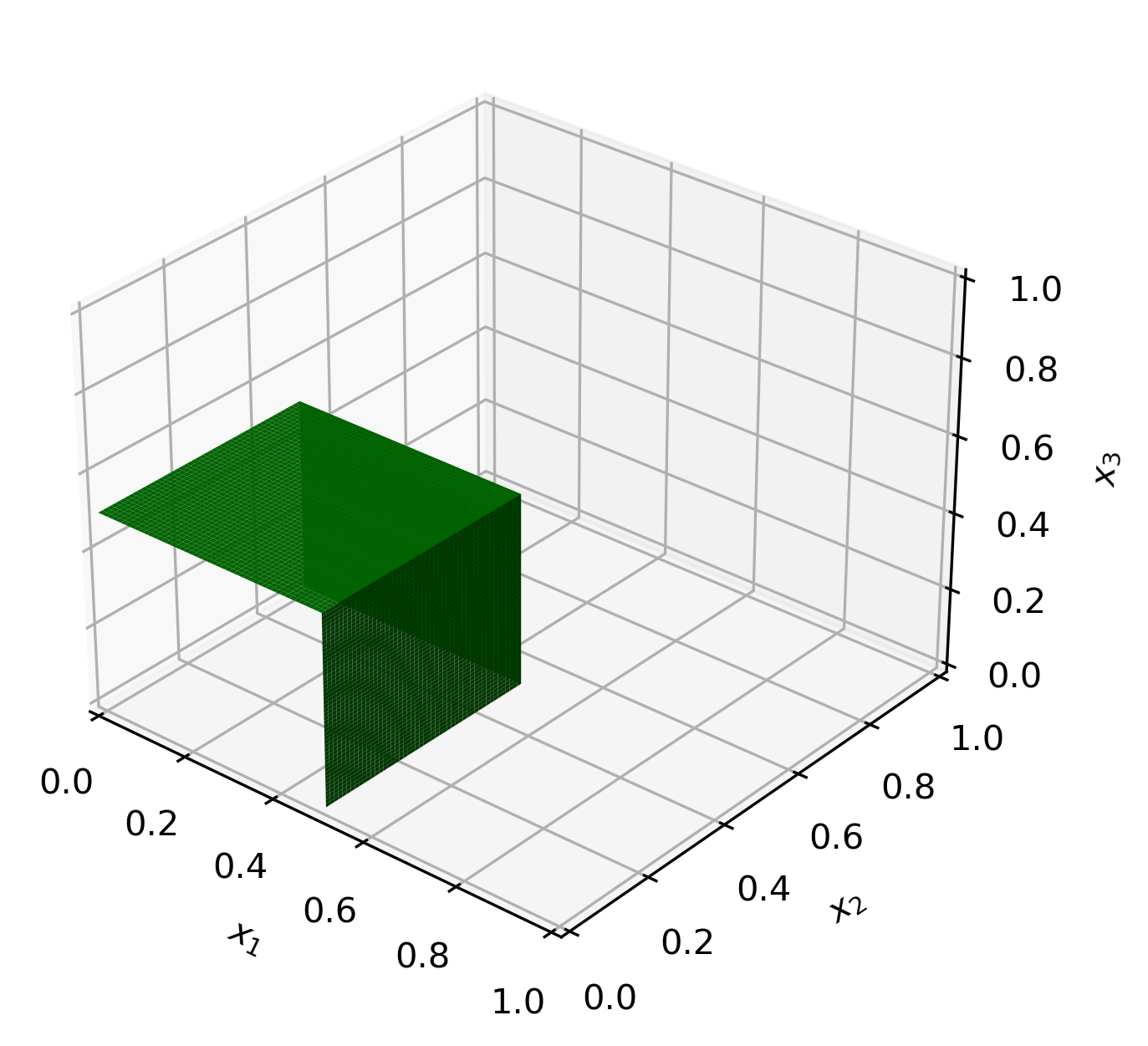}
\end{minipage}%
}%
\\
\subfigure[An approximation of $\chi(\mathbf{x})$ by the $10000$--$10000$--1--1 ReLU NN function in \eqref{convex analytic} with $\varepsilon=1/20$ on $x_i=0.255$ for $i=3,\ldots, 10000$\label{nd approx graph e20}]{
\begin{minipage}[t]{0.4\linewidth}
\centering
\includegraphics[width=1.8in]{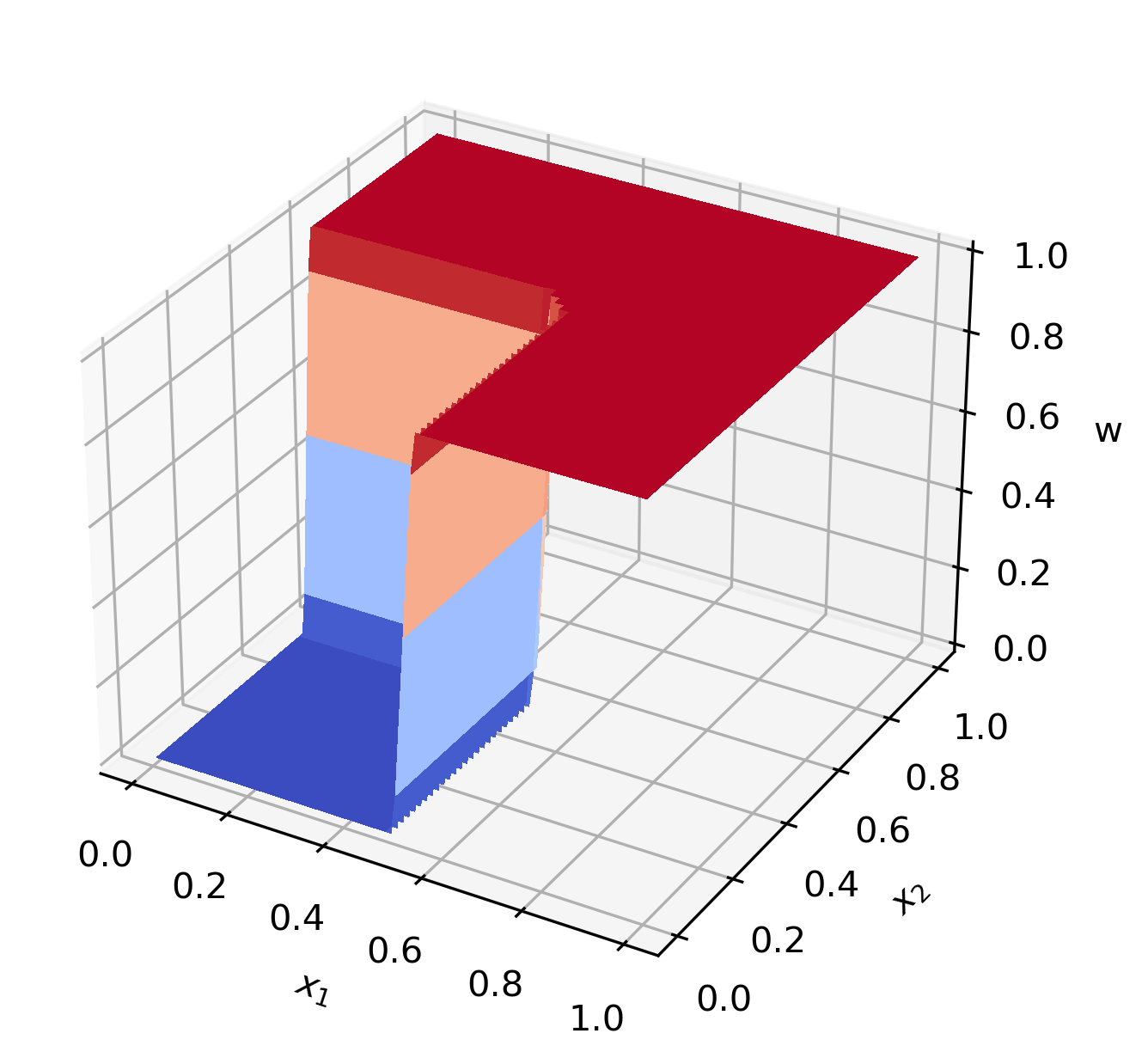}
\end{minipage}%
}%
\hspace{0.2in}
\subfigure[An approximation of $\chi(\mathbf{x})$ by the $10000$--$10000$--1--1 ReLU NN function in \eqref{convex analytic} with $\varepsilon=1/200$ on $x_i=0.255$ for $i=3,\ldots, 10000$\label{nd approx graph e200}]{
\begin{minipage}[t]{0.4\linewidth}
\centering
\includegraphics[width=1.8in]{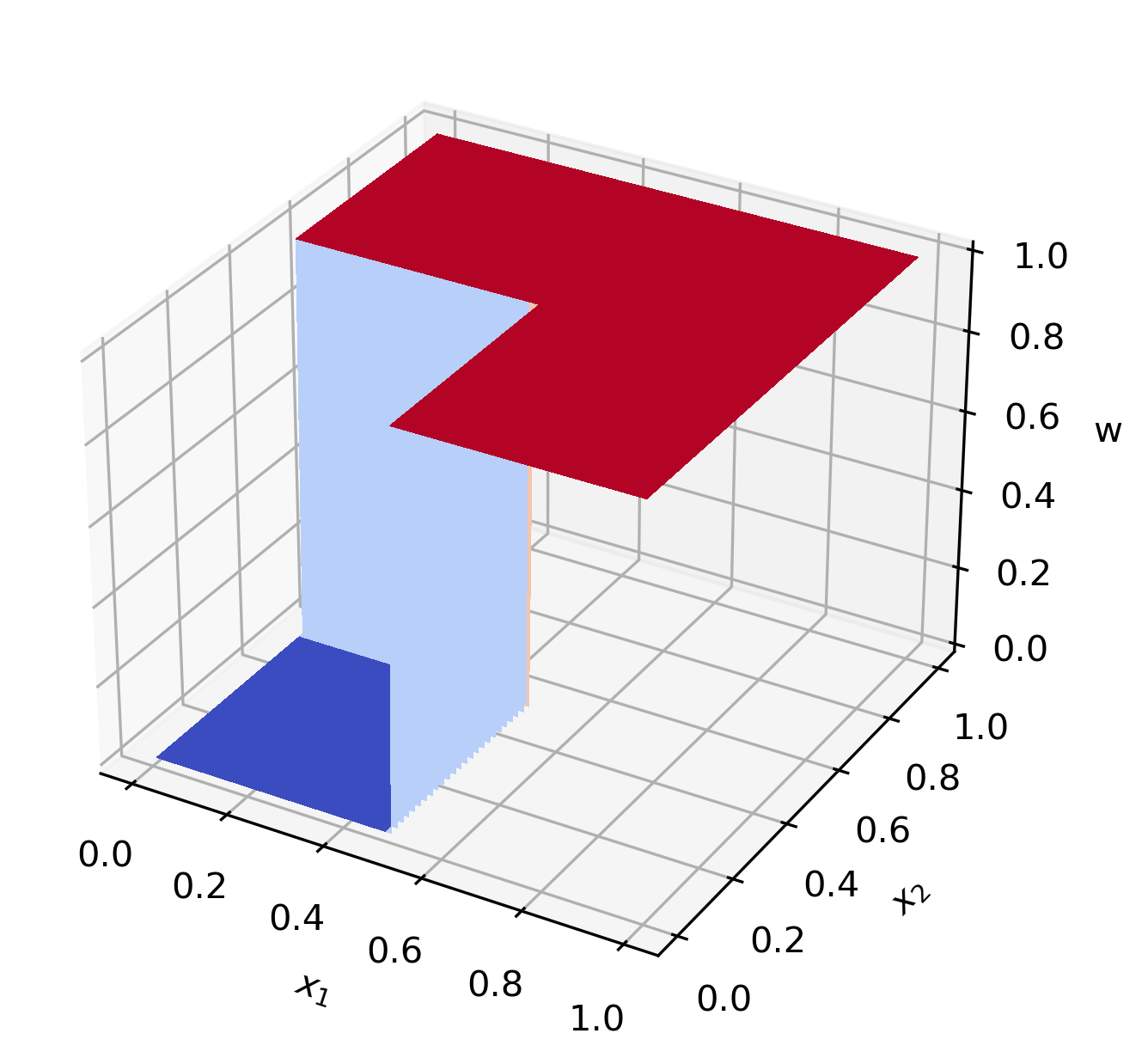}
\end{minipage}%
}%
\\
\subfigure[The breaking hyperplanes of the approximation in Figure \ref{nd approx graph e20} on $x_i=0.255$ for $i=3,\ldots, 10000$\label{nd_breaking_e20}]{
\begin{minipage}[t]{0.4\linewidth}
\centering
\includegraphics[width=1.8in]{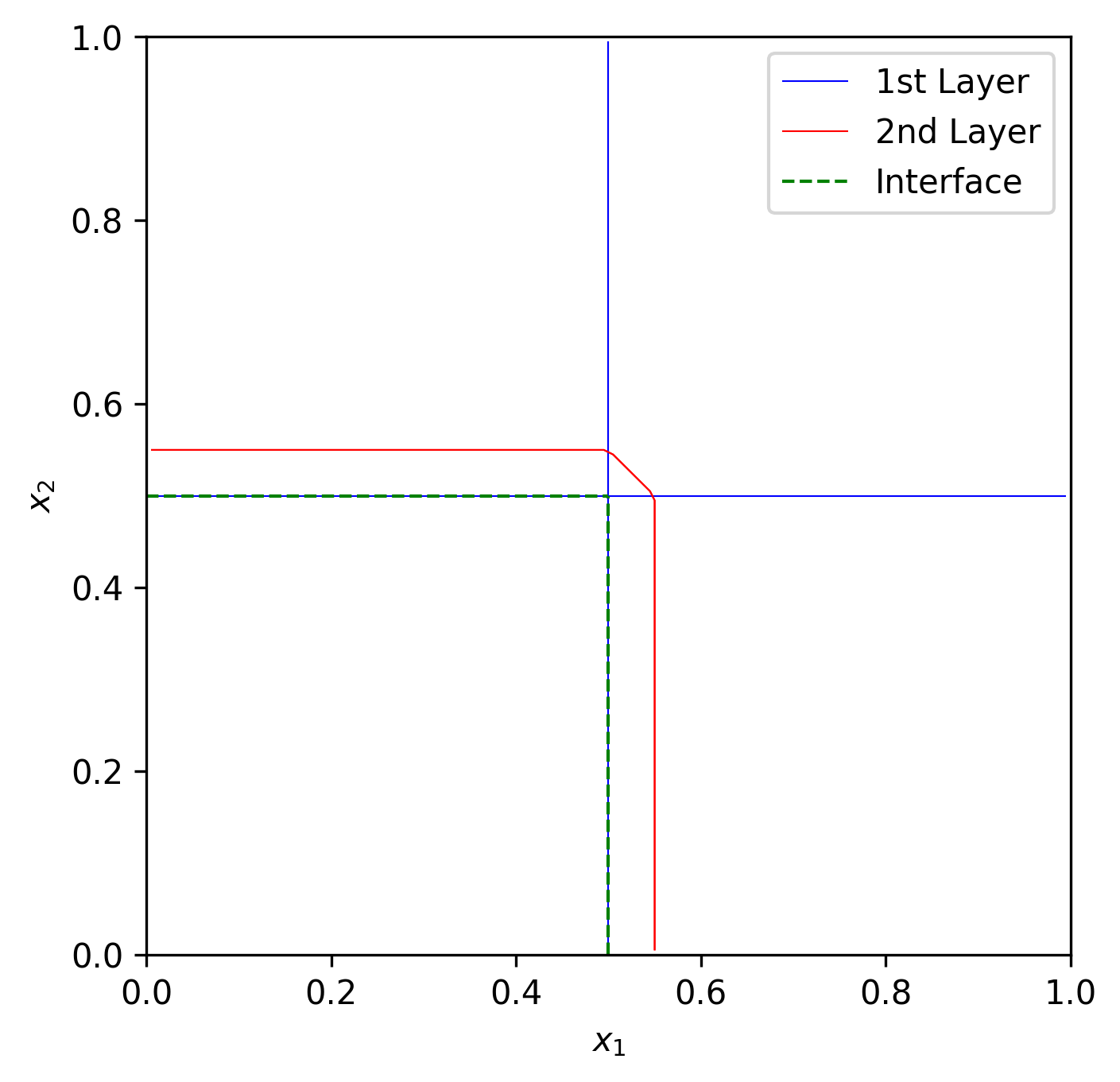}
\end{minipage}%
}%
\hspace{0.2in}
\subfigure[The breaking hyperplanes of the approximation in Figure \ref{nd approx graph e200} on $x_i=0.255$ for $i=3,\ldots, 10000$\label{nd_breaking_e200}]{
\begin{minipage}[t]{0.4\linewidth}
\centering
\includegraphics[width=1.8in]{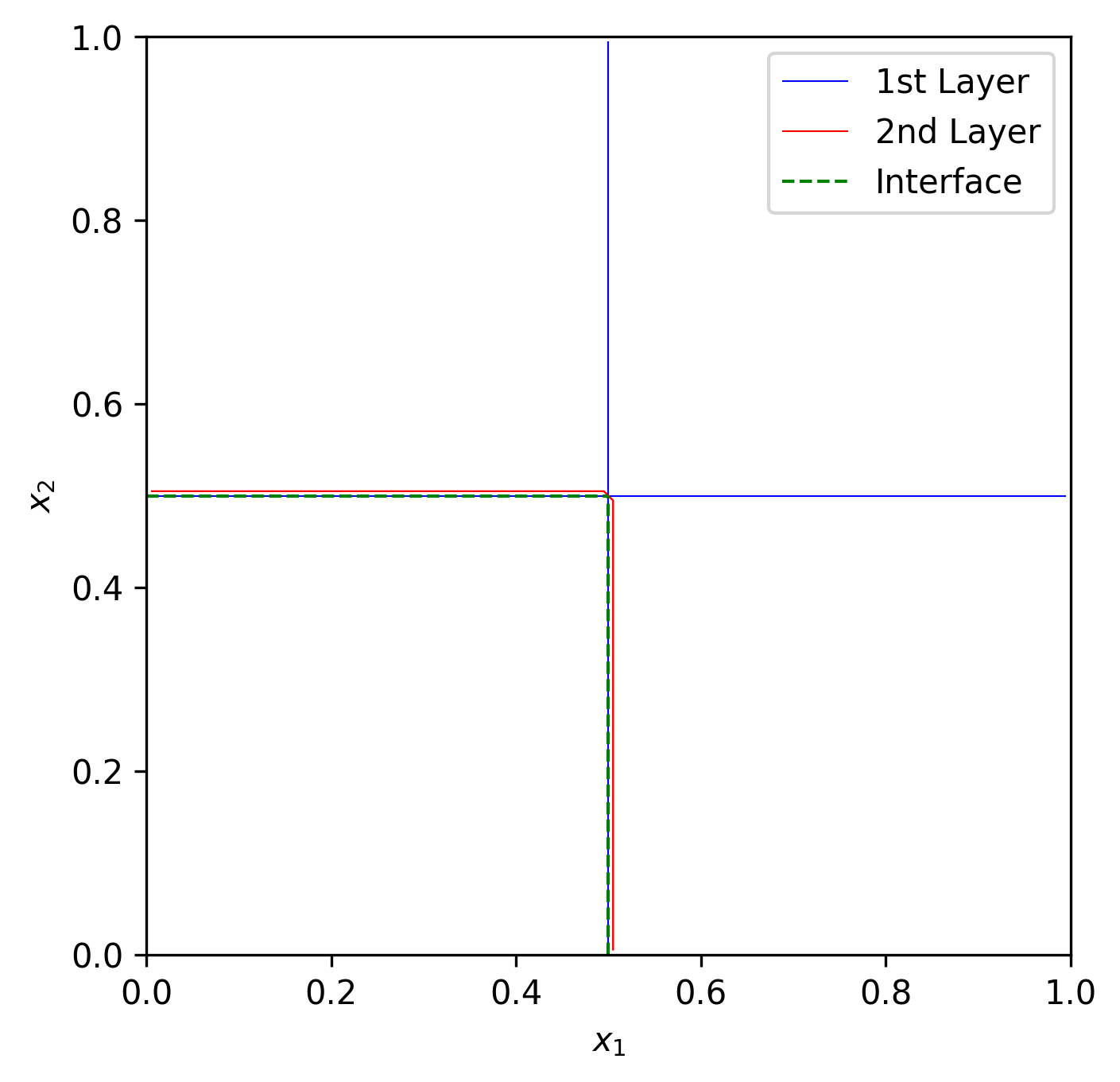}
\end{minipage}%
}%
\caption{A convex example to illustrate Theorem \ref{general theorem} for the case $d=10000$}
\label{nd example}
\end{figure}

Let $d=10000$, $\Omega=(0,1)^d$, 
\[
\Omega_1=\{\bx=(x_1,\ldots,x_d)\in \Omega: x_1< 1/2, \ldots, x_d < 1/2\}, \quad\mbox{and}\quad \Omega_2=\Omega\setminus\Omega_1. 
\]
The intersection between the piecewise constant function $\chi(\mathbf{x})$ and the hyperplanes $x_i=0.255$ for $i=3,\ldots, 10000$ is shown in Figure \ref{nd example function}. The interface $\Gamma$ is the boundary of a hypercube in $\Omega$ (see Figure \ref{nd interface} for a three-dimensional section of it):
\[
\Gamma=\bigcup_{i=1}^d\{\bx=(x_1,\ldots,x_d)\in \Omega: x_i = 1/2, \text{ and } 0\le x_j\le 1/2 \text{ for }j\neq i\}.
\]
In this example, we can simply take $\hat{\chi}=\chi$. Noting that the hypercube consists of 10000 hyperplanes of the form $x_i-1/2=0$ for $i=1,\ldots,10000$, the corresponding NN approximation is
\begin{equation}
\mathcal{N}(\mathbf{x})=1-\sigma\left(1-\frac{1}{\varepsilon}\sum_{i=1}^d\sigma(x_i-1/2)\right).
\end{equation}
Two sectional views of $\mathcal{N}(\bx)$ are shown in Figures \ref{nd approx graph e20} and \ref{nd approx graph e200} with $\varepsilon=1/20$ and $1/200$, respectively. Figures \ref{nd_breaking_e20} and \ref{nd_breaking_e200} plot the corresponding breaking hyperplanes.

\subsection{A two-dimensional non-convex example}

\begin{figure}[htbp]
\centering
\subfigure[The piecewise constant function $\chi(\mathbf{x})$\label{general example function}]{
\begin{minipage}[t]{0.4\linewidth}
\centering
\includegraphics[width=1.8in]{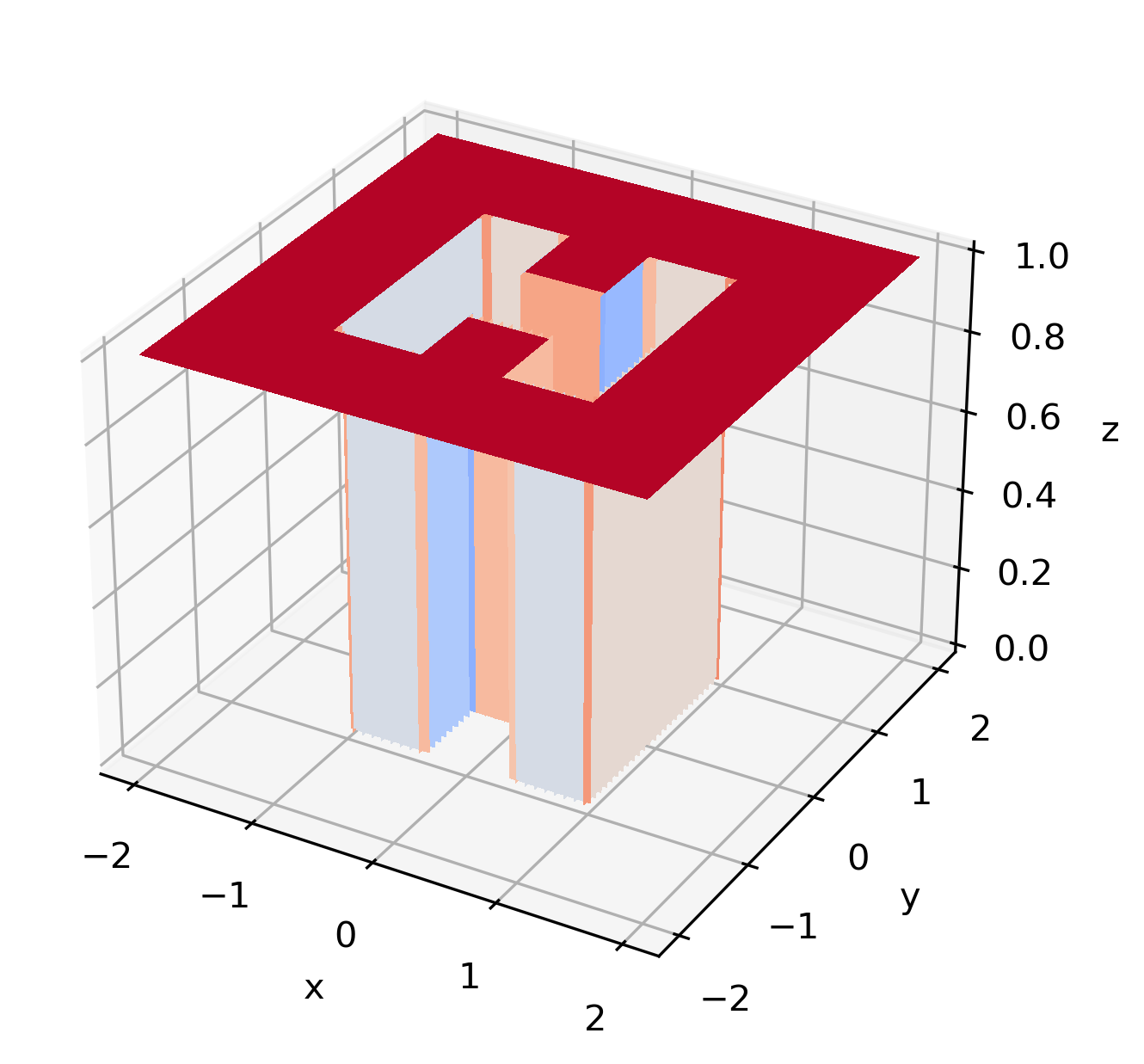}
\end{minipage}%
}%
\hspace{0.2in}
\subfigure[The interface\label{general interface}]{
\begin{minipage}[t]{0.4\linewidth}
\centering
\includegraphics[width=1.8in]{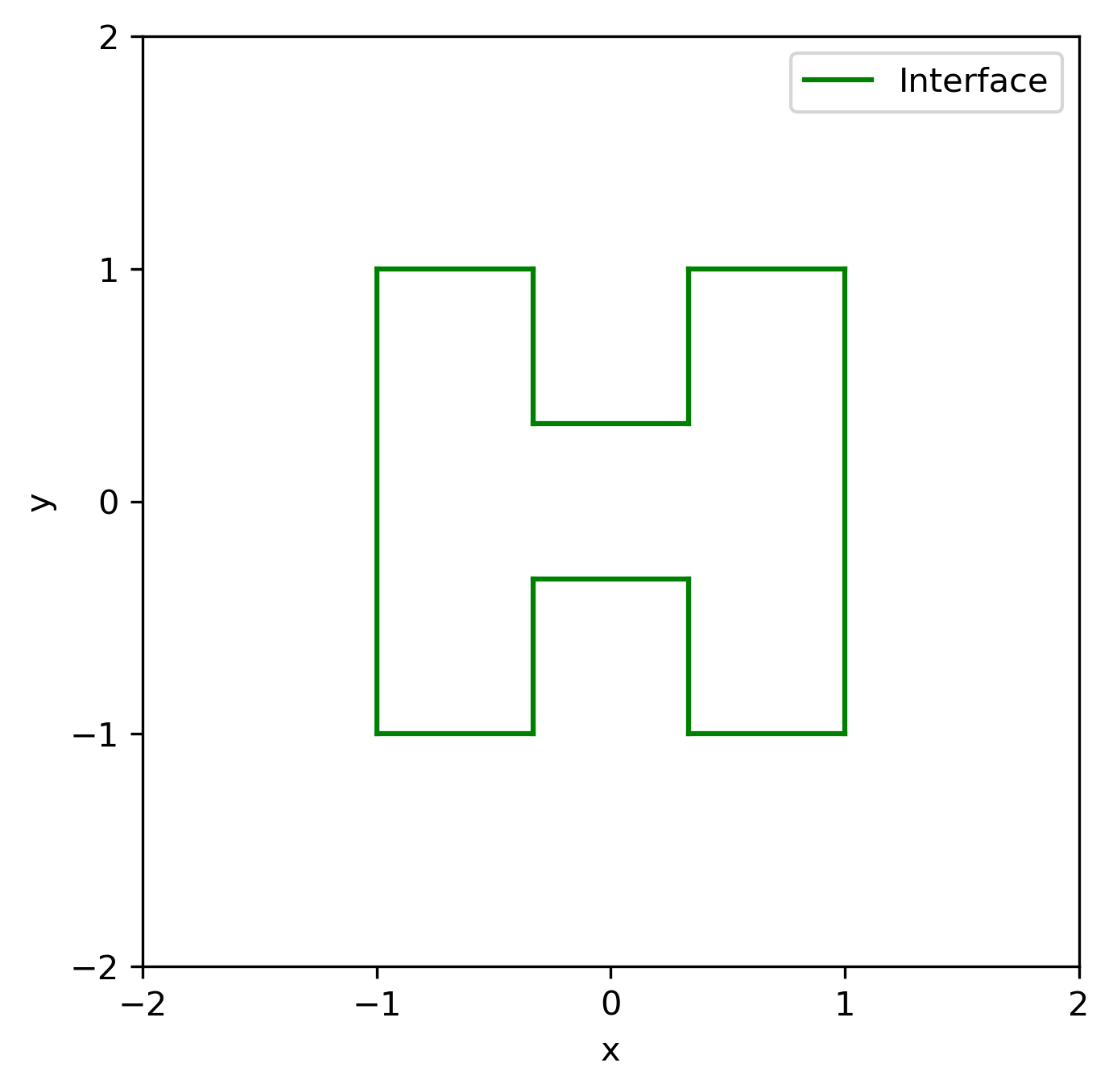}
\end{minipage}%
}%
\\
\subfigure[The convex hull of $\hat{\Omega}_1$\label{general_hull}]{
\begin{minipage}[t]{0.4\linewidth}
\centering
\includegraphics[width=1.8in]{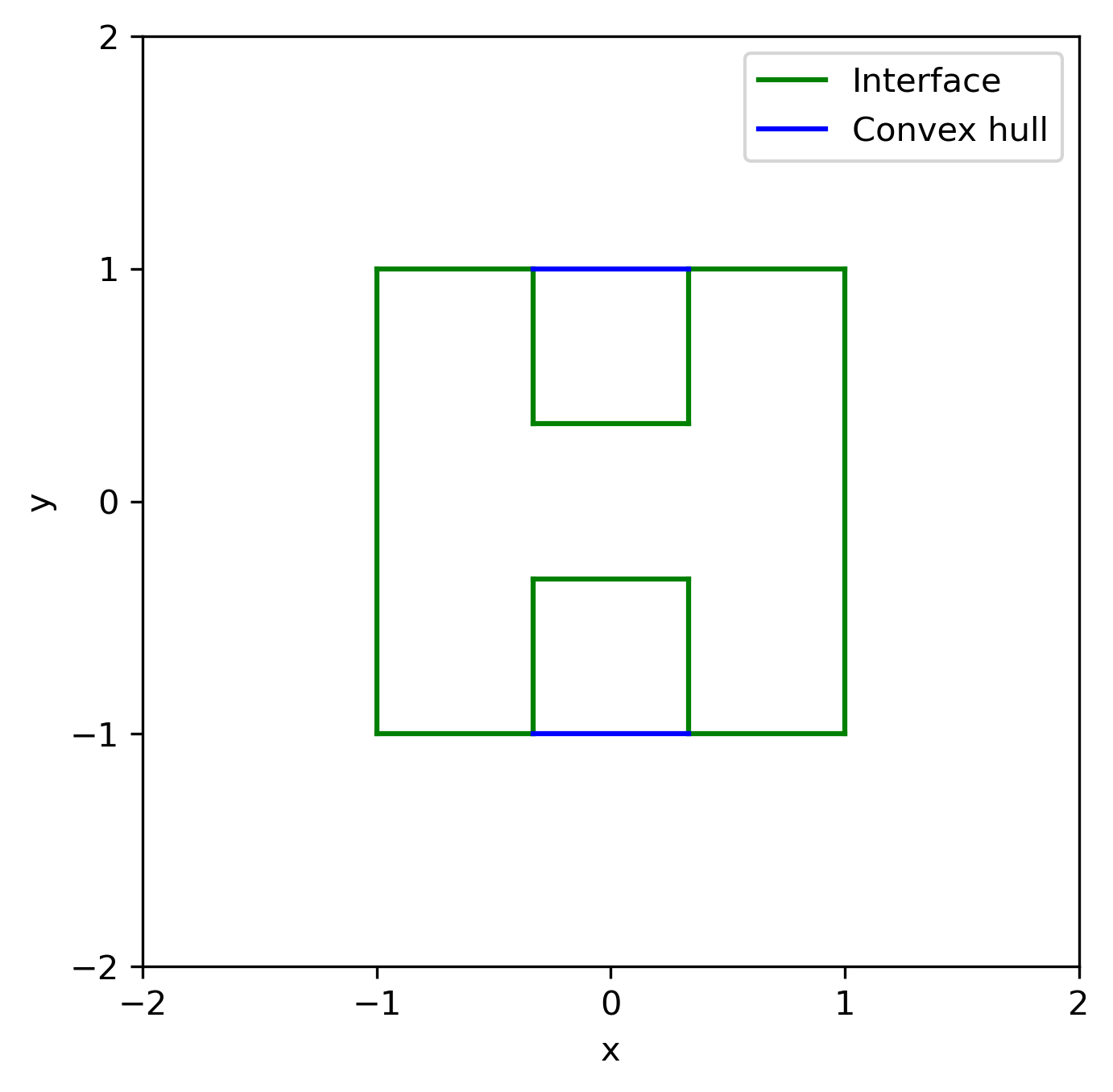}
\end{minipage}%
}%
\hspace{0.2in}
\subfigure[An approximation of $\chi(\mathbf{x})$ by the 2--12--3--1 ReLU NN function with $\varepsilon=1/12$\label{general_hull_graph_12}]{
\begin{minipage}[t]{0.4\linewidth}
\centering
\includegraphics[width=1.8in]{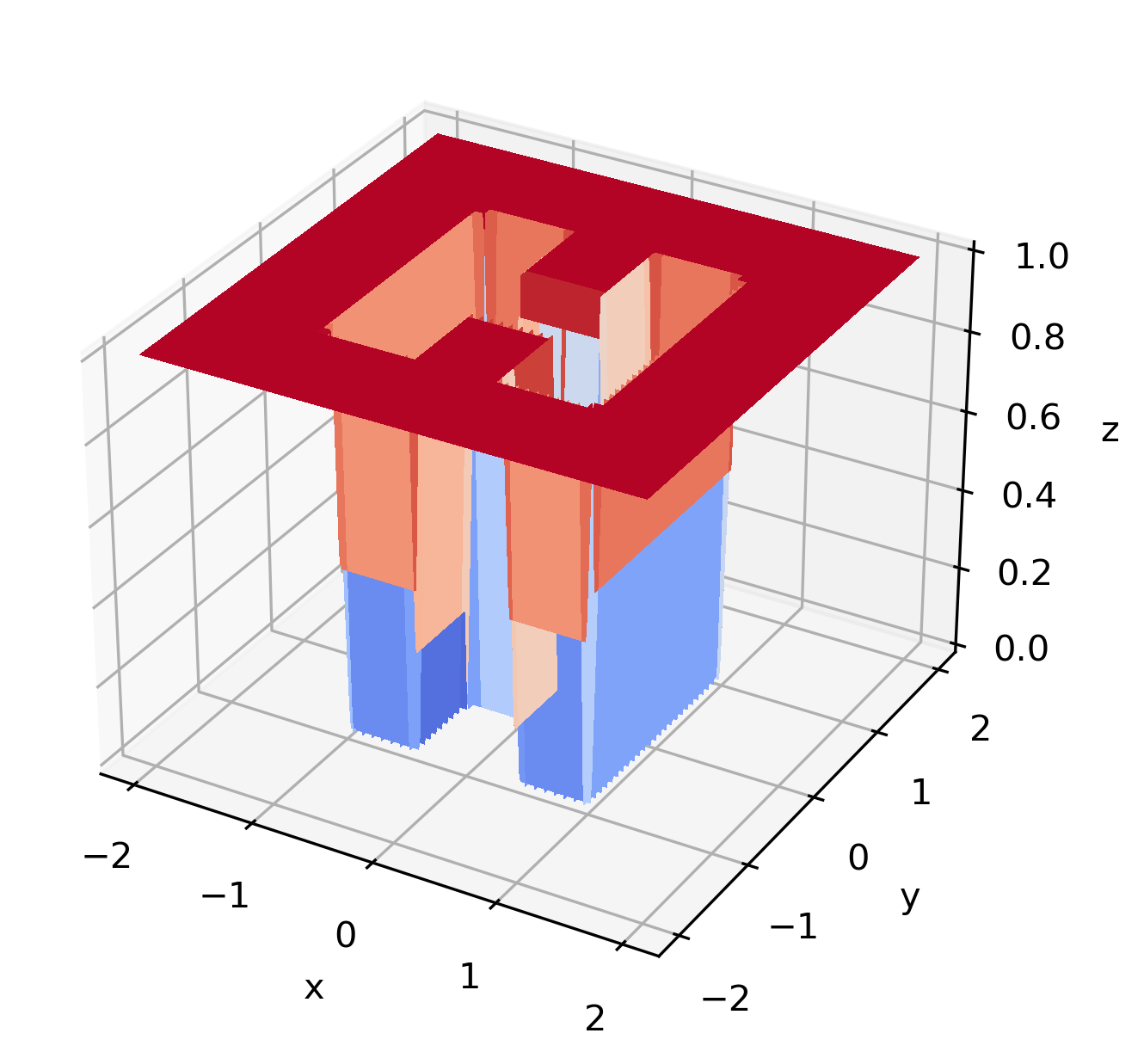}
\end{minipage}%
}%
\\
\subfigure[An approximation of $\chi(\mathbf{x})$ by the 2--12--3--1 ReLU NN function with $\varepsilon=1/200$\label{general_hull_graph_200}]{
\begin{minipage}[t]{0.4\linewidth}
\centering
\includegraphics[width=1.8in]{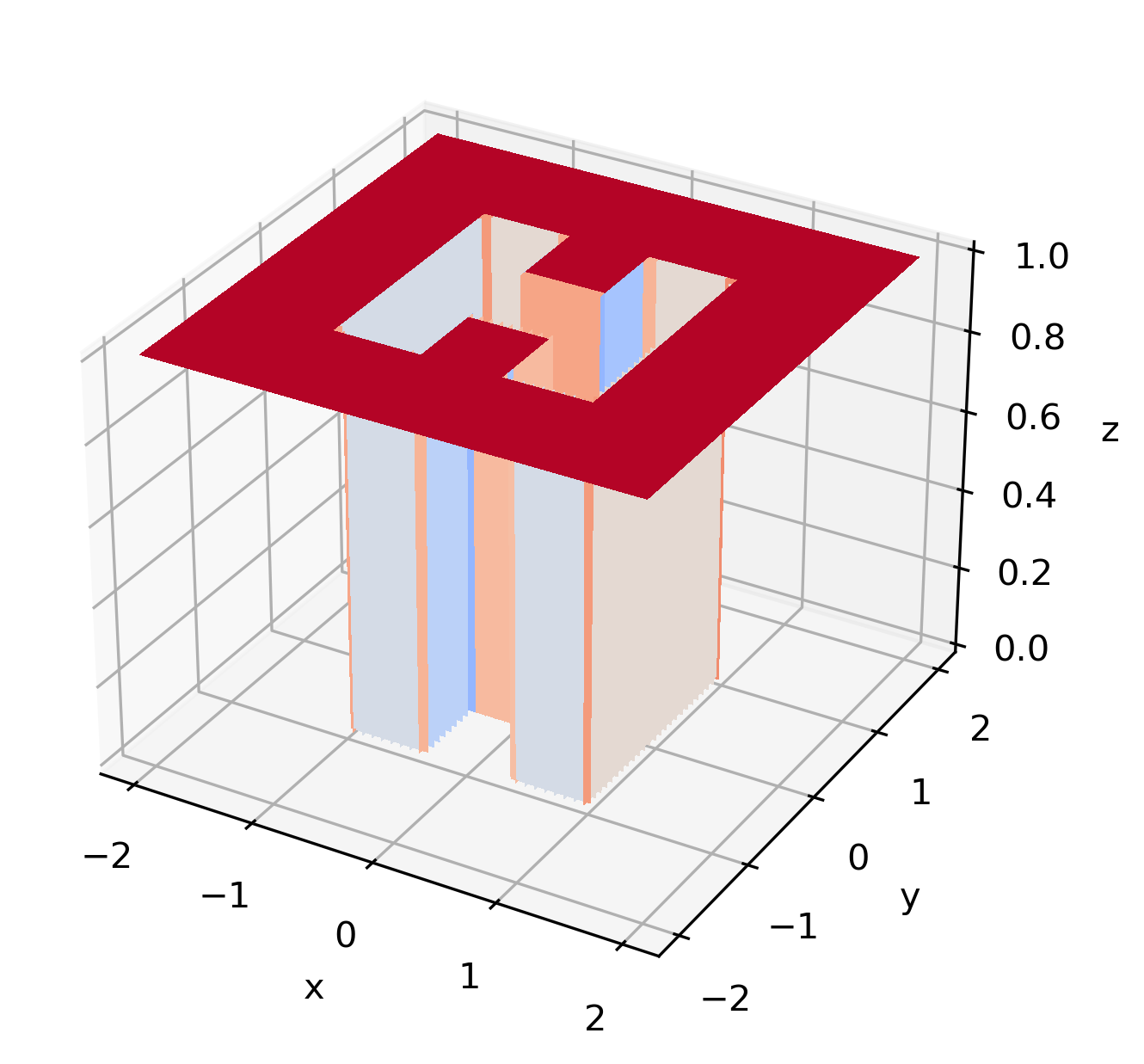}
\end{minipage}%
}%
\hspace{0.2in}
\subfigure[The breaking hyperplanes of the approximation in Figure \ref{general_hull_graph_12}\label{general_hull_breaking_12}]{
\begin{minipage}[t]{0.4\linewidth}
\centering
\includegraphics[width=1.8in]{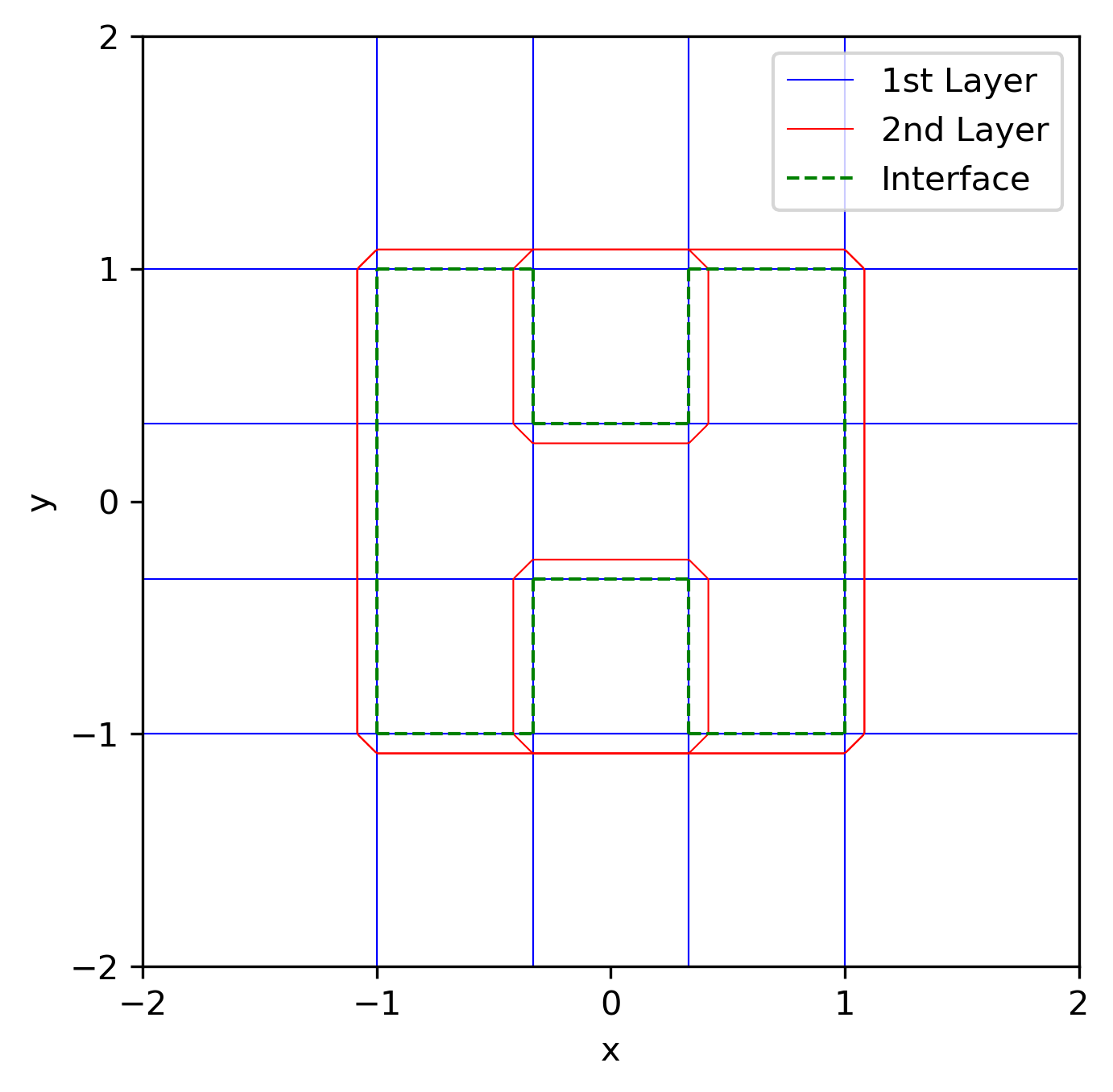}
\end{minipage}%
}%
\\
\subfigure[The breaking hyperplanes of the approximation in Figure \ref{general_hull_graph_200}\label{general_hull_breaking_200}]{
\begin{minipage}[t]{0.4\linewidth}
\centering
\includegraphics[width=1.8in]{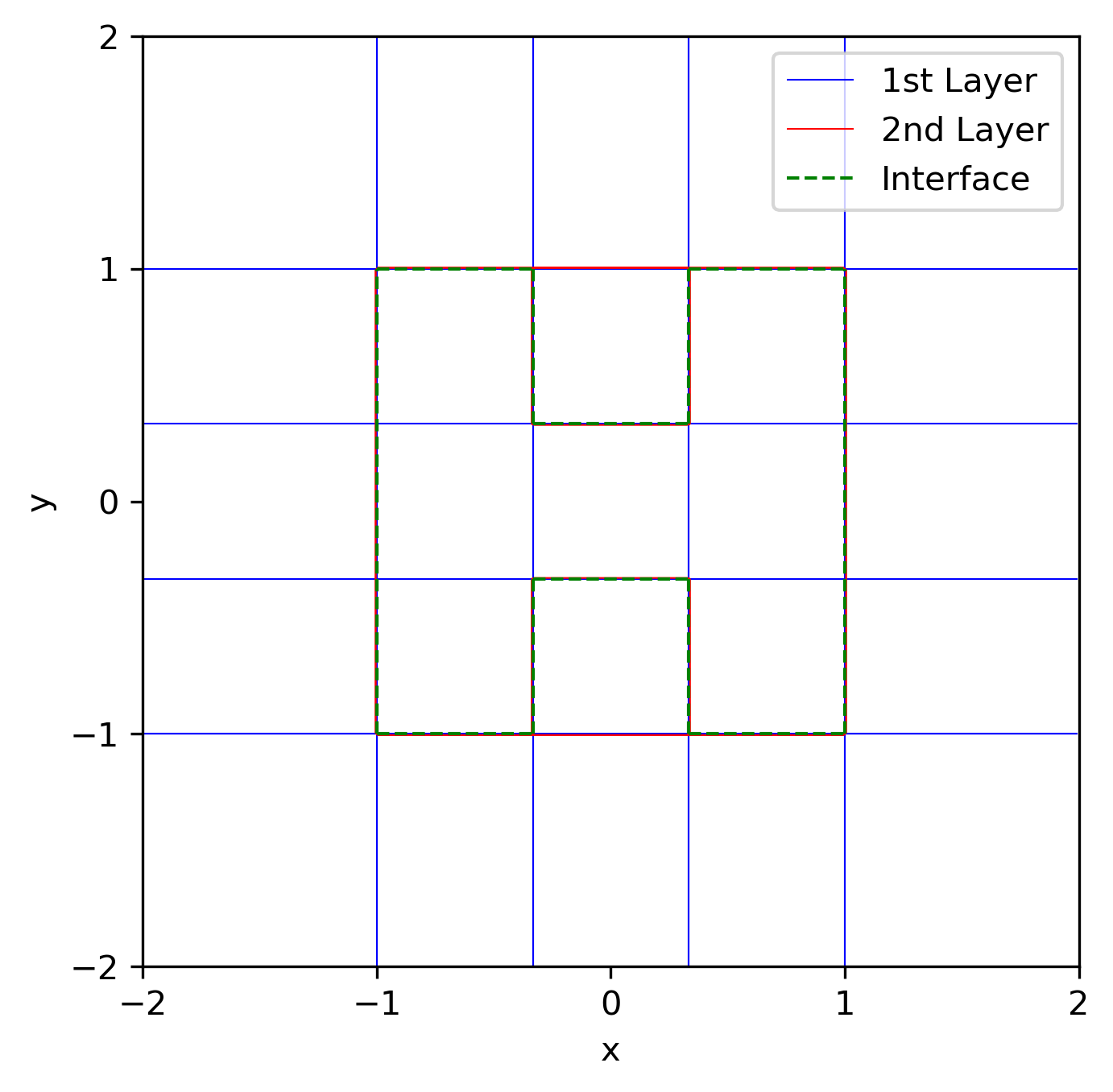}
\end{minipage}%
}%
\caption{A non-convex example to illustrate Theorem \ref{general theorem} for the case $d=2$ (convex hull)}
\end{figure}

\begin{figure}[htbp]
\centering
\subfigure[A convex decomposition of $\hat{\Omega}_1$\label{general_decomposition}]{
\begin{minipage}[t]{0.4\linewidth}
\centering
\includegraphics[width=1.8in]{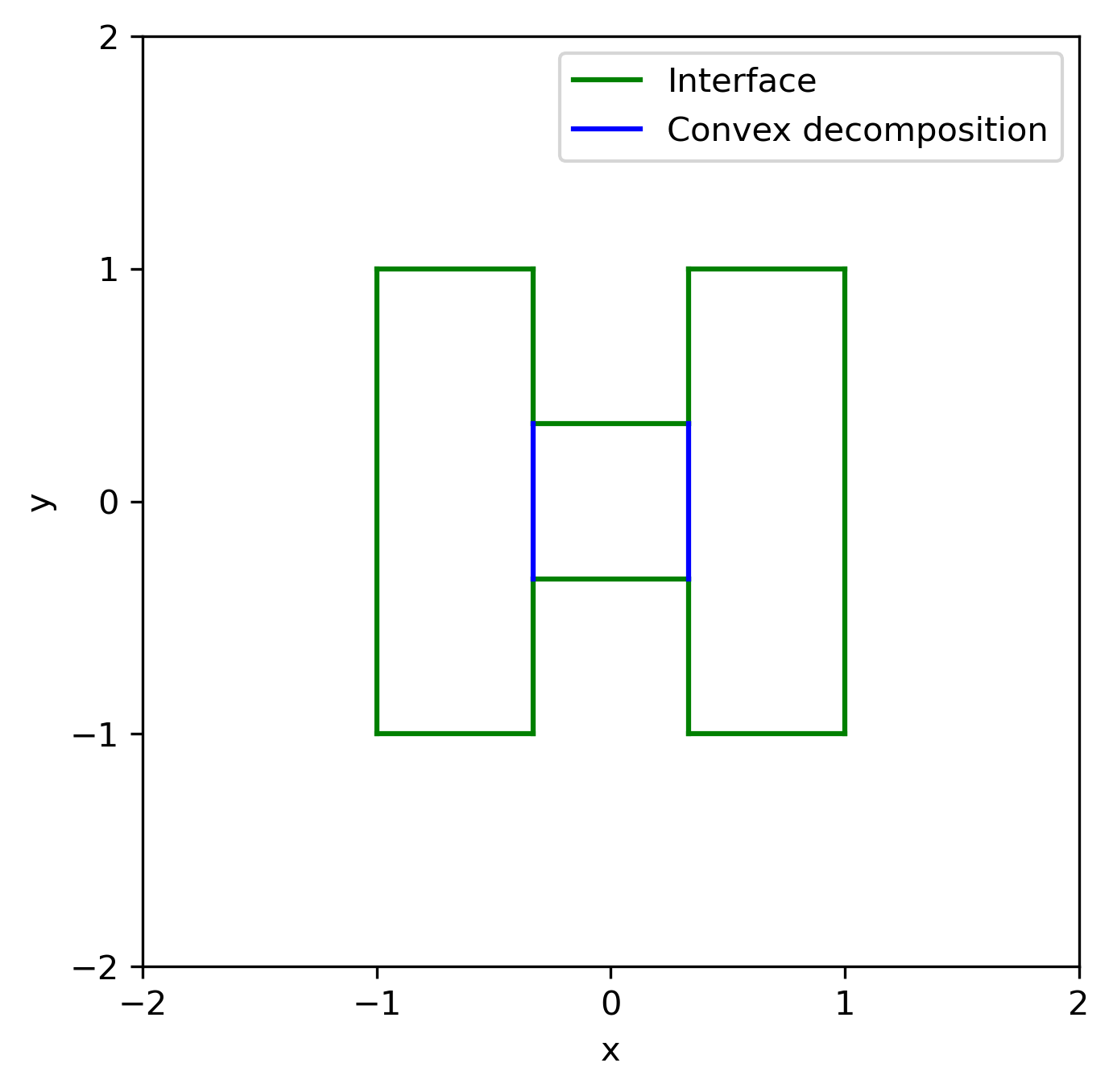}
\end{minipage}%
}%
\hspace{0.2in}
\subfigure[An approximation of $\chi(\mathbf{x})$ by the 2--12--3--1 ReLU NN function with $\varepsilon=1/12$\label{general_decomposition_graph_12}]{
\begin{minipage}[t]{0.4\linewidth}
\centering
\includegraphics[width=1.8in]{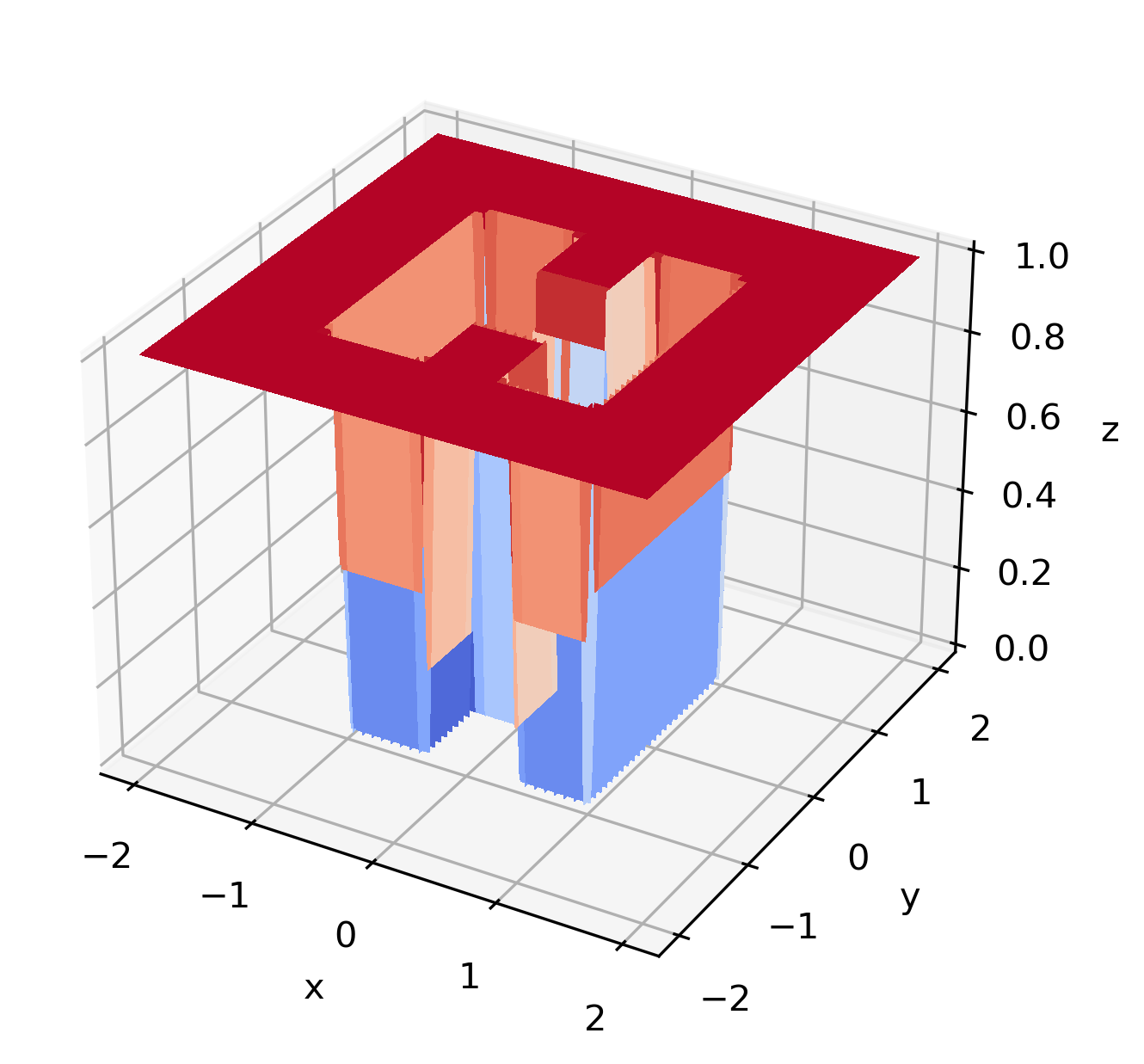}
\end{minipage}%
}%
\\
\subfigure[An approximation of $\chi(\mathbf{x})$ by the 2--12--3--1 ReLU NN function with $\varepsilon=1/200$\label{general_decomposition_graph_200}]{
\begin{minipage}[t]{0.4\linewidth}
\centering
\includegraphics[width=1.8in]{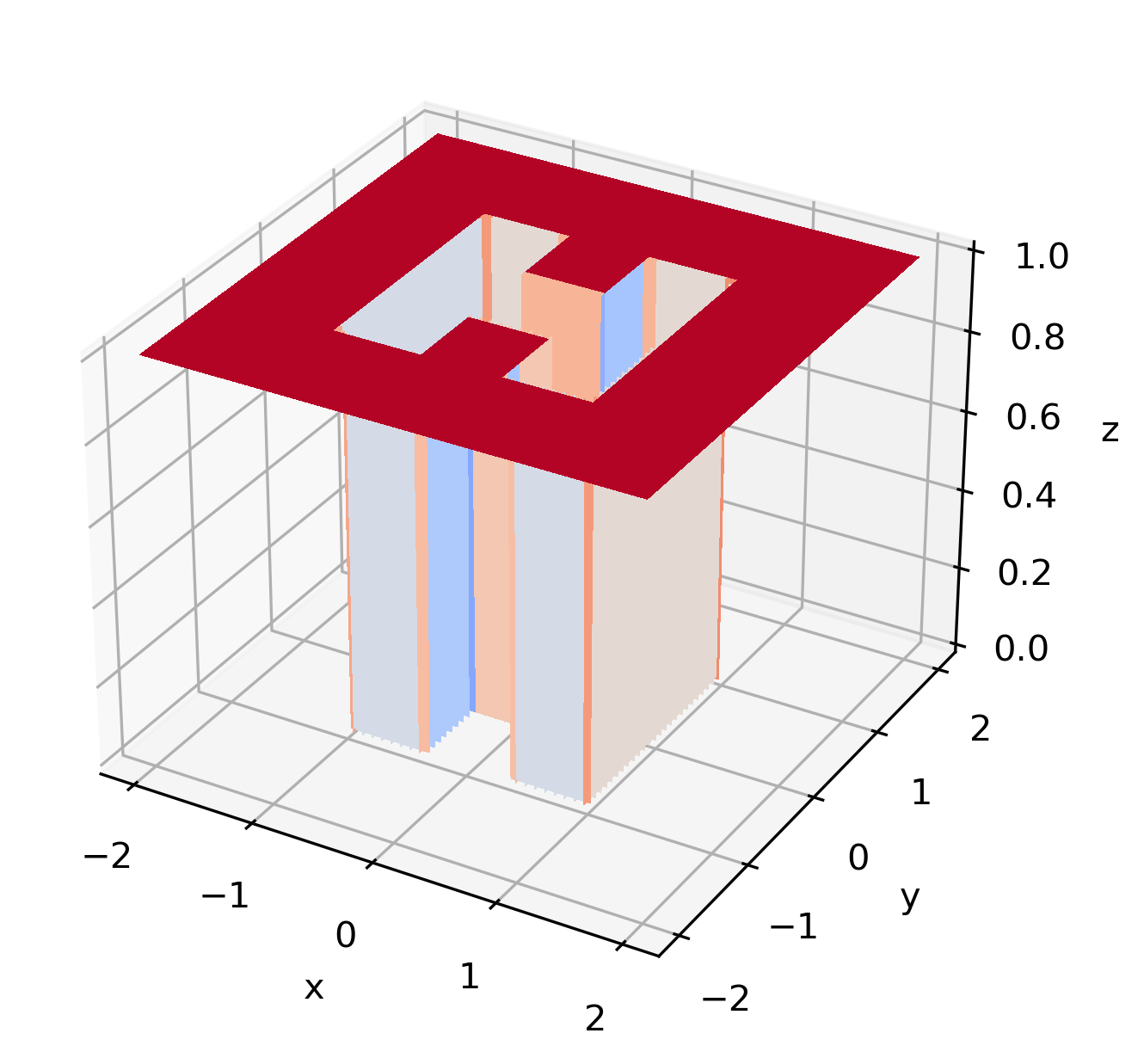}
\end{minipage}%
}%
\hspace{0.2in}
\subfigure[The breaking hyperplanes of the approximation in Figure \ref{general_decomposition_graph_12}\label{general_decomposition_breaking_12}]{
\begin{minipage}[t]{0.4\linewidth}
\centering
\includegraphics[width=1.8in]{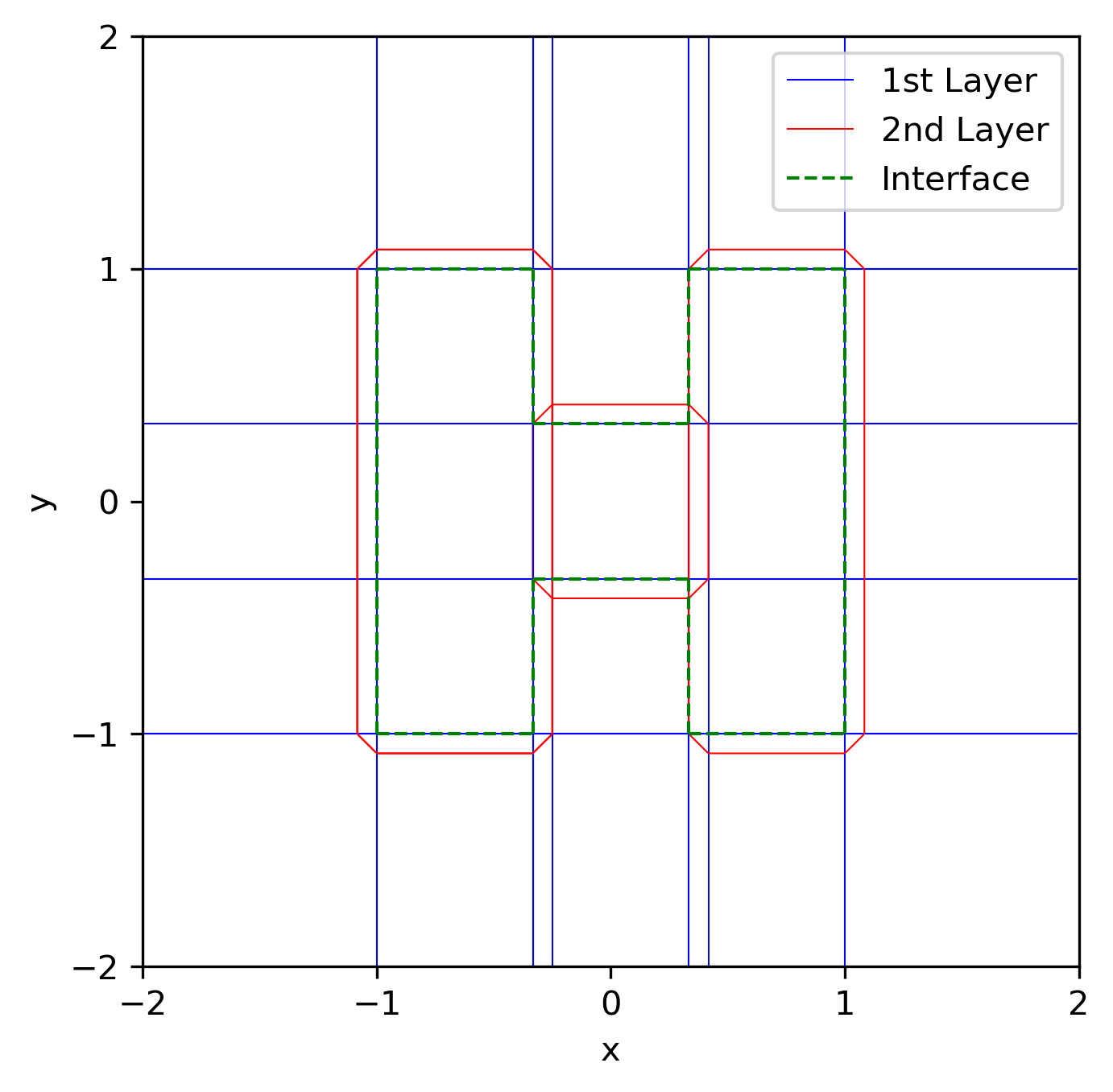}
\end{minipage}%
}%
\\
\subfigure[The breaking hyperplanes of the approximation in Figure \ref{general_decomposition_graph_200}\label{general_decomposition_breaking_200}]{
\begin{minipage}[t]{0.4\linewidth}
\centering
\includegraphics[width=1.8in]{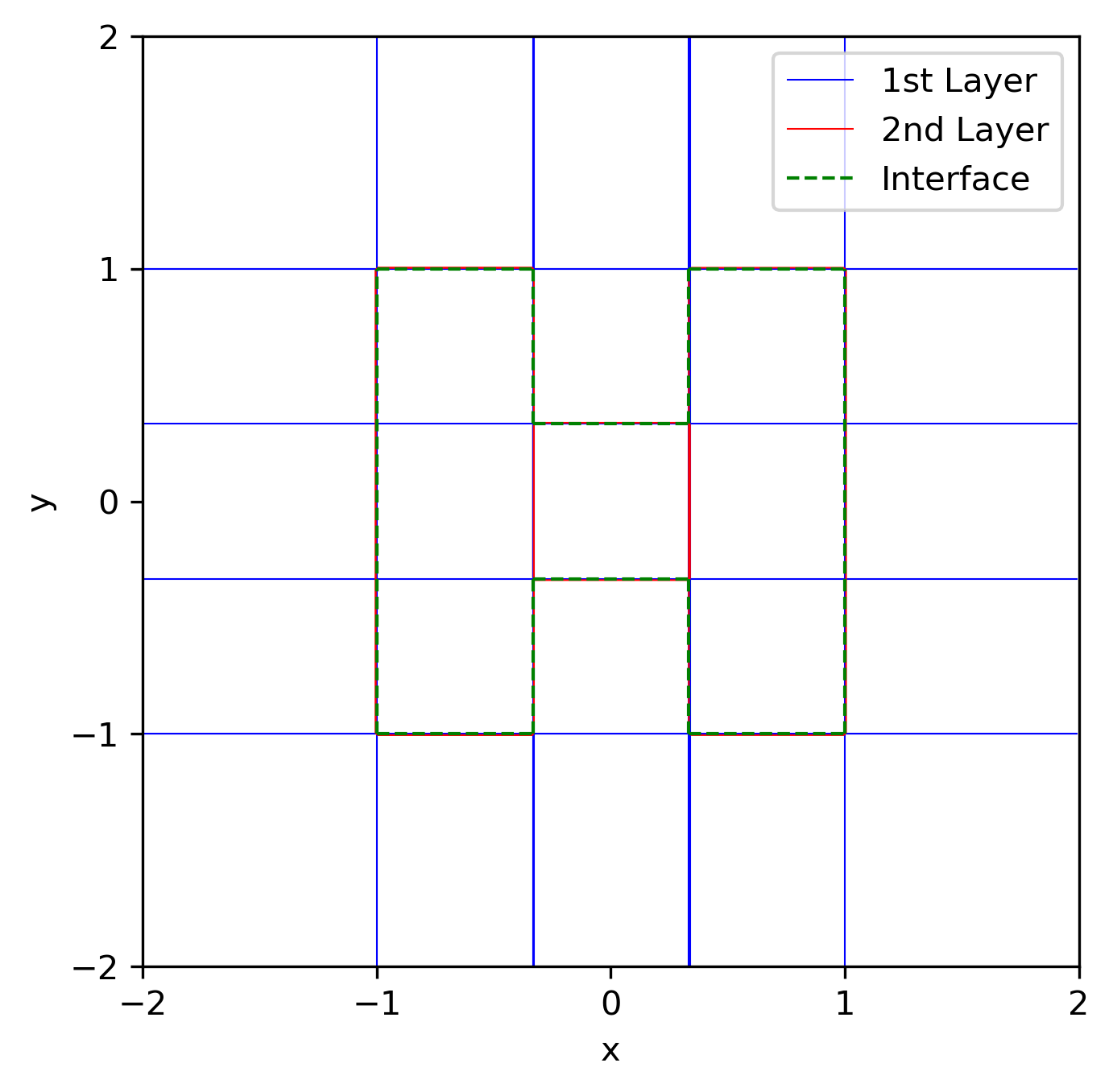}
\end{minipage}%
}%
\caption{A non-convex example to illustrate Theorem \ref{general theorem} for the case $d=2$ (convex decomposition)}
\end{figure}

Let $\Omega=(-2,2)^2$ and $\Omega_1$ be the H-shaped region depicted in Figure \ref{general interface} whose boundary is the interface $\Gamma=\partial\Omega_1\cap \partial\Omega_2 = \partial\Omega_1$. The unit step function $\chi(\mathbf{x})$ is depicted in Figure \ref{general example function}. Again, in this example, we can simply take $\hat{\chi}=\chi$. We construct 2--12--3--1 ReLU NN functions using 2--4--1--1 ReLU NN functions as discussed in Subsections \ref{convex hull subsubsection} and \ref{convex decomposition subsubsection} (see Figures \ref{general_hull} and \ref{general_decomposition}). The approximations with $\varepsilon=1/12$ and $1/200$ are depicted in Figures \ref{general_hull_graph_12} and \ref{general_hull_graph_200} for the convex hull and Figures \ref{general_decomposition_graph_12} and \ref{general_decomposition_graph_200} for the convex decomposition, respectively. Their corresponding breaking lines are plotted in Figures \ref{general_hull_breaking_12} and \ref{general_hull_breaking_200} for the convex hull and Figures \ref{general_decomposition_breaking_12} and \ref{general_decomposition_breaking_200} for the convex decomposition, respectively.

\bibliographystyle{siamplain}
\bibliography{references}
\end{document}